%

\magnification=\magstep1
\input amstex
\documentstyle{amsppt}
\pagewidth{6.5truein}
\pageheight{8.9truein}
\ifx\refstyle\undefinedZQA\else\refstyle{C}\fi
\loadbold
\loadeusm

\define\cf{\operatorname{cf}}
\define\open{\text{\rm open}}
\define\clopen{\text{\rm clopen}}
\define\Borel{\text{\rm Borel}}
\define\Umble{\text{\rm $U$\snug-measurable}}

\define\Fr{\text{\rm Fr}}

\define\rk{\text{\rm rk}}
\define\QED{ \null\nobreak\hfill$\blacksquare$}
\define\restrict{{\restriction}}
\define\nullset{\varnothing}
\define\nullseq{{\langle\rangle}}
\define\setdiff{\backslash}
\define\forces{\Vdash}
\define\DELTA{{\boldsymbol\Delta}}
\define\Suslin{{\Cal A}}
\define\treerestrict{\restriction^*}
\define\wkrestrict{{\restrict}}
\define\cfrestrict{{\restrict}}
\define\subalggen#1#2{{H_{#1}({#2})}}
\define\NNC{{NNC}}
\define\scriptX{{\eusm X}}
\define\scriptC{{\eusm C}}
\define\eps{\varepsilon}
\define\image#1{[#1]}
\define\concat{{}^\cap}
\define\Pname#1{{\dot{#1}}}
\define\An{{\Pname A}}
\define\Tn{{\Pname T}}
\define\Zn{{\Pname Z}}
\define\Yv{{\vec Y}}
\define\Zv{{\vec Z}}

\define\support#1{{\text{\rm Supp}({#1})}}
\define\domain{\operatorname{domain}}
\define\range{\operatorname{range}}
\define\SIGMA{{\boldsymbol\Sigma}}
\define\PI{{\boldsymbol\Pi}}
\define\R{{\bold R}}
\define\substitute(#1,#2,#3){{#1}({#2}/{#3})}
\define\suchthat{\,\colon\allowbreak\;}
\define\funcfrom{\,\colon\;}
\define\pair(#1,#2){({#1},{#2})}
\define\Cohentheta{{\text{\rm Fn}(\theta,2)}}
\define\Nbhd(#1){N_{#1}}
\define\setfuncs#1#2{{}^{#1}#2}
\define\setfuncsb#1#2{{}^{#1}\!#2}

\define\Bagemihl{1}
\define\Barnes{2}
\define\BaumgartnerHajnal{3}
\define\Devlin{4}
\define\DevlinParis{5}
\define\Dougherty{6}
\define\Gitik{7}
\define\Jech{8}
\define\Koepke{9}
\define\KuratowskiP{10}
\define\Louveau{11}
\define\Maharam{12}
\define\Moschovakis{13}
\define\MrowkaNMC{14}
\define\MrowkaSID{15}
\define\Prikry{16}
\define\Ristow{17}
\define\Shelah{18}
\define\Silver{19}
\define\Simms{20}
\define\Solovay{21}

\topmatter
\title Narrow coverings of $\omega$-ary product spaces\endtitle
\author Randall Dougherty\endauthor
\affil Ohio State University\endaffil
\date July 1, 1996 \enddate
\thanks
   This work was partially supported by grants from the
   National Science Foundation and the Alfred P. Sloan
   Foundation.
\endthanks
\address Department of Mathematics, Ohio State University,
Columbus, OH 43210\endaddress
\email rld\@math.ohio-state.edu \endemail
\abstract
Results of Sierpi\'nski and others have shown that certain
finite-dimensional product sets can be written as unions
of subsets, each of which is ``narrow'' in a corresponding
direction; that is, each line in that direction intersects
the subset in a small set.  For example, if the
set $\omega \times \omega$ is partitioned into two pieces
along the diagonal, then one piece meets every horizontal line
in a finite set, and the other piece meets each vertical line in a
finite set.  Such partitions or coverings can exist
only when the sets forming the product are of limited size.

This paper considers such coverings for products of infinitely
many sets (usually a product of $\omega$~copies of the same
cardinal~$\kappa$).  In this case, a covering of the product by
narrow sets, one for each coordinate direction, will exist no
matter how large the factor sets are.  But if one restricts the
sets used in the covering (for instance, requiring them to
be Borel in a product topology), then the existence of
narrow coverings is related to a number of large cardinal properties:
partition cardinals, the free subset problem, nonregular ultrafilters,
and so on.

One result given here is a relative consistency proof for a
hypothesis used by S.~Mr\'owka to construct a counterexample
in the dimension theory of metric spaces.
\endabstract
\endtopmatter

\document

%

\head 1. Introduction \endhead

     The set $\omega \times \omega $ can be partitioned along the
diagonal into two pieces $\{(m,n)\suchthat  m < n\}$ and $\{(m,n)\suchthat
m \ge  n\}$.  The first of these pieces has a property which might be
called ``narrowness in the first coordinate'':  for each $n$, there
are only finitely many $m$\snug's such that $(m,n)$ is in the set.
(In other words, each ``line in the direction of the first coordinate
axis'' has a relatively small intersection with the set.)  And the second
piece is ``narrow in the second coordinate.''  Similarly, if $\omega_1
\times \omega_1$ is divided into two pieces in this way, then each piece
contains only countably many points along each line in the corresponding
coordinate direction.  But $\omega_1 \times \omega_1$ turns out
to be too large to partition into two pieces which are
narrow in the finite sense.

By a more complicated construction, one can partition the set
$\omega_1\times \omega_1\times \omega_1$ into three pieces, each of
which is narrow in one of the three coordinates, in the sense of only
containing finitely many points on each line in the corresponding
coordinate direction.  If one allows the narrow sets to contain
countably many points on each such line, then a suitable partition
exists for the set $\omega_2 \times \omega_2 \times \omega_2$. The
$\omega_1$ and~$\omega_2$ are largest possible for the respective
partitions to exist. This is part of a large collection of results
proven by many people over the past eighty years.  A few more details
are given in section~2 of this paper; for a much more thorough
presentation of the subject, see Simms~\cite{\Simms}.

     The purpose of the present paper is to investigate the
problem of expressing an infinitary product as a union
of subsets, each of which is narrow in some coordinate direction.
More specifically, given a set~$X$ and a cardinal~$\lambda$, can the
$\omega $\snug-dimensional product $\setfuncsb{\omega}X$ be covered by
(written as a union
of) sets $A_n$ ($n < \omega$), where $A_n$ is $\lambda$\snug-narrow in
the $n$\snug'th coordinate direction (i.e., each line parallel to the
$n$\snug'th coordinate axis meets~$A_n$ in fewer than $\lambda$ points)?
Stated this way, the answer turns out to be `yes' no matter
how large $X$~is, for any $\lambda \ge 2$.  But if one puts further
restrictions on the sets~$A_n$ (e.g., that they be Borel in the
product topology on $\setfuncsb\omega X$ with $X$~discrete), then
one gets a number of interesting questions related to several
other well-known concepts --- partition cardinals, the free subset
problem, nonregular ultrafilters, and so on.

     The problem arose from a construction in dimension theory:
S.~Mr\'owka~\cite{\MrowkaSID,\MrowkaNMC} has shown that a hypothesis
called $(A_{\aleph_0})$ or $S(\aleph_0)$ implies the existence of a
metrizable space with zero inductive dimension whose completions (under
all possible metrics) have nonzero inductive dimension. (A topological
space has zero inductive dimension iff its topology has a basis of
clopen sets.)  The statement of $S(\aleph_0)$ is: if $X$~has
size~$2^{\aleph_0}$, then $\setfuncsb \omega X$ cannot be written as a
union of sets~$A_n$ ($n < \omega$) where $A_n$ is $\aleph_1$\snug-narrow
in the $n$\snug'th coordinate and is $F_\sigma$ in the product topology
on~$\setfuncsb \omega X$ with $X$~discrete. Here we will show that
$S(\aleph_0)$ is consistent relative to a large cardinal (the partition
cardinal $\kappa \to (\omega_1+\omega)^{<\omega}$), and that,
conversely, consistency of $S(\aleph_0)$ implies consistency of a
slightly smaller large cardinal ($\kappa \to (\omega)^{<\omega}$). So a
large cardinal well below the level of a measurable cardinal suffices
for the construction of Mr\'owka's example.

The organization of the present paper is as follows.  Section 2 gives
notational conventions, the main definitions of terms including those
used informally above, and some basic results. Section~3 gives
connections between narrow coverings, indiscernibles, and the free
subset problem, thus showing that large cardinals are necessary to get
the nonexistence of narrow coverings, and that slightly larger cardinals
are sufficient. Section~4 shows that some of these nonexistence results
are preserved under forcing which adds Cohen or random reals; this
suffices to prove the relative consistency of Mr\'owka's hypothesis
$S(\aleph_0)$.  Section~5 gives a method for using ultrafilters to prove
results about Borel sets (an approach previously taken by
Louveau~\cite{\Louveau}), and Section~6 uses this method to get results
about narrow coverings using suitably nonregular ultrafilters. Section~7
considers the question of how complicated a clopen narrow covering has
to be when it does exist; this leads to the study of ranks of trees of
finite free sequences. Section~8 lists some of the more interesting
questions which remain open. Sections~4 through~7 are independent of
each other, except that Section~6 depends on Section~5.

Much of this paper comes from my doctoral dissertation~\cite{\Dougherty};
however, other parts, such as the consistency proof for~$S(\aleph_0)$,
are new.

I would like to (and hereby do) thank Professors J.~Silver and
J.~Addison for many enlightening discussions, and
T.~Carlson and H.~Friedman for helpful comments.

\head 2.  Definitions and Basic Results \endhead

     Throughout this paper we will be working in ZFC, the
usual axioms of set theory including the axiom of choice. 
Cardinals will be initial ordinals; the cardinal $\aleph_{\alpha}$ will
be denoted by $\omega_{\alpha}$ when its set or ordinal nature is being
emphasized.  Since each cardinal is a set of its own
cardinality, we will not lose generality by stating many
results for cardinals rather than arbitrary sets. 
Natural numbers are finite ordinals, and each ordinal is
the set of its predecessors.  The immediate successor
(cardinal) of a cardinal $\lambda $ is denoted by $\lambda^+$.  The
cardinality of a set $S$ is denoted by $|S|$.

     For any function $f$ and any set $S$, $f\image S$ and $f^{-1}\image
     S$ denote
the image and inverse image of~$S$ under~$f$, respectively.  The
collection of all functions from $X$ to $Y$ is denoted by $\setfuncs XY$.
A sequence is a function~$s$ whose domain is an
ordinal; this ordinal is called~$\ell(s)$, the length of~$s$.  The
symbol~$\concat $ denotes
concatenation of sequences.  A sequence may be denoted by a list of its
members between angle brackets: $\langle \alpha,\beta,\gamma \rangle$,
$\langle \gamma \rangle$, $\nullseq$, $\langle a_n\suchthat  n < \omega
\rangle$, etc.

     If sets $S(\beta)$ are defined for all $\beta  < \alpha $, then
     $S({<}\alpha)$
will denote $\bigcup_{\beta<\alpha}S(\beta)$.  Variants such as $S({\le}
\alpha)$ are defined similarly.

A tree is a set $T$ of sequences such that any initial segment of a
member of\/ $T$ is a member of\/ $T$.  If\/ $T$ is a tree of finite
sequences, we define $[T]$ to be the set of sequences $s$ of length
$\omega $ such that $s\restrict n \in  T$ for all $n < \omega $.

\definition{Definition 2.1} Let $\scriptX$ be a product of sets.  A {\it
line parallel to the $n\!$\snug'th coordinate axis} in~$\scriptX$ is a
subset of\/~$\scriptX$ obtained by allowing the $n$\snug'th coordinate
of a point to vary while holding all other coordinates fixed. In other
words, the line parallel to the $n$\snug'th coordinate axis
in~$\scriptX$ through the point~$x$ is the set of $y \in \scriptX$ such
that $y(i) = x(i)$ for all $i \ne n$. \enddefinition

\definition{Definition 2.2}

(a) A subset $A$ of a product set~$\scriptX$ is
{\it $\lambda$\snug-narrow in the $n\!$\snug'th coordinate}
if every line parallel to the
$n$\snug'th coordinate axis in~$\scriptX$ meets~$A$ in fewer than
$\lambda$ points.

(b) A {\it $\lambda$\snug-narrow covering} of\/~$\scriptX$ is a
collection of sets~$A_n$, one for each coordinate~$n$, such that
$\bigcup_n A_n = \scriptX$ and, for each~$n$, $A_n$~is
$\lambda$\snug-narrow in the $n$\snug'th coordinate. \enddefinition

In particular, $\aleph_0$\snug-narrow means that each line in the
relevant direction contains only finitely many points of the set, while
$\aleph_1$-narrow means each such line contains countably many points of
the set. A $\lambda$\snug-narrow covering of\/~$\scriptX$ can easily be
converted into a {\sl partition} of\/~$\scriptX$ by replacing the
sets~$A_n$ with the sets $B_n = A_n \setdiff \bigcup_{m < n} A_m$, which
will still be $\lambda$\snug-narrow.

Clearly, for given $\lambda$ and $d$, the existence of
$\lambda$\snug-narrow coverings of the product $\setfuncsb dX$ depends
only on the cardinality of the set $X$.  Furthermore, if such a covering
exists for~$\setfuncsb dX$ (using sets $A_n \subseteq \setfuncsb dX$),
then one exists for $\setfuncs dY$ for any $Y \subseteq X$ (using the
sets $A_n \cap \setfuncs dY$). So, if such a covering does {\sl not}
exist for $\setfuncsb dX$, then one also does not exist for $\setfuncsb
d{X'}$ whenever $|X'|\ge|X|$.

The existence of narrow coverings for finite products of an infinite
set~$X$ has been studied by a number of authors; see
Simms~\cite{\Simms} for a full survey.  The main result along this
line is Theorem~2.149 of that survey, which comes from
Kuratowski~\cite{\KuratowskiP}.

\proclaim{Theorem 2.3 {\rm (Kuratowski)}} For any natural number $n >
0$, ordinal $\alpha$, and set $X$, there exists an
$\aleph_\alpha$\snug-narrow covering of\/~$\setfuncsb nX$ if and only
if\/ $|X| < \aleph_{\alpha+n-1}$. \QED\endproclaim

For the sake of completeness, we can consider the case of
$\lambda$-narrow coverings for {\sl finite}~$\lambda$ as well.

\proclaim{Proposition 2.4} For any natural numbers $n,m > 0$
and any set $X$, there exists an $m$\snug-narrow
covering of\/~$\setfuncsb nX$ if and only if\/ $|X| \le (m-1)n$.
\QED\endproclaim

\demo{Proof} Let $k = (m-1)n$.  It will suffice to show that
$\setfuncs nk$ has an $m$\snug-narrow covering, but
$\setfuncs n{(k+1)}$ does not.

Define sets $A_j \subseteq \setfuncs nk$ for $j < n$ as follows:
$$x \in A_j \iff (m-1)j \le \biggl(\sum_{i=0}^{n-1} x(i)\biggr) \bmod k
< (m-1)(j+1).$$
It is easy to check that the sets $A_j$ form an $m$\snug-narrow
covering of\/ $\setfuncs nk$.

On the other hand, a subset of\/ $\setfuncs n{(k+1)}$ which is
$m$\snug-narrow in any coordinate must contain at most
$(k+1)^{n-1}(m-1)$ points, so the union of $n$~such sets contains at
most $(k+1)^{n-1}k$ points, and hence is not all of\/ $\setfuncs
n{(k+1)}$.  Therefore, $\setfuncs n{(k+1)}$ has no $m$\snug-narrow
covering. \QED\enddemo

We now move on to products of infinitely many sets, specifically
products of the form~$\setfuncsb \omega X$.  The preceding results
would suggest that a $\lambda$\snug-narrow covering
of\/~$\setfuncsb \omega X$ exists if $X$~is sufficiently small, but not
if $X$~is too large.  The following result shows that
the actual situation is rather different.  This result was proved
for $X = \R$ by Bagemihl~\cite{\Bagemihl} using methods
of Davies; see Theorem~3.60 of Simms~\cite{\Simms}.

\proclaim{Theorem 2.5} For any $X$, there is a\/ $2$\snug-narrow
covering of\/~$\setfuncsb \omega X$. \endproclaim

\demo{Proof} Define an equivalence relation $\sim $ on $\setfuncsb
{\omega}X$ by:  $x  \sim  y $ iff $\{i\suchthat  x (i) \ne  y (i)\}$ is
finite.  For each $x$, let $[x]$ be the equivalence class of $x$.
Choose a representative $r(c) \in  c$ for each equivalence class $c$.
Let
     $$A_n = \{x  \in  \setfuncsb {\omega}X\suchthat  x(n) = r([x])(n)\}.$$
If $x\in  \setfuncsb {\omega}X$, $y \in  A_n$, and $y$~is on the line
parallel to the $n$\snug'th coordinate axis through~$x$, then $y \sim x$, so
     $$y(n) = r([y])(n) = r([x])(n),$$
so $y(n)$ is uniquely determined; hence, $A_n$ is $2$\snug-narrow
in the $n$\snug'th coordinate.
Any $x \in  \setfuncsb {\omega}X$ is in $A_n$ for all but finitely
many $n$, since $x \sim  r([x])$, so $\setfuncsb {\omega}X =
\bigcup_{n<\omega}A_n$ and we are done.\QED\enddemo

This sort of proof is commonly referred to as a ``blatant application
of the Axiom of Choice.''  (The proof also involves a blatant
application of the Axioms of Separation, but people tend to be
less concerned about that.)  The usual reaction to such a construction
is ``But is there an example using `reasonable' sets?''  This leads
to the following definition, which is stated negatively because
we will usually be considering circumstances under which
narrow coverings do not exist.

\definition{Definition 2.6} Given a set $X$, a cardinal $\lambda$,
and a property (or collection)~$P$ of subsets of\/~$\setfuncsb \omega X$,
we say that $\NNC(X,\lambda,P)$ holds iff there does not
exist a $\lambda$\snug-narrow covering of\/~$\setfuncsb \omega X$
using sets satisfying (or in)~$P$.
\enddefinition

The property~$P$ will often be `open' or `Borel' or some
other property from topology; in these cases, we will assume that
the topology on~$\setfuncsb \omega X$ is the product topology
with $X$~discrete.

As noted before, the existence of narrow coverings of\/~$\setfuncsb
\omega X$ depends only on the cardinality of~$X$; hence, we will usually
just consider the case where $X$~is itself a cardinal. A narrow covering
of\/~$\setfuncsb \omega X$ can be cut down to give a narrow covering
of\/~$\setfuncs \omega Y$ for any $Y \subseteq X$. Also, if the
condition~$P$ and the narrowness requirement on the sets in the covering
are relaxed, then any narrow coverings that worked for the strict
conditions will still work for the relaxed conditions.  These two
trivial monotonicity properties can be stated together as follows.

\proclaim{Lemma 2.7} If $\neg \NNC(X_0,\lambda_0,P_0)$, $X_1 \subseteq
X_0$, $\lambda_1 \ge  \lambda_0$, and $\{A\cap \setfuncsb
{\omega}{X_1}\suchthat  A \in  P_0\} \subseteq  P_1$, then $\neg
\NNC(\kappa_1,\lambda_1,P_1)$. \QED\endproclaim

If we have a $\lambda$\snug-narrow covering $\langle A_i \suchthat i < n
\rangle$ of a finitary product~$\setfuncsb n X$, then we can convert it
into a $\lambda$\snug-narrow covering $\langle B_i \suchthat i < \omega
\rangle$ of\/~$\setfuncsb \omega X$ by letting $B_i = \nullset$ for $i
\ge n$ and $B_i = \{x\suchthat x\restrict n \in A_i\}$ for $i < n$.
Since membership of a point~$x$ in the sets~$B_i$ depends only on the
first $n$ coordinates of~$X$, these sets are clopen in~$\setfuncsb
\omega X$.  Therefore, Theorem~2.3 (with $n = m+2$) gives the following
consequence.

\proclaim{Corollary 2.8} For any ordinal $\alpha$ and any
$m < \omega$, $\neg\NNC(\aleph_{\alpha+m},\aleph_\alpha,\clopen)$.
\QED\endproclaim

     There is no way to extend this result to get
$\neg\NNC(\aleph_{\alpha+\omega},\aleph_{\alpha},\clopen)$, as we will see
in the next section.

Sometimes the following slight variant of $\NNC(X,\lambda,P)$ is useful.

\definition{Definition 2.9} Given a set $X$, a cardinal $\lambda$, and a
property (or collection)~$P$ of subsets of\/~$\setfuncsb \omega X$, we
say that $\NNC(X,{<}\lambda,P)$ holds iff there do not exist sets $A_n
\subseteq \setfuncsb \omega X$ with property (or in collection)~$P$ such
that $\bigcup_{n<\omega}A_n = \setfuncsb \omega X$ and, for each $n$,
$A_n$ is $\lambda'_n$\snug-narrow in the $n$\snug'th coordinate for some
$\lambda'_n < \lambda$. \enddefinition

So $\NNC(X,\lambda,P)$ implies $\NNC(X,{<}\lambda,P)$, which in turn
implies $\NNC(X,\lambda',P)$ for all $\lambda' < \lambda$.  In fact, if
$\cf \lambda > \omega$, then $\NNC(X,{<}\lambda,P)$ is equivalent to
$(\forall \lambda'{<}\lambda)\,\NNC(X,\lambda',P)$ (because the supremum
of the cardinals~$\lambda'_n$ from the definition will be a cardinal
$\lambda' < \lambda$).  But if $\cf \lambda = \omega$, then
$\NNC(X,{<}\lambda,P)$ says a little more.

%

\head 3.  Indiscernibles and the Free Subset Problem \endhead

     In this section, we will show that the statement
$\NNC(\kappa,\mu^+,\open)$ is equivalent to a more familiar
assertion, namely that every structure on $\kappa$ with $\mu
$ operations has an infinite free subset.  In particular, this
will show that $\NNC(\kappa,\aleph_1,\open)$ implies the large
cardinal property $L\models  \kappa\to (\omega)^{<\omega}$.  On the
other hand, a similar but stronger property will be shown to imply
$\NNC(\kappa,\lambda,\Borel)$.  We will start with the latter result,
the idea for which was suggested to me by J.~Silver.

     Recall some definitions from partition theory.
For any set $S$ and any natural number~$n$, let $[S]^n = \{a \subseteq
S\suchthat  |a| = n\}$; let $[S]^{<\omega} = \bigcup_{n < \omega} [S]^n$.
If $\kappa$ and $\lambda $ are cardinals and $\alpha
$ is a limit ordinal, then $\kappa\to (\alpha)^{<\omega}_\lambda$ denotes
the assertion that, for any $F\funcfrom  [\kappa]^{<\omega}\to \lambda
$, there is a set $S \subseteq  \kappa$ of order type $\alpha $ such
that, for each $n < \omega $, $F$ is constant on $[S]^n$.  (We will omit
the $\lambda $ in the case $\lambda  = 2$.)  Jech \cite{\Jech,
pp.~392-396} gives a number of facts about this property, among which
is the result of Rowbottom
that $\kappa\to (\alpha)^{<\omega}$ implies $\kappa\to
(\alpha)^{<\omega}_{2^{\aleph_0}}$.

\proclaim{Theorem 3.1} Let $\kappa$,~$\lambda $, and~$\mu $ be cardinals,
and let $S$~be the collection of subsets of\/~$\setfuncs {\omega}\kappa$
which can be expressed as Boolean combinations of $\mu $~open subsets
of\/~$\setfuncs {\omega}\kappa$.
If $\kappa\to (\lambda +\omega)^{<\omega}_{2^{\mu}}$ (here $+$~is ordinal
addition), then $\NNC(\kappa,\lambda,S)$.
If $\lambda$~is infinite and $\kappa\to (\lambda)^{<\omega}_{2^{\mu}}$,
then $\NNC(\kappa,{<}\lambda,S)$.
\endproclaim

\demo{Proof} The case $\mu  = 0$ is trivial, so, by the preceding
remark, we may assume that $\mu $ is infinite.  Let $\langle
A_n\suchthat  n < \omega \rangle$ be any sequence of sets in $S$ such
that $\bigcup_{n<\omega}A_n = \setfuncs {\omega}\kappa$.  For the first
implication, assume $\kappa\to (\lambda +\omega)^{<\omega}_{2^{\mu}}$;
we must find an $n$ such that $A_n$ is not $\lambda$\snug-narrow in the
$n$\snug'th coordinate.

Each $A_n$ is a Boolean combination of $\mu $ open sets, so there is a
sequence $\langle G_{\alpha}\suchthat  \alpha  < \mu \rangle$ of open
subsets of\/ $\setfuncs {\omega}\kappa$ such that each $A_n$ is a Boolean
combination of these open sets.  Define a function $F\funcfrom
[\kappa]^{<\omega}\to \setfuncs {\mu}2$ as follows:  for any strictly
increasing sequence $\sigma  \in  \setfuncs {<\omega}\kappa$ and any
$\alpha  < \mu $, let $F(\range(\sigma))(\alpha) = 1$ iff $\{\sigma
\concat s\suchthat s \in  \setfuncs {\omega}\kappa\} \subseteq
G_{\alpha}$.  Since $\kappa\to (\lambda +\omega)^{<\omega}_{2^{\mu}}$,
there is a strictly increasing function $g\funcfrom  \lambda +\omega \to
\kappa$ such that $F$ is constant on $[\range(g)]^n$ for each $n <
\omega $.

Now, suppose $s$ and~$s'$ are strictly increasing sequences of elements
of $\range(g)$ of length~$\omega $, and $\alpha  < \mu $.  If $s \in
G_{\alpha}$, then there is $n < \omega $ such that $\{(s\restrict
n)\concat t\suchthat  t \in  \setfuncs {\omega}\kappa\} \subseteq
G_{\alpha}$, since $G_{\alpha}$ is open.  This gives $F(s\image
n)(\alpha) = 1$, so $F(s'\image n)(\alpha) = 1$, so $\{(s'\restrict
n)\concat t\suchthat  t \in \setfuncs {\omega}\kappa\} \subseteq
G_{\alpha}$, so $s' \in  G_{\alpha}$. Conversely, if $s' \in
G_{\alpha}$, then $s \in  G_{\alpha}$ by the same argument.  Therefore,
$s \in  G_{\alpha}$ iff $s' \in  G_{\alpha}$ for each $\alpha  < \mu $,
so, since $A_n$ is a Boolean combination of the sets $G_{\alpha}$, $s
\in  A_n$ iff $s' \in  A_n$ for each $n < \omega $.

There is at least one $n$ such that $s \in  A_n$, so, since $s$ and $s'$
are arbitrary, there is an $n < \omega $ such that, for all strictly
increasing $s \in  \setfuncs {\omega}{(\range(g))}$, $s \in  A_n$.
In particular, if we let
     $$s_\beta = (g\restrict n)\concat \langle g(n+\beta)\rangle\concat
     \langle g(n+\lambda +m)\suchthat  m < \omega \rangle$$
for $\beta  < \lambda $, we will have $s_\beta \in  A_n$ for all $\beta
< \lambda $; since $s_\beta(m) \ne  s_\gamma(m)$ only if $m = n$, $A_n$
is not $\lambda$\snug-narrow in the $n$\snug'th coordinate.

This completes the proof of the first implication.  The proof of
the second is similar:
Define~$F$ as before, and let $g\funcfrom \lambda \to \kappa$ be
increasing with $F$~constant on $[\range(g)]^n$ for each~$n$.
Find~$n$ such that all increasing $\omega$\snug-sequences from
$[\range(g)]^n$ are in $A_n$.
For any $\lambda' < \lambda$, we can find in~$\range(g)$
an increasing sequence of $n$~elements followed by $\lambda'$~elements
followed by $\omega$~elements; use these elements to form
sequences~$s_\beta$ for $\beta < \lambda'$ in~$A_n$ which differ
only at the $n$\snug'th coordinate.  This shows that
$A_n$~is not $\lambda'$\snug-narrow in the $n$\snug'th coordinate.
\QED\enddemo

\proclaim{Corollary 3.2} If $\kappa\to (\lambda +\omega)^{<\omega}$,
then $\NNC(\kappa,\lambda,\Borel)$.  If $\lambda$ is infinite
and $\kappa\to (\lambda)^{<\omega}$,
then $\NNC(\kappa,{<}\lambda,\Borel)$. \QED\endproclaim

     Now we give the relation between $\NNC(\kappa,\lambda,\open)$ and the
free subset problem, which has been considered in papers by
Devlin~\cite{\Devlin, \S4}, Devlin and~Paris~\cite{\DevlinParis},
Shelah~\cite{\Shelah}, and Koepke~\cite{\Koepke}, among others.
The relevant definitions are as follows.  If $S$ is a subset of (the
domain of) a structure $M$, let $\subalggen MS$ be the substructure of~$M$
generated by $S$. Such a set $S$ is said to be {\it free for $M$} iff,
for every $S' \subseteq  S$, $(\subalggen M{S'})\cap S = S'$. If $\kappa$,
$\lambda $, and $\mu $ are cardinals, then $\Fr_{\mu}(\kappa,\lambda)$
means that every structure of cardinality $\kappa$ with $\mu $ operations
(possibly including $0$\snug-ary operations, i.e., constants) has a free
subset of cardinality $\lambda $.

\proclaim{Theorem 3.3} For any infinite cardinals $\kappa$ and $\mu $,
$\Fr_{\mu}(\kappa,\aleph_0)$ iff $\NNC(\kappa,\mu^+,\open)$. \endproclaim

\demo{Proof} First suppose that $\Fr_{\mu}(\kappa,\aleph_0)$ fails, and
let $M$ be a structure with $\mu $ operations and universe $\kappa$
which has no infinite free subset.  Define subsets $A_n$ of\/ $\setfuncs
{\omega}\kappa$ for $n < \omega $ as follows:  for any $s \in  \setfuncs
{\omega}\kappa$, put $s \in  A_n$ iff $s(n) \in  \subalggen
M{\{s(m)\suchthat  m \ne  n\}}$. If $s(m) = s(n)$ but $m \ne  n$, then
$s \in  A_n$; if $s$ is one-to-one, then $s \in  A_n$ for some $n$ since
$M$ has no infinite free subset. Therefore, $\bigcup_{n<\omega}A_n =
\setfuncs {\omega}\kappa$.  Since $\subalggen MS$ is the union of
$\subalggen Ma$ over all finite $a \subseteq  S$, the sets $A_n$ are
open.  Since $M$ has only $\mu $ operations, $|\subalggen MS| \le  \mu $
for any countable $S$, so $A_n$ is $\mu^+$\snug-narrow in the
$n$\snug'th coordinate. Therefore, $\NNC(\kappa,\mu^+,\open)$ fails.

     For the converse, suppose $\NNC(\kappa,\mu^+,\open)$ fails. Let
$\{A_n\suchthat  n < \omega \}$ be a collection of open sets with union
$\setfuncs {\omega}\kappa$ such that $A_n$ is $\mu^+$\snug-narrow in the
$n$\snug'th coordinate.  For each triple $(\alpha,m,n)$ with $\alpha <
\mu $ and $m < n < \omega $, we will define a function $f_{\alpha
mn}\funcfrom  \setfuncs {n-1}\kappa\to \kappa$. Given $g$, $x$, and $y$ such
that $g$ is a function with $x$ in its domain, let $\substitute(g,x,y)$
be the function obtained from~$g$ by replacing the value at~$x$
with~$y$; that is, $\substitute(g,x,y) = (g\setdiff \{(x,g(x))\})\cup
\{(x,y)\}$.  Now suppose $m < n < \omega $ and $\sigma  \in  \setfuncs
n\kappa$.  Let $\sigma' \in \setfuncs {n-1}\kappa$ be the sequence
obtained from~$\sigma$ by deleting the $m$\snug'th coordinate.
Since $A_m$ is $\mu^+$\snug-narrow in the $m$\snug'th
coordinate, we can choose a sequence $\langle \beta_{\alpha}\suchthat
\alpha  < \mu \rangle$ of elements of~$\kappa$ (depending only
on~$\sigma'$ and~$m$, not on~$\sigma(m)$) which includes every
$\beta  < \kappa$ such that $\{\substitute(\sigma,m,\beta)\concat
s\suchthat  s \in  \setfuncs {\omega}\kappa\} \subseteq  A_m$.  Let
$f_{\alpha mn}(\sigma') = \beta_{\alpha}$.  Now let $M$ be the structure
$(\kappa,(f_{\alpha mn})_{\alpha<\mu,\,m<n<\omega})$; clearly $M$ has
$\mu $ operations and cardinality $\kappa$.  Let $S$ be any infinite
subset of~$\kappa$, and choose a one-to-one $s \in \setfuncs {\omega}S$.
There is an $m < \omega $ such that $s \in  A_m$; since $A_m$ is open,
there is an $n < \omega $ such that $\{(s\restrict n)\concat t\suchthat
t \in  \setfuncs {\omega}\kappa\} \subseteq  A_m$, and we may assume $n
> m$.  Let $\sigma'$ be $s \restrict n$ with coordinate~$m$ deleted.
By the definition of $f_{\alpha mn}$, there
must be an $\alpha  < \mu $ such that $f_{\alpha
mn}(\sigma') = s(m)$. But
$\sigma'$~is a sequence of elements of
$S\setdiff\{s(m)\}$, so $S$ cannot be free for $M$.  Therefore, $M$ has
no infinite free subset, so $\Fr_{\mu}(\kappa,\aleph_0)$ fails.
\QED\enddemo

     This equivalence allows us to translate several
results of Devlin and Paris on the free subset problem into
results about $\NNC$:

\proclaim{Corollary 3.4}

{\rm (a)} If $\kappa$ is real-valued measurable, then
$\NNC(\kappa,\lambda,\open)$ for all $\lambda  < \kappa$.

{\rm (b)} The statement $\NNC(\kappa,\mu^+,\open)$ (as an assertion
about~$\kappa$ and~$\mu$) is absolute downward
for transitive models of\/~{\rm ZFC}, and is preserved
under forcing extensions which satisfy the countable chain condition.

{\rm (c)} If $\kappa\to (\omega)^{<\omega}_{2^\mu}$,
then $\NNC(\kappa,\mu^+,\open)$.

{\rm (d)} If $\kappa$
is the least cardinal such that $\kappa\to (\omega)^{<\omega}$, then
$\NNC(\kappa,\lambda,\open)$ for all $\lambda  < \kappa$.

{\rm (e)} If\/ $V = L$ or $V = L[D]$ where $D$ is a normal ultrafilter
over a measurable cardinal, then $\NNC(\kappa,\aleph_1,\open)$
iff $\kappa\to (\omega)^{<\omega}$.
\endproclaim

\demo{Proof}
(a) Devlin~\cite{\Devlin, p.~315}.
(b) Devlin~\cite{\Devlin, pp.~314--316}.
(c) Any homogeneous set for a structure
is free for that structure \cite{\Devlin, p.~314}.
(d) This follows from~(c) and the fact that this~$\kappa$ is
a strong limit cardinal satisfying $\kappa \to (\omega)^{<\omega}_\mu$
for all $\mu < \kappa$ (Silver; see Jech~\cite{\Jech, Lemma~32.9}).
(e) Devlin and Paris~\cite{\DevlinParis, pp.~334-335}.
\QED\enddemo

Therefore, $\NNC(\kappa,\aleph_1,\open)$ implies
$L\models  \kappa\to (\omega)^{<\omega}$.  So the consistency
strength of $\NNC(\kappa,\aleph_1,\open)$ is the same as that of
$\kappa\to (\omega)^{<\omega}$, while the consistency strength of
$\NNC(\kappa,\aleph_1,\Borel)$ lies somewhere between that of $\kappa\to
(\omega)^{<\omega}$ and that of $\kappa\to (\omega_1+\omega)^{<\omega}$.

     Koepke~\cite{\Koepke} uses a measurable cardinal to construct a
model in which $\Fr_{\aleph_0}(\aleph_{\omega},\aleph_0)$ (equivalently,
$\NNC(\aleph_{\omega},\aleph_1,\open)$) holds.
In fact, the properties he proves about this model imply a
stronger result:

\proclaim{Theorem 3.5} If ``there is
a measurable cardinal'' is consistent with\/~{\rm ZFC}, then so is
$\NNC(\aleph_\omega,{<}\aleph_\omega,\Borel)$.
\endproclaim

\demo{Proof} Let $\kappa = \aleph_\omega$.
In the generic extension constructed
by Koepke~\cite{\Koepke}, the following property holds:
for any $f\funcfrom  [\kappa]^{<\omega}\to
2$, there is a sequence $\langle C_i\suchthat  i < \omega \rangle$ such
that $C_i$~is a cofinal subset of $\omega_{2i+2}$ and, for any finite
sequences $\langle i_n\suchthat  n < N\rangle$, $\langle \alpha_n\suchthat
n < N\rangle$, and $\langle \beta_n\suchthat  n < N\rangle$ such
that $i_0 < i_1 < \dots  < i_{n-1} < \omega $ and $\alpha_m,
\beta_m \in  C_{i_m}$, we have $f(\alpha_0,\alpha_1,\dots,\alpha_{n-1})
= f(\beta_0,\beta_1,\dots,\beta_{n-1})$.  The same argument as for
Rowbottom's result that
$\kappa\to (\alpha)^{<\omega}$ implies $\kappa\to
(\alpha)^{<\omega}_{2^{\aleph_0}}$ \cite{\Jech, Lemma~32.8}
can be used to show that the above property actually
holds for any $f\funcfrom  [\kappa]^{<\omega}\to \setfuncs {\omega}2$.

Now suppose we have Borel sets $A_n \subseteq \setfuncs \omega\kappa$
for $n < \omega$ with union $\setfuncs \omega\kappa$, and natural
numbers~$k_n$ for $n < \omega$; we must show that, for some~$n$,
$A_n$~is not $\aleph_{k_n}$\snug-narrow in the $n$\snug'th coordinate.
We may assume $k_0 < k_1 < k_2 < \dotso$.  There is a sequence $\langle
G_m\suchthat m < \omega \rangle$ of open subsets of\/~$\setfuncs
{\omega}\kappa$ such that each of the sets $A_n$ is a Boolean
combination of the sets $G_m$, $m < \omega $.  Define $f\funcfrom
[\kappa]^{<\omega}\to \setfuncs {\omega}2$ by:  $f(\sigma)(m) = 1$ iff
$\{\sigma \concat s\suchthat  s \in \setfuncs {\omega}\kappa\} \subseteq
G_m$. Since $G_m$~is open, for any $s \in \setfuncs\omega\kappa$, we
have $s \in G_m$ iff there is an~$n$ such that $f(s\restrict n)(m) = 1$.
Find $\langle C_i\suchthat  i < \omega \rangle$ as in the preceding
paragraph.  Then, if $s$ and $s'$ are sequences of length $\omega $ such
that $s(i), s'(i) \in  C_{k_i}$ for each $i < \omega $, then
$f(s\restrict n) = f(s'\restrict n)$ for all $n$, so $\{m\suchthat  s
\in  G_m\} = \{m\suchthat s' \in  G_m\}$, and since the $A_n$\snug's are
Boolean combinations of the $G_m$\snug's, $\{n\suchthat s \in  A_n\} =
\{n\suchthat  s' \in  A_n\}$.  Hence, there is a fixed~$n$ such that $s
\in  A_n$ for all such $s$; since there is a collection of
$\aleph_{2k_n+2} > \aleph_{k_n}$ such~$s$\snug's which differ only at
coordinate~$n$, $A_n$ is not $\aleph_{k_n}$\snug-narrow in the
$n$\snug'th coordinate, and we are done. \QED\enddemo

Note that the argument here actually gives
$\NNC(\aleph_\omega,{<}\aleph_\omega,S)$ where $S$~is the
collection of sets which are expressible as Boolean combinations
of countably many open sets; this collection includes the Borel
sets and many other sets as well.

By the way, standard chain-condition and closure arguments
(see Shelah's version~\cite{\Shelah}) show that $\aleph_{\omega}$~is a
strong limit cardinal in this model.

\head 4.  Forcing and Narrow Coverings \endhead

In this section, we will show that, at least for most $\kappa$
and $\lambda$, the properties $\NNC(\kappa,\lambda,\allowbreak
F_\sigma)$ and
$\NNC(\kappa,\lambda,\Borel)$ are preserved under forcing to add any
number of Cohen reals or random reals.
This will prove the
consistency of Mr\'owka's hypothesis $S(\aleph_0)$,
given a suitable large cardinal.

\proclaim{Theorem 4.1} Let $M[G]$ be a generic extension of a ground
model~$M$ of\/~{\rm ZFC}, obtained by the standard forcing to add either
any number of Cohen reals or any number of random reals. Let $\kappa$
and~$\lambda$ be cardinals in~$M$, with $\cf \lambda > \omega$. If
$\NNC(\kappa,\lambda,F_\sigma)$ is true in~$M$, then it is true
in~$M[G]$.  The same holds for $\NNC(\kappa,\lambda,\Borel)$.
\endproclaim

\proclaim{Corollary 4.2} If $(\exists\kappa)(
\kappa\to(\omega_1+\omega)^{<\omega})$ is consistent with\/~{\rm ZFC},
then so are $\NNC(2^{\aleph_0},\allowbreak \aleph_1,F_\sigma)$
(i.e.,~$S(\aleph_0)$) and $\NNC(2^{\aleph_0},\aleph_1,\Borel)$.
\endproclaim

\demo{Proof} Start with a model where $\kappa$~has the specified
partition property, so that Corollary~3.2 applies, and
add $\kappa$~Cohen or random reals. \QED\enddemo

Note that, if we start with a measurable cardinal~$\kappa$ and
add $\kappa$~random reals, we get a model where $\kappa$~is
real-valued measurable and $\NNC(\kappa,{<}\kappa,\Borel)$ holds.
It is still open whether $\NNC(\kappa,{<}\kappa,\Borel)$ actually follows
from real-valued measurability of~$\kappa$.

\proclaim{Corollary 4.3} If ``there is a measurable cardinal'' is
consistent with\/~{\rm ZFC}, then so is
$(2^{\aleph_0}=\aleph_{\omega+1}) +
\NNC(\aleph_\omega,{<}\aleph_\omega,\Borel)$. \endproclaim

\demo{Proof} Start with a model obtained from Theorem~3.5, and
add $\aleph_{\omega+1}$~Cohen or random reals. \QED\enddemo

So we have a model where $S(\aleph_0)$~holds and
$2^{\aleph_0}=\aleph_{\omega+1}$.
Note that $\aleph_{\omega+1}$ is the smallest
possible value for~$2^{\aleph_0}$ in a model of~$S(\aleph_0)$, since, by
Corollary~2.8, $\NNC(\kappa,\aleph_1,\clopen)$
cannot hold for $\kappa < \aleph_\omega$ (and since K\"onig's
theorem implies that $2^{\aleph_0}$ cannot be equal to~$\aleph_\omega$).

The proof of Theorem~4.1 for random reals is somewhat simpler
than that for Cohen reals, so it will be given first.
In both cases the $F_\sigma$~version is given separately because
the full Borel version requires additional work.

All of the arguments below are carried out within the ground model~$M$.
The forcing partial orders will be written so that
$p \le q$ means that $p$~is a stronger condition than~$q$.

The idea of the proof is to show that a counterexample to
$\NNC(\kappa,\lambda,S)$ (where $S$~is `$F_\sigma$' or `Borel')
in the generic extension can be turned into a counterexample in
the ground model.  To say that there is a counterexample in the
extension means that there exist names~$\An_n$ for $n < \omega$
and a forcing condition $p_0$ (in the generic filter) such that
$$p_0 \forces \bigcup_{n < \omega} \An_n =
\setfuncs \omega\kappa\tag 4.1$$
and, for each $n<\omega$,
$$p_0 \forces \text{$\An_n$ has property $S$ and is $\lambda$\snug-narrow 
in the $n$\snug'th coordinate}.\tag 4.2$$

One could get a narrow covering of the~$\setfuncs \omega\kappa$ of the
ground model by simply restricting the sets~$\An_n$ to this space,
but the resulting sets would probably not be in the ground model.
However, given a name~$\An$, we {\sl can} define in the
ground model a set which will definitely include the set named by~$\An$:

\definition{Definition 4.4} Given a name~$\An$ and a
forcing condition~$p_0$, the set of {\it potential members of~$\An$
(assuming~$p_0$)} is the set of all~$x$ (in the ground model) such
that there exists $p \le p_0$ such that $p \forces x \in \An$.
\enddefinition

The ``(assuming~$p_0$)'' will usually be omitted since $p_0$~will
be clear from the context.

Suppose we have $p_0$ and~$\An_n$ satisfying \thetag{4.1}
and~\thetag{4.2}.  Let $B_n$ be the set of potential members
of~$\An_n$.  Then $B_n \subseteq \setfuncs \omega\kappa$ for each~$n$.
Also, for any $s \in \setfuncs \omega\kappa$, we have
$p_0 \forces (\exists n)\,s \in \An_n$, so, for some~$n$
and some $p \le p_0$, $p \forces s \in \An_n$.  Therefore,
$\bigcup_{n < \omega} B_n = \setfuncs \omega\kappa$.
We next show that the set~$B_n$ is $\lambda$\snug-narrow
in the $n$\snug'th coordinate.

\proclaim{Lemma 4.5} Let $P$ be a notion of forcing (partial ordering)
with the countable chain condition, and let $\kappa$ and~$\lambda$
be cardinals such that $\cf \lambda > \omega$.  Suppose that
$p_0 \in P$ and $\An$ is a $P$\snug-name such that
$p_0 \forces \An \subseteq \setfuncs \omega\kappa$.  If
$$p_0 \forces \text{$\An$ is $\lambda$\snug-narrow in the
$n\!$\snug'th coordinate},$$ then the set of potential members
of~$\An$ is $\lambda$\snug-narrow in the
$n\!$\snug'th coordinate. \endproclaim

\demo{Proof} Let $B$ be the set of potential members of~$\An$. Let $s$
be a member of\/ $\setfuncs \omega\kappa$; we must see that $B$ contains
fewer than $\lambda$ points on the line $$\{s' \in \setfuncs
\omega\kappa\suchthat \text{$s(m) = s'(m)$ for $m \ne n$}\}.$$ In other
words, letting $\substitute(s,n,\alpha)$ denote the sequence $s$ with
entry number $n$ replaced with~$\alpha$ (as in Section~3), we must show
that $\{\alpha < \kappa\suchthat \substitute(s,n,\alpha) \in B_n\}$ has
size less than $\lambda$.

Since $p_0 \forces \text{$\An$ is $\lambda$\snug-narrow in the $n$\snug'th
coordinate}$, there exist $P$-names $\Pname \beta$ and $\Pname f$ such
that $p_0$ forces that $\Pname \beta < \lambda$ and $\Pname f$ is a
function with domain $\Pname\beta$ enumerating the ordinals $\Pname\alpha$
such that $\substitute(s,n,{\Pname\alpha}) \in \An$.  By the usual
countable chain condition argument (choosing a maximal antichain of
conditions below $p_0$ which decide the value of $\Pname\beta$), there
is a countable set $S$ of ordinals less than $\lambda$ such that $p_0
\forces \Pname\beta \in S$.  Let $\beta_0$ be the least upper bound of
$S$; since $\lambda$ has uncountable cofinality, $\beta_0 < \lambda$.

By the same argument, for each $\gamma < \beta_0$, there is a countable
set $W_\gamma \subset \kappa$ such that $p_0$~forces $\Pname f(\gamma)$,
if it exists, to be in $W_\gamma$.  Let $W = \bigcup_{\gamma < \beta_0}
W_\gamma$.  Then, for any ordinal $\alpha < \kappa$, if
$\alpha \notin W$, then $p_0$~forces that $\alpha$ is not in the
range of $\Pname f$, so $p_0 \forces \substitute(s,n,\alpha) \notin \An$,
so $\substitute(s,n,\alpha) \notin B$.  Therefore,
$\{\alpha < \kappa\suchthat \substitute(s,n,\alpha) \in B\} \subseteq W$;
since $|W| \le |\beta_0|\cdot\aleph_0 < \lambda$, we are done.
\QED\enddemo

So the sets $B_n$ form a $\lambda$\snug-narrow covering of\/~$\setfuncs
\omega\kappa$ (in the ground model).  If we can show that $$(p_0 \forces
\text{$\An_n$ has property~$S$}) \implies \text{$B_n$~has property~$S$}
\tag 4.3$$ (where $S$~is `$F_\sigma$' or `Borel'), then we will have
completed the proof that a counterexample to $\NNC(\kappa,\lambda,S)$ in
the generic extension gives a counterexample to $\NNC(\kappa,\lambda,S)$
in the ground model.

We first consider the case of random real forcing.  Actually, the
argument applies more generally to any forcing notion which is a measure
algebra.  (A measure algebra is a complete Boolean algebra with an
associated nonzero $\sigma$\snug-additive probability function; see Jech
\cite{\Jech, p.~421} for details.  In particular, random real forcing is
given by a measure algebra, and any measure algebra has the countable
chain condition.)  But Maharam~\cite{\Maharam} has shown that this is
not much of a generalization.

Since we are using a complete Boolean algebra as the forcing
notion, every sentence~$\varphi$ of the forcing language has an
associated Boolean value~$\| \varphi \|$.

In the usual way, any closed set $F \subseteq \setfuncs \omega\kappa$
can be expressed in the form~$[T]$, the set of infinite branches
through some tree $T \subseteq \setfuncs {<\omega}\kappa$:
given~$F$, let $T$~be the set of finite sequences~$\sigma$ such that
some member of~$F$ starts with~$\sigma$.  Conversely, any set of
the form~$[T]$ is closed.

\proclaim{Lemma 4.6} Let $P$ be a notion of forcing obtained from a
measure algebra.  Suppose that $p_0 \in P$ and $\An$ is a $P$\snug-name
such that $p_0 \forces \An \subseteq \setfuncs \omega\kappa$.  If $p_0
\forces \text{$\An$ is $F_\sigma$}$, then the set of potential members
of~$\An$ is $F_\sigma$. \endproclaim

\demo{Proof} Let $B$ be the set of potential members of~$\An$, and let
$\mu$ be the probability function for the measure algebra. Since every
nonzero member of the Boolean algebra is given nonzero measure by~$\mu$,
we can rewrite the definition of~$B$ as follows: $$B = \{s \in \setfuncs
\omega\kappa \suchthat \mu(p_0 \cdot \|s \in \An \|) > 0 \}.$$ We must
see that this set is~$F_\sigma$.

Any $F_\sigma$ subset of\/~$\setfuncs \omega\kappa$ is a countable union
of closed sets, each of which can be expressed in the form~$[T]$ for
some tree $T \subseteq \setfuncs {<\omega}\kappa$; furthermore, we may
assume that the union is an increasing union.  Therefore, there are
$P$\snug-names~$\Tn_m$ for $m<\omega$ such that $$p_0 \forces
\text{$\Tn_m$ is a tree, $[\Tn_m] \subseteq [\Tn_{m+1}]$, and $\An =
\bigcup_{m<\omega} [\Tn_m]$}.$$ The Boolean value $p_0 \cdot \| s \in
\An\|$ is the sum (least upper bound) of the Boolean values $p_0 \cdot
\| s \in [\Tn_m]\|$, which form an increasing sequence; since the
$\sigma$\snug-additivity of~$\mu$ implies continuity with respect to
increasing limits, we get $$\mu(p_0 \cdot \| s \in \An\|) = \sup_m
\mu(p_0 \cdot \| s \in [\Tn_m]\|).$$ Similarly, $p_0 \cdot \| s \in
[\Tn_m]\|$ is the decreasing limit of the Boolean values $p_0 \cdot \|
s\restrict k \in \Tn_m\|$, so $$\mu(p_0 \cdot \| s \in [\Tn_m]\|) =
\inf_k \mu(p_0 \cdot \| s\restrict k \in \Tn_m\|).$$ Therefore, $$\align
s \in B &\iff \mu(p_0 \cdot \|s \in \An \|) > 0 \cr &\iff \sup_m \inf_k
\mu(p_0 \cdot \| s\restrict k \in \Tn_m\|) > 0 \cr &\iff (\exists m)\,\,
\inf_k \mu(p_0 \cdot \| s\restrict k \in \Tn_m\|) > 0 \cr &\iff (\exists
m)(\exists\eps)(\forall k)\,\, \mu(p_0 \cdot \| s\restrict k \in
\Tn_m\|) > \eps, \endalign$$ where $\eps$~varies over the positive
rational numbers. Since the condition $\mu(p_0 \cdot \| s\restrict k \in
\Tn_m\|) > \eps$ depends only on~$s\restrict k$, the set of~$s$
satisfying this condition is clopen. Therefore, $B$~is~$F_\sigma$, as
desired. \QED\enddemo

This completes the proof of Theorem~4.1 for the case of $F_\sigma$~sets
and random real forcing, which suffices for the relative consistency
of~$S(\aleph_0)$.

In order to do the case of Borel sets and random real forcing (in fact,
measure algebra forcing), we will need to work with codes of Borel sets,
and it will be convenient to work with these codes in a slightly more
restrictive way than usual.

Define the Borel hierarchy as usual: $\SIGMA^0_1$~sets are open sets,
$\PI^0_1$~sets are closed sets, $\SIGMA^0_\alpha$~sets for $1 < \alpha
< \omega_1$ are countable unions of\/ $\PI^0_\beta$~sets for (possibly
varying) $\beta < \alpha$, and $\PI^0_\alpha$~sets for $1 < \alpha
< \omega_1$ are countable intersections of\/ $\SIGMA^0_\beta$~sets for
$\beta < \alpha$.  So a $\PI^0_\alpha$~set is just the complement of
a $\SIGMA^0_\alpha$~set.

Every closed set $F \subseteq \setfuncs \omega\kappa$ is a countable
intersection of clopen sets: if $F = [T]$, then $F = \bigcap_{n <
\omega} C_n$ where $C_n = \{ s \in \setfuncs \omega\kappa \suchthat s
\restrict n \in T\}$. Similarly, every open set is a countable union of
clopen sets~$C_n$.  From these facts, one can inductively prove the
usual inclusions: $\SIGMA^0_\beta \cup \PI^0_\beta \subseteq
\SIGMA^0_\alpha \cap \PI^0_\alpha$ for $\beta < \alpha$. Also, the
collections $\SIGMA^0_\alpha$ and~$\PI^0_\alpha$ are closed under finite
unions and intersections.

Let $n,m \mapsto \pair(n,m)$ be a one-to-one function from
$\omega \times \omega$ to $\omega\setdiff\{0\}$ such that
$\pair(n,m)$ increases with~$m$ for each fixed~$n$.  As usual, this
allows us to code up infinitely many $\omega$\snug-sequences into
one, and conversely extract from one sequence~$s$ the infinitely
many subsequences~$(s)_n$ defined by $(s)_n(m) = s(\pair(n,m))$.

\definition{Definition 4.7}

{\rm (a)} A {\it Borel code} (of level~$\alpha$) is a sequence
$c \in\setfuncs \omega\omega$ such that either $c(0) \ge 2$ or,
for all~$n$, $(c)_n$ is a Borel code (of level $<\alpha$).

{\rm (b)} Given a Borel code~$c$ and a sequence of sets
$\langle Z_m \suchthat m < \omega \rangle$, define the
set $c(\langle Z_m \suchthat m < \omega \rangle)$ as follows:
if $c(0) \ge 2$, then
$$c(\langle Z_m \suchthat m < \omega \rangle) = Z_{c(1)};$$
if $c(0) = 0$, then
$$c(\langle Z_m \suchthat m < \omega \rangle) = \bigcup_{n < \omega}
(c)_n(\langle Z_m \suchthat m < \omega \rangle);$$
if $c(0) = 1$, then
$$c(\langle Z_m \suchthat m < \omega \rangle) = \bigcap_{n < \omega}
(c)_n(\langle Z_m \suchthat m < \omega \rangle).$$

{\rm (c)} Given $c$ and $\langle Z_m \suchthat m < \omega \rangle$
as above, where $Z_m \subseteq \setfuncs \omega\kappa$, say that
$\langle Z_m \suchthat m < \omega \rangle$ is {\it good} for~$c$
iff:
\roster
\item for all $m$, $Z_m$ is a clopen set for which membership
depends only on the first $m$~coordinates (i.e., if $s \in Z_m$
and $s\restrict m = s' \restrict m$, then $s' \in Z_m$); and
\item during the recursive computation of
$c(\langle Z_m \suchthat m < \omega \rangle)$,
all unions are increasing unions and
all intersections are decreasing intersections.
\endroster
\enddefinition

Now a very slight variation of a standard argument gives:

\proclaim{Lemma 4.8} For each nonzero $\alpha < \omega_1$, there is
a universal\/ $\SIGMA^0_\alpha$~code, i.e.,
a Borel code~$c$ such that every $\SIGMA^0_\alpha$~subset
of\/~$\setfuncs \omega\kappa$ is of the form
$c(\langle Z_m \suchthat m < \omega \rangle)$ for some sequence
$\langle Z_m \suchthat m < \omega \rangle$ which is good for~$c$
(and the converse: $c(\langle Z_m \suchthat m < \omega \rangle)$
is $\SIGMA^0_\alpha$ for any clopen sets~$Z_m$).
Similarly, for each~$\alpha$ there is a universal\/ $\PI^0_\alpha$~code.
\endproclaim

\demo{Proof}
To get a universal $\PI^0_1$~code, just define~$c$ so that
$c(0) = 1$, $(c)_n(0) = 2$, and $(c)_n(1) = n$ for all~$n$;
this gives $c(\langle Z_m \suchthat m < \omega \rangle) =
\bigcap_{m < \omega} Z_m$.  This works because, given any
closed set~$[T]$, we can let $Z_m = \{s\suchthat s\restrict m \in T\}$
to generate~$[T]$ from~$c$.  A similar argument with the
complements gives a universal $\SIGMA^0_1$~code --- just let
$c(0)$ be~$0$ instead of\/~$1$.

Now suppose $\alpha > 1$.  If $\alpha$~is a limit ordinal, let
$\alpha_0,\alpha_1,\dotsc$ be a strictly increasing sequence of
ordinals converging to~$\alpha$; if $\alpha = \beta+1$, let
$\alpha_n = \beta$ for all~$n$.  Apply the inductive hypothesis
to get a universal $\PI^0_{\alpha_n}$~code~$c_n$ for each~$n$.
Let $c'_n$ be $c_n$ with all references to the $m$\snug'th
given clopen set replaced with references to the $\pair(n,m)$\snug'th
clopen set, so that $$c'_n(\langle Z_m \suchthat m < \omega \rangle)
= c_n(\langle Z_{\pair(n,m)} \suchthat m < \omega \rangle)$$
for any sets~$Z_m$.  Now we can find~$c$ so that
$c(0) = 0$ and $(c)_n = c'_n$ for all~$n$.

This~$c$ is a universal $\SIGMA^0_\alpha$ code.  Given any
$\SIGMA^0_\alpha$~set~$A$, find sets~$B_j$ for $j <\omega$
with union~$A$ so that each~$B_j$ is $\PI^0_\beta$ for some
$\beta < \alpha$.  Then $A$~is the increasing union of
the sets $B'_k = \bigcup_{j < k} B_j$, and each~$B'_k$ is
also $\PI^0_\beta$ for some $\beta<\alpha$ (and $B'_0 = \nullset$).
We can find a nondecreasing sequence $k_0,k_1,k_2,\dotsc$
of natural numbers tending to infinity so slowly that
$B'_{k_n}$~is a $\PI^0_{\alpha_n}$~set for all~$n$.
For each~$n$, choose a sequence $\langle Z_m^{(n)} \suchthat m <
\omega \rangle$ which is good for~$c_n$ so that
$c_n(\langle Z_m^{(n)} \suchthat m < \omega \rangle) = B'_{k_n}$.
Define $\langle Z_m \suchthat m < \omega \rangle$ so that
$Z_{\pair(n,m)} = Z_m^{(n)}$ for all $m$ and~$n$, and
$Z_k = \nullset$ for all remaining~$k$; then
$\langle Z_m \suchthat m < \omega \rangle$ is good for~$c$
(here we use the fact that $\pair(n,m)$ increases with~$m$,
so that $\pair(n,m) \ge m$) and
$c(\langle Z_m \suchthat m < \omega \rangle) = A$,
as desired.

The argument for~$\PI^0_\alpha$ is the same.
\QED\enddemo

Note that the construction of the universal $\SIGMA^0_\alpha$
or $\PI^0_\alpha$~code~$c$
is very absolute, once one has chosen a cofinal $\omega$\snug-sequence
for each limit ordinal~$\le\alpha$.  In particular, if $c$~is
constructed for~$\alpha$ in a ground model~$M$, then the same~$c$
will work in any extension~$M[G]$
of~$M$, although there will probably be more good
sequences to apply it to.

\proclaim{Lemma 4.9} Let $P$ be a notion of forcing obtained from
a measure
algebra, with associated probability function~$\mu$.  Suppose that
$p_0 \in P$ and $\An$ is a $P$\snug-name such that
$p_0 \forces \An \subseteq \setfuncs \omega\kappa$.  If
$p_0 \forces \text{$\An$ is Borel}$, then the function
$s \mapsto \mu(p_0 \cdot \|s \in \An\|)$ is a Borel-measurable
function from\/~$\setfuncs \omega\kappa$ to~$[0,1]$. \endproclaim

\demo{Proof} We know that $$p_0 \forces (\text{$\An$ is
$\SIGMA^0_{\Pname\alpha}$ for some $\Pname\alpha < \omega_1$}).$$ By the
usual countable chain condition argument, the set of $\beta < \omega_1$
such that $(\exists p{\le}p_0)\, p \forces \Pname\alpha = \beta$ is
countable, and if we choose $\gamma < \omega_1$ to be greater than all
such~$\beta$, then we will have $p_0 \forces \text{$\An$ is
$\SIGMA^0_\gamma$}$. Let $c$ be a universal $\SIGMA^0_\gamma$~code (in
the ground model); then there exist names~$\Zn_m$ for $m < \omega$ such
that $$p_0 \forces (\text{$\langle \Zn_m \suchthat m < \omega \rangle$
is good for $c$ and $c(\langle \Zn_m \suchthat m < \omega \rangle) =
\An$}). \tag 4.4$$

So we must show: if we have a Borel code~$c$ and names~$\Zn_m$
so that \thetag{4.4}~holds, then the function~$f$ defined by
$f(s) = \mu(p_0 \cdot \|s \in \An\|)$ is Borel-measurable.
The proof of this is by induction on the complexity of~$c$.

If $c(0) \ge 2$, then
$\An = c(\langle \Zn_m \suchthat m < \omega \rangle)$ is just
$\Zn_{c(1)}$.  By the goodness assumption, membership of~$s$
in~$\Zn_{c(1)}$ depends only on $s \restrict (c(1))$,
so $f(s)$ depends only on~$s\restrict(c(1))$
and hence is a Borel-measurable (even clopen-measurable) function
of~$s$.

If $c(0) = 0$, then $\An$ is the increasing union of the sets
$\An_n = (c)_n(\langle \Zn_m \suchthat m < \omega \rangle)$, so
the Boolean value $p_0 \cdot\|s \in \An\|$ is the increasing limit of
the Boolean values $p_0 \cdot\|s \in \An_n\|$.  Hence,
$f(s)$ is the increasing limit of the numbers
$f_n(s) = \mu(p_0 \cdot\|s \in \An_n\|)$; the functions~$f_n$
are Borel-measurable by the inductive hypothesis, so
$f$~is Borel-measurable.

Similarly, if $c(0) = 1$, then $f$~is a decreasing limit
of a sequence of Borel-measurable functions,
so $f$~is Borel measurable.
\QED\enddemo

In particular, the set of potential members of~$\An$ is Borel, since
this set is just $\{s\suchthat \mu(p_0 \cdot \|s \in \An\|)>0\}$. This
shows that \thetag{4.3} holds for $S = \Borel$, which completes the
proof of the random real version of Theorem~4.1.

(If one keeps track of the Borel levels in Lemma~4.9, one finds:
if $p_0 \forces \text{$\An_n$ is $\SIGMA^0_\alpha$}$, then
$\{s\suchthat \mu(p_0 \cdot \|s \in \An\|)>r\}$ is~$\SIGMA^0_\alpha$;
if $p_0 \forces \text{$\An_n$ is $\PI^0_\alpha$}$, then
$\{s\suchthat \mu(p_0 \cdot \|s \in \An\|)>r\}$ is~$\SIGMA^0_{\alpha+1}$.
Hence, the property $\NNC(\kappa,\lambda,\SIGMA^0_\alpha)$
is preserved by measure algebra forcing if $\cf\lambda>\omega$.)

Now let $P \in M$ be the forcing notion for adding a certain number of
Cohen reals.  We may take $P$ to be the set of all finite partial
functions from some ordinal $\theta \in M$ to $\{0,1\}$, where, given
two such functions $p,q$, we have $p \le q$ iff $q \subseteq p$.
This~$P$ is called $\Cohentheta$.

Again, for the proof that $\NNC(\kappa,\lambda,S)$ (where $S$~is
`$F_\sigma$' or `Borel') is preserved under forcing with~$P$, suppose
that we have a condition~$p_0$ and names~$\An_n$ for $n < \omega$ such
that \thetag{4.1} and~\thetag{4.2} hold.  Let $B_n$ be the set of
potential members of~$A_n$.  Then the sets~$B_n$ form a
$\lambda$\snug-narrow covering of\/~$\setfuncs \omega\kappa$ as before,
and it remains to show that \thetag{4.3}~holds in order to get a
counterexample to $\NNC(\kappa,\lambda,S)$ in the ground model.

For any statement $\varphi$ in the forcing language for~$P$, one can find
a maximal antichain~$D$ of conditions in $P$ which either force
$\varphi$ or force $\neg\varphi$.  Since $P$ has the countable chain
condition, $D$ is countable.  Let $C$ be the union of the domains
of the members of $D$; then $C$~is a countable subset of $\theta$.
Now, for any condition $q$, $q \forces \varphi$ if and only if
$q$ is incompatible with all members of $D$ which force $\neg\varphi$;
it follows that $q \forces \varphi$ iff $q\restrict C \forces \varphi$.
Call a set $C \subseteq \theta$
with this property a {\it support} of $\varphi$.  Note that,
if $C \subseteq C' \subseteq \theta$ and $S$ is a support of $\varphi$,
then $C'$ is a support of $\varphi$ (since $q \restrict C' \le
q \restrict C$).

For each $\varphi$, let $\support\varphi$ be a countable support
of $\varphi$; it does not matter which one is chosen.  (One can
just take the first one in some fixed well-ordering of the power
set of $\theta$.  Or, in fact, one can show that, for the case of
this particular forcing notion, each $\varphi$ has a unique
{\sl minimal} support, which can be chosen as $\support\varphi$.)
We will assume that, if $\varphi$ and~$\varphi'$ are equivalent
(i.e., for all~$p$, $p \forces \varphi \leftrightarrow \varphi'$),
then $\support\varphi = \support{\varphi'}$.

For any set $C \subseteq \theta$, let $P \cfrestrict C$ be the set
of members of $P$ whose domains are subsets of~$C$.  Note that,
if $C$ is countable, then $P \cfrestrict C$ is countable.

For any $m < \omega$, the sets $\{ s \in \setfuncs \omega\kappa
\suchthat \tau \subseteq s\}$ for $\tau \in \setfuncs m\kappa$ form a
partition of\/~$\setfuncs \omega\kappa$ into clopen pieces.  Hence, a
subset of\/~$\setfuncs \omega\kappa$ is closed if and only if its
intersection with each of these pieces is closed, and the same holds
for~$F_\sigma$. In other words, if we define $X \wkrestrict \tau$ (for
$X \subseteq \setfuncs \omega\kappa$) to be $\{s \in \setfuncs
\omega\kappa \suchthat \tau \concat s \in X\}$, then $X$~is~$F_\sigma$
if and only if $X\wkrestrict\tau$ is~$F_\sigma$ for all
$\tau\in\setfuncs m\kappa$.

\proclaim{Lemma 4.10} Let $P = \Cohentheta$.  Suppose that $p_0 \in P$
and $\An$ is a $P$\snug-name such that $p_0 \forces \An \subseteq
\setfuncs \omega\kappa$.  If $p_0 \forces (\text{$\An$ is $F_\sigma$})$,
then the set of potential members of~$\An$ is $F_\sigma$. \endproclaim

\demo{Proof} Let $B$ be the set of potential members of~$\An$.
As in Lemma~4.6, there are
$P$\snug-names~$\Tn_m$ for $m<\omega$ such that
$$p_0 \forces \text{$\Tn_m$ is a tree, $[\Tn_m] \subseteq [\Tn_{m+1}]$,
and $\An = \bigcup_{m<\omega} [\Tn_m]$}.$$

For each finite sequence $\sigma \in \setfuncs {<\omega}\kappa$, define
a set $S_\sigma \subseteq \theta$ as follows:
$$S_\sigma = \domain(p_0) \cup \bigcup_{k \le \ell(\sigma)}
\bigcup_{m < \omega} \support{\sigma\restrict k \in \Tn_{m}}.$$
So $S_\sigma$ is countable, and $S_\sigma \subseteq S_\tau$ if
$\sigma \subseteq \tau$.

If $s \in B$, then there is a condition $p' \le p_0$ such that $p'
\forces s \in \An$.  Then there must exist $p \le p'$ and a specific $m
< \omega$ such that $p \forces s \in [\Tn_{m}]$. Equivalently, $p
\forces s\restrict k \in \Tn_{m}$ for all $k < \omega$.  Now, if $C$ is
the countable set $\bigcup_{j < \omega} S_{s \restrict j}$, then $C$~is
a support of $(s\restrict k \in \Tn_{m})$ for all~$k$, so $p\restrict C
\forces s\restrict k \in \Tn_{m}$ for all $k$.  Since $C$ is the
increasing union of the sets $S_{s \restrict j}$, and the domain of $p$
is finite, we actually have $p \restrict C = p \restrict S_{s \restrict
j}$ for some $j$. Also, since $S_{s \restrict j}$ includes the domain of
$p_0$, we still have $p \restrict S_{s \restrict j} \le p_0$. Therefore,
if $s \in B$, then there exist $m,j < \omega$ and $p \in P \cfrestrict
S_{s \restrict j}$ such that $p \le p_0$ and $(\forall k)\,p \forces
s\restrict k \in \Tn_{m}$. The converse of this statement is clearly
true as well. So $B = \bigcup_{m,j < \omega} \hat B_{m,j}$, where $$\hat
B_{m,j} = \{s \suchthat (\exists p{\in}P \cfrestrict S_{s \restrict j})
\,\, p \le p_0 \,\&\, (\forall k)\, p \forces s \restrict k \in
\Tn_{m}\}.$$

Now, given $m < \omega$ and $\tau \in \setfuncs {<\omega}\kappa$, let
$$\hat B_{m}^{(\tau)} = \{s \suchthat (\exists p{\in}P \cfrestrict
S_\tau) \,\, p \le p_0 \,\&\, (\forall k)\, p \forces s \restrict k \in
\Tn_{m}\}.$$ Then $\hat B_{m}^{(\tau)}$ is explicitly a countable union
(over $p$) of a countable intersection (over $k$) of clopen sets, so it
is an $F_\sigma$ set.  But, for all $\tau \in \setfuncs j\kappa$, we
have $\hat B_{m,j} \wkrestrict \tau = \hat B_{m}^{(\tau)} \wkrestrict
\tau$, so $\hat B_{m,j} \wkrestrict \tau$ is~$F_\sigma$; hence, $\hat
B_{m,j}$ is $F_\sigma$.  Therefore, $B$ is $F_\sigma$. \QED\enddemo

So \thetag{4.3} holds for $S=F_\sigma$.

For the Borel case, it will be convenient to change the Borel coding
definitions given earlier (Definition~4.7) so as to use intersections
and complements instead of intersections and unions. This means that,
when $c(0) = 0$, we will have $$c(\langle Z_m \suchthat m < \omega
\rangle) = \setfuncs \omega\kappa \setdiff (c)_0(\langle Z_m \suchthat m
< \omega \rangle).$$ The results proved earlier about Borel codes, such
as the existence of universal $\SIGMA^0_\alpha$~codes, go through as
before.

As we did for $F_\sigma$ sets, we can show that, for any
set $X \subseteq \setfuncs \omega\kappa$ and any $\alpha$ and~$n$,
$X$~is~$\SIGMA^0_\alpha$ if and only if $X \restrict \tau$
is~$\SIGMA^0_\alpha$ for all $\tau \in \setfuncs n\kappa$; the same holds
for~$\PI^0_\alpha$.  (This is proved by induction on~$\alpha$,
with a little care at limit stages.  Alternatively, one can show easily
by induction on Borel codes~$c$ that
$$c(\langle Z_m \suchthat m < \omega \rangle) \wkrestrict \tau
= c(\langle Z_m \wkrestrict \tau \suchthat m < \omega \rangle)$$
for any $\tau$ and any sets~$Z_m$; then apply this
to the case of a universal $\SIGMA^0_\alpha$ or $\PI^0_\alpha$~code.)

Just as for the random real case, we see that, if $$p_0 \forces
\text{$\An$ is a Borel subset of\/ $\setfuncs \omega\kappa$},$$ then we
can find a Borel code~$c$ (in the ground model) and a sequence of
names~$\Zn_m$ such that $$p_0 \forces \text{$\langle \Zn_m \suchthat m <
\omega \rangle$ is good for $c$ and $c(\langle \Zn_m \suchthat m <
\omega \rangle) = \An$}.$$ In fact, we can ensure, by modifying the
names~$\Zn_m$ if necessary, that $\nullset$ (the weakest condition
in~$P$) forces ``$\langle \Zn_m \suchthat m < \omega \rangle$ is good
for~$c$.''  It follows that $\support{s \in \Zn_m}$ depends only on~$s
\restrict m$, not on the rest of~$s$.

\proclaim{Lemma 4.11} Let $P = \Cohentheta$.  Suppose that
$c$~is a Borel code (in terms of intersections and complements,
as above) and $\langle \Zn_m \suchthat m < \omega \rangle$
is a sequence of names for subsets of\/~$\setfuncs \omega\kappa$
such that $$\nullset \forces
\text{$\langle \Zn_m \suchthat m < \omega \rangle$ is good
for~$c$}.$$
Then:

{\rm (a)} For any $s \in \setfuncs \omega\kappa$, $C_s = \bigcup_{j <
\omega} S_{s\restrict j}$ is a support for $(s \in c(\langle \Zn_m
\suchthat m < \omega \rangle))$, where $$S_{s\restrict j} = \bigcup_{m
\le j} \support{s \in \Zn_m}.$$

{\rm (b)} 
There is an ordinal
$\alpha < \omega_1$ such that, for each $p \in P$, the set\/
$\{s \in \setfuncs \omega\kappa\suchthat p \forces s \in c(\langle \Zn_m
\suchthat m < \omega \rangle)\}$ is $\PI^0_\alpha$.
\endproclaim

(The notation $S_{s \restrict j}$ makes sense because $\support{s \in
\Zn_m}$ depends only on $s \restrict m$; in other words, $S_\tau$~is
well-defined for $\tau\in\setfuncs {<\omega}\kappa$.)

\demo{Proof} Induct on~$c$.  If $c(0) \ge 2$,
then $c(\langle \Zn_m \suchthat
m < \omega \rangle)$ is just~$\Zn_m$ for $m = c(1)$, so (a)~is obvious;
for (b), the specified set is actually clopen (and hence~$\PI^0_1$)
since membership of~$s$ in~$\Zn_m$ depends only on~$s\restrict m$.

If $c(0) = 1$, then $c(\langle \Zn_m \suchthat
m < \omega \rangle)$ is the intersection of the sets
$(c)_n(\langle \Zn_m \suchthat m < \omega \rangle)$,
so $p \forces s \in c(\langle \Zn_m \suchthat
m < \omega \rangle)$ if and only if $p \forces s \in 
(c)_n(\langle \Zn_m \suchthat m < \omega \rangle)$ for all~$n$.
Now (a) and~(b) for~$c$ follow easily from the corresponding
facts for~$(c)_n$.  (The~$\alpha$ for~$c$ is the supremum of the
corresponding ordinals for~$(c)_n$.)

Now suppose $c(0) = 0$, so $c(\langle \Zn_m \suchthat
m < \omega \rangle)$ is the complement of
$(c)_0(\langle \Zn_m \suchthat m < \omega \rangle)$.
The induction hypothesis states that (a) and~(b) hold
for~$(c)_0$.  We now get
$$\align &p \forces s \in c(\langle \Zn_m \suchthat
m < \omega \rangle) \\
&\qquad\qquad\iff (\forall q{\le}p) \,\, q \not\forces s \in
(c)_0(\langle \Zn_m \suchthat m < \omega \rangle) \\
&\qquad\qquad\iff (\forall q{\le}p)
\,\, q\restrict C_s \not\forces s \in
(c)_0(\langle \Zn_m \suchthat m < \omega \rangle) \\
&\qquad\qquad\iff (\forall q{\le}p\restrict C_s)
\,\, q\restrict C_s \not\forces s\in
(c)_0(\langle \Zn_m \suchthat m < \omega \rangle)\\
&\qquad\qquad\iff (\forall q{\le}p\restrict C_s) \,\, q \not\forces s
\in (c)_0(\langle \Zn_m \suchthat m < \omega \rangle) \\
&\qquad\qquad\iff p\restrict C_s \forces s \in c(\langle \Zn_m \suchthat
m < \omega \rangle). \endalign$$
So (a)~holds for~$c$.

Let $\alpha_0$ be the ordinal given by~(b) for~$(c)_0$, and let $\alpha
= \alpha_0+1$.  Then (b)~holds for~$c$ for this value of~$\alpha$.  To
see this, let $B$ be the desired set $\{s \in \setfuncs
\omega\kappa\suchthat p \forces s \in c(\langle \Zn_m \suchthat m <
\omega \rangle)\}$. By the inductive hypothesis, $C_s$~is a support for
$(s \in (c)_0(\langle \Zn_m \suchthat m < \omega \rangle))$. Since
conditions in~$P$ are finite and $C_s$~is the increasing union of the
sets~$S_{s \restrict j}$, we have $P \cfrestrict C_s = \bigcup_{j <
\omega} P \cfrestrict S_{s \restrict j}$. Hence, $$\align &p \forces s
\in c(\langle \Zn_m \suchthat m < \omega \rangle) \\ &\qquad\qquad\iff
(\forall q{\le}p)\,\, q \not\forces s \in (c)_0(\langle \Zn_m \suchthat
m < \omega \rangle) \\ &\qquad\qquad\iff (\forall q\text{ compatible
with }p)\,\, q \not\forces s \in (c)_0(\langle \Zn_m \suchthat m <
\omega \rangle) \\ &\qquad\qquad\iff (\forall q\text{ compatible with
}p)\,\, q\restrict C_s \not\forces s \in (c)_0(\langle \Zn_m \suchthat m
< \omega \rangle) \\ &\qquad\qquad\iff (\forall q{\in}P\cfrestrict
C_s\text{ compatible with }p)\,\, q \not\forces s \in (c)_0(\langle
\Zn_m \suchthat m < \omega \rangle) \\ &\qquad\qquad\iff (\forall
j)(\forall q{\in}P\cfrestrict S_{s \restrict j} \text{ compatible with
}p)\,\, q \not\forces s \in (c)_0(\langle \Zn_m \suchthat m < \omega
\rangle). \endalign$$ So $B$ is the intersection over $j$ of the sets
$$\hat B_j = \{s\suchthat (\forall q{\in}P\cfrestrict S_{s \restrict j}
\text{ compatible with }p)\,\, q \not\forces s \in (c)_0(\langle \Zn_m
\suchthat m < \omega \rangle)\}.$$ If we let $$\hat B^{(\tau)} =
\{s\suchthat (\forall q{\in}P\cfrestrict S_\tau \text{ compatible with
}p)\,\, q \not\forces s \in (c)_0(\langle \Zn_m \suchthat m < \omega
\rangle)\},$$ then $\hat B^{(\tau)}$ is a countable intersection (over
$q$) of complements of sets that are $\PI^0_{\alpha_0}$ by the induction
hypothesis, so $\hat B^{(\tau)}$ is $\PI^0_\alpha$. But $\hat
B_j\wkrestrict\tau = \hat B^{(\tau)}\wkrestrict \tau$ for all $\tau \in
\setfuncs j\kappa$, so $\hat B_j$ is $\PI^0_\alpha$. Therefore, $B$ is
$\PI^0_\alpha$, as desired. \QED\enddemo

We can now prove~\thetag{4.3} for $S=\Borel$.  Given $p_0$ and~$\An_n$,
find $c$ and $\langle \Zn_m \suchthat m < \omega \rangle$ as above
for the {\sl complement} of~$\An_n$.  Then we find that
the set~$B_n$ of potential members of~$\An_n$ is just
$\{ s \suchthat p_0 \not\forces s \in c(\langle \Zn_m \suchthat
m < \omega \rangle)\}$.  By Lemma~4.11, the complement
of~$B_n$ is Borel, so $B_n$~is Borel.  This completes the
proof of Theorem~4.1.

Again, more careful accounting of Borel levels shows that,
if $p_0 \forces (\text{$\An$ is $\SIGMA^0_\alpha$})$,
then the set of potential members of~$\An$ is $\SIGMA^0_\alpha$.
Hence, the property $\NNC(\kappa,\lambda,\SIGMA^0_\alpha)$
is preserved by forcing to add Cohen reals
(assuming $\cf\lambda>\omega$).

The proof of Theorem~4.1 does not go through for arbitrary forcing
notions with the countable chain condition. In fact, one can show that
an arbitrary subset of\/~$\setfuncs \omega\kappa$ in the ground model
can be expressed as the set of ``potential members'' of a closed subset
of\/~$\setfuncs \omega\kappa$ in a c.c.c.\ forcing extension. (For $s
\in \setfuncs \omega\kappa$, let~$s^*$ be the set of finite initial
segments of~$s$. Given $B \subseteq \setfuncs \omega\kappa$, let~$P$ be
the poset of partial functions~$p$ from $\setfuncs {<\omega}\kappa$ to
$\{0,1\}$ such that the domain of~$p$ is the union of a finite set and
finitely many sets~$s^*$ for $s \in B$, and $p(\sigma) =0$ for only
finitely many~$\sigma$. If $G\funcfrom \setfuncs {<\omega}\kappa \to
\{0,1\}$ is the resulting generic function and $\An$~is a name for the
closed set $\{s \in \setfuncs \omega\kappa\suchthat (\forall n)\,
G(s\restrict n)=1\}$, then $B$~is the set of potential members
of~$\An$.) So it is still open whether $\NNC(\kappa,\lambda,\Borel)$ is
always preserved by c.c.c.\ forcing.

%

\head 5.  $U$\snug-measurability \endhead

     Throughout this section and the next, the letter~$U$
will denote an ultrafilter, usually over the cardinal~$\kappa$.
We recall several definitions pertaining to ultrafilters: given cardinals
$\kappa$,~$\lambda $, and~$\mu $, an ultrafilter~$U$ over $\kappa$ is
{\it uniform} iff every member of~$U$ has cardinality~$\kappa$; $U$~is
{\it $\lambda$\snug-complete} iff the intersection of any collection
of fewer than~$\lambda $ members of~$U$ is a member of~$U$; $U$~is {\it
$\lambda$\snug-indecomposable} iff every set of cardinality $\lambda $
whose union is in~$U$ has a subset of cardinality less than $\lambda
$ whose union is in~$U$; $U$~is {\it $(\lambda,\mu)$\snug-regular}
iff there is a collection $\{Y_{\alpha}\suchthat  \alpha  < \lambda
\}$ of elements of~$U$ such that, for any $S \subseteq  \lambda $
of cardinality~$\mu $, $\bigcap_{\alpha\in S}Y_{\alpha} = \nullset
$; and $U$~is {\it $(\lambda,\mu)$\snug-nonregular} iff $U$~is not
$(\lambda,\mu)$\snug-regular.

As in the preceding section, given $X \subseteq  \setfuncs
{\omega}\kappa$ and $\sigma  \in \setfuncs {<\omega}\kappa$, define
$X\wkrestrict \sigma$ to be $\{s \in\nobreak\setfuncs
{\omega}\kappa\suchthat  \sigma \concat s \in  X\}$. Again recall that
any closed subset of\/~$\setfuncs \omega\kappa$ can be expressed in the
form~$[T]$, the set of infinite branches through some tree $T \subseteq
\setfuncs {<\omega}\kappa$.

\definition{Definition 5.1} Let $U$ be an ultrafilter over $\kappa$.

(a)  A tree $T \subseteq  \setfuncs {<\omega}\kappa$ is {\it
$U$\snug-branching} iff $\nullseq \in  T$ and, for any $\sigma  \in  T$,
$\{\alpha  \in  \kappa\suchthat  \sigma \concat \langle \alpha \rangle
\in  T\} \in  U$.

(b)  A set $X \subseteq  \setfuncs {\omega}\kappa$ is {\it $U$\snug-large}
({\it $U$\snug-small}) iff there is a
$U$\snug-branching tree $T$ such that $[T] \subseteq  X$ ($[T]\cap
X = \nullset)$.

(c)  A set $X \subseteq  \setfuncs {\omega}\kappa$ is {\it
$U$\snug-determined\/} iff $X$~is either $U$\snug-large or $U$\snug-small.

(d)  A set $X \subseteq  \setfuncs {\omega}\kappa$ is {\it
$U$\snug-null\/} iff, for each $\sigma  \in  \setfuncs {<\omega}\kappa$,
$X\wkrestrict \sigma $ is $U$\snug-small.

(e)  A set $X \subseteq  \setfuncs {\omega}\kappa$ is {\it
$U$\snug-measurable} iff, for each $\sigma  \in  \setfuncs
{<\omega}\kappa$, $X\wkrestrict \sigma $ is
$U$\snug-determined.\enddefinition

     The remainder of this section is devoted to results on
$U$\snug-null and $U$\snug-measurable sets; many of these results are
analogous to facts about the standard notion of measurability for
subsets of, say, the Cantor space, or the real line.  In the next section we
will use these results to obtain further information about the
property~$\NNC$.

Louveau~\cite{\Louveau} gives definitions equivalent to these, for the
case $\kappa = \omega$, and uses them to give an alternate proof of
the theorem of Silver~\cite{\Silver} that all analytic subsets
of\/~$\setfuncs \omega\omega$ are Ramsey.  Much of the rest of this section
appears in another form in Louveau's paper.  (Carlson and Galvin have
also done unpublished work along these lines.)

\proclaim{Proposition 5.2} The intersection of two $U$\snug-branching
trees is $U$\snug-branching.  Hence, the $U$\snug-large sets form a
filter over\/~$\setfuncs {\omega}\kappa$, and the $U$\snug-small sets
form the dual ideal. \endproclaim

\demo{Proof} Easy. \QED\enddemo

\proclaim{Lemma 5.3} If $X_n \subseteq  \setfuncs {\omega}\kappa$ for $n <
\omega $ and, for each $n < \omega $ and $\sigma  \in  \setfuncs n\kappa$,
$X_n\wkrestrict \sigma $ is $U$\snug-small, then $\bigcup_{n<\omega}X_n$
is $U$\snug-small. \endproclaim

\demo{Proof} Let $X = \bigcup_{n<\omega}X_n$. For each $n \in  \omega $
and each $\sigma  \in  \setfuncs n\kappa$, choose a $U$\snug-branching
tree~$T'_\sigma$ such that $(X_n\wkrestrict \sigma)\cap [T'_\sigma] =
\nullset $.  Let $T_n = \setfuncs {<n}\kappa\cup \{\sigma \concat \tau
\suchthat  \sigma  \in  \setfuncs n\kappa,\, \tau  \in T'_\sigma\}$; then
$T_n$ is a $U$\snug-branching tree, $\setfuncs m\kappa \subseteq  T_n$
for $m \le  n$, and $[T_n]\cap X_n = \nullset $.  Let $T =
\bigcap_{n<\omega}T_n$; then $[T]\cap X = \nullset $.  Clearly $\nullseq
\in  T$, and for any $\sigma  \in  T$, if $\sigma  \in  \setfuncs
n\kappa$, then
     $$\align\{\alpha  \in  \kappa\suchthat  \sigma \concat \langle
     \alpha \rangle \in  T\} &= \{\alpha  \in  \kappa\suchthat  (\forall
     m{\in} \omega)\, \sigma \concat \langle \alpha \rangle \in  T_m\}\\
     &= \{\alpha  \in  \kappa\suchthat  (\forall m{\le} n) \, \sigma
     \concat \langle \alpha \rangle \in  T_m\}\\
     &\qquad\qquad\qquad\text{(since $\setfuncs {n+1}\kappa \subseteq
     T_m$ for $m > n$)}\\ &= \bigcap_{m\le n}\{\alpha  \in
     \kappa\suchthat \sigma \concat \langle \alpha \rangle \in
     T_m\},\endalign$$
and since each tree $T_m$ is $U$\snug-branching, each set $\{\alpha
\in\nobreak\kappa\suchthat  \sigma \concat \langle \alpha \rangle \in
T_m\}$ is in $U$, so $\{\alpha \in\nobreak\kappa\suchthat  \sigma
\concat \langle \alpha \rangle \in  T\} \in  U$.  Therefore, $T$ is
$U$\snug-branching, so $X$ is $U$\snug-small. \QED\enddemo

\proclaim{Theorem 5.4} The $U$\snug-null sets form a $\sigma $\snug-ideal.
\endproclaim

\demo{Proof} Clearly any subset of a $U$\snug-null set is $U$\snug-null.
Now suppose that we have $U$\snug-null sets $X_n$, $n \in  \omega $, and
let $X = \bigcup_{n<\omega}X_n$; we must see that $X$ is $U$\snug-null.
For each $\sigma  \in  \setfuncs {<\omega}\kappa$, we have $X\wkrestrict
\sigma = \bigcup_{n<\omega}(X_n\wkrestrict \sigma)$.  If $n \in  \omega
$ and $\tau  \in  \setfuncs n\kappa$, then $(X_n\wkrestrict
\sigma)\wkrestrict \tau  = X_n\wkrestrict (\sigma \concat \tau)$ is
$U$\snug-small by hypothesis; therefore, by Lemma 5.3, $X\wkrestrict
\sigma $ is $U$\snug-small. Since $\sigma $ was arbitrary, $X$ is
$U$\snug-null. \QED\enddemo

\proclaim{Theorem 5.5} Every open subset of\/ $\setfuncs {\omega}\kappa$ is
$U$\snug-measurable. \endproclaim

\demo{Proof} Let $X \subseteq  \setfuncs {\omega}\kappa$ be open; then
$X\wkrestrict \sigma $ is open for each $\sigma  \in  \setfuncs
{<\omega}\kappa$, so it will suffice to show that $X$ is
$U$\snug-determined. Let $S = \{\sigma  \in  \setfuncs
{<\omega}\kappa\suchthat \text{$X\wkrestrict \sigma$  is
$U$\snug-large}\}$. If $\nullseq \in\nobreak S$, we are done, so assume
$\nullseq \notin  S$.  If $\sigma  \in  \setfuncs {<\omega}\kappa$ and
$A = \{\alpha  \in  \kappa\suchthat  \sigma \concat \langle \alpha
\rangle \in  S\} \in  U$, then choose a $U$\snug-branching tree
$T_{\alpha}$ for each $\alpha  \in  A$ such that $[T_{\alpha}] \subseteq
X\wkrestrict (\sigma \concat \langle \alpha \rangle)$, and let $$T' =
\{\nullseq\}\cup \{\langle \alpha \rangle\concat \tau \suchthat  \alpha
\in  A,\, \tau  \in T_{\alpha}\};$$ it is easy to see that $T'$ is a
$U$\snug-branching tree and $[T'] \subseteq  X\wkrestrict \sigma $, so
$\sigma  \in  S$. Hence, for any $\sigma  \in  \setfuncs
{<\omega}\kappa\setdiff S$, $\{\alpha  \in \kappa\suchthat  \sigma
\concat \langle \alpha \rangle \notin  S\} \in  U$. Now let $$T =
\{\sigma  \in  \setfuncs {<\omega}\kappa\suchthat  (\forall m{\le} \ell
(\sigma))\, \sigma \restrict m \notin  S\};$$ since $\sigma  \in  T$ and
$\sigma \concat \langle \alpha \rangle \notin  S$ imply $\sigma \concat
\langle \alpha \rangle \in  T$, $T$ is a $U$\snug-branching tree.  If $s
\in  [T]$, then $s\restrict n \notin  S$ for all $n \in  \omega $, so
$X\wkrestrict (s\restrict n) \ne  \setfuncs {\omega}\kappa$ for all $n
\in  \omega $; since $X$ is open, this implies $s \notin  X$.
Therefore, $[T]\cap X = \nullset $, so $X$ is $U$\snug-small.
\QED\enddemo

\proclaim{Theorem 5.6} The $U$\snug-measurable sets form a $\sigma
$\snug-algebra of subsets of\/ $\setfuncs {\omega}\kappa$. \endproclaim

\demo{Proof} Clearly the complement of a $U$\snug-measurable subset
of\/ $\setfuncs {\omega}\kappa$ is $U$\snug-measurable.  Now suppose we have
$U$\snug-measurable sets $X_n$, $n \in  \omega $; we must see that $X =
\bigcup_{n<\omega}X_n$ is $U$\snug-measurable.  Again it will suffice to
show that $X$ is $U$\snug-determined, since the same will
apply to $X\wkrestrict \sigma  = \bigcup_{n<\omega}(X_n\wkrestrict \sigma)$
for any $\sigma  \in  \setfuncs {<\omega}\kappa$.  Let
     $$Y = \{s \in  \setfuncs {\omega}\kappa\suchthat  (\exists
     n{,}m{\in} \omega)\, X_n\wkrestrict (s\restrict m)\text{ is
     $U$\snug-large}\}.$$
Then $Y$ is open, so by Theorem 5.5 there is a $U$\snug-branching
tree $T$ such that either $[T] \subseteq  Y$ or $[T]\cap Y = \nullset $.

Suppose $[T] \subseteq  Y$, and let
$$
S = \{\sigma  \in  T\suchthat  (\exists n{\in} \omega)(X_n\wkrestrict \sigma
\text{ is $U$\snug-large})\text{ and }
(\forall n{\in} \omega)(\forall m{<}\ell (\sigma))(X_n\wkrestrict (\sigma
\restrict m)\text{ is $U$\snug-small})\};
$$
then $S$ is an antichain in $T$, and since $[T] \subseteq  Y$, for
each $s \in  [T]$ there is $m \in  \omega $ such that $s\restrict m
\in  S$.  For each $\sigma  \in  S$, choose a $U$\snug-branching tree
$T_{\sigma}$ such that $[T_{\sigma}] \subseteq  X_n\wkrestrict \sigma $
for some~$n$.  Let
     $$T' = \{\sigma \concat \tau \suchthat  \sigma  \in  S,\, \tau  \in
     T_{\sigma}\}\cup \{\sigma \restrict m\suchthat  \sigma  \in  S,\, m
     < \ell (\sigma)\}.$$
Clearly $\{\sigma \restrict m\suchthat  \sigma \in S,\, m < \ell (\sigma)\}$
is a subtree of\/ $T$ which does not meet $S$; since every infinite branch
through $T$ meets $S$, $\{\sigma \restrict m\suchthat \sigma \in  S,\, m <
\ell (\sigma)\}$ has no infinite branches.  Therefore, for any $s \in
[T']$, there are $\sigma  \in  S$ and $b \in  [T_{\sigma}]$ such that
$s = \sigma \concat b$; since $[T_{\sigma}] \subseteq  X\wkrestrict \sigma
$ for each $\sigma  \in  S$, $[T'] \subseteq  X$.  Now, if $\sigma  \in
S$ and $\tau  \in  T_{\sigma}$, then
     $$\{\alpha  \in  \kappa\suchthat  \sigma \concat \tau \concat \langle
     \alpha \rangle \in  T'\} = \{\alpha  \in  \kappa\suchthat  \tau \concat
     \langle \alpha \rangle \in  T_{\sigma}\} \in  U;$$
and if $\sigma  \in  S$ and $m < \ell (\sigma)$, then
     $$\{\alpha  \in  \kappa\suchthat  (\sigma \restrict m)\concat \langle
     \alpha \rangle \in  T'\} = \{\alpha  \in  \kappa\suchthat  (\sigma
     \restrict m)\concat \langle \alpha \rangle \in  T\} \in  U.$$
(To see this, let $\tau = (\sigma \restrict m)\concat \langle
\alpha \rangle$.  Since $S$~is an antichain, no proper initial segment
of~$\sigma$ is in~$S$, so no proper initial segment
of~$\tau$ is in~$S$.  This means that,
if $\tau \in T'$, then $\tau$ must be a
member of~$S$ or an initial segment of one, so $\tau \in T$.
On the other hand, if $\tau \in T$, then $\tau$~can be extended to some
$s \in [T]$, and some initial segment of~$s$ must be in~$S$, so
$\tau$~must be a member of~$S$ or an initial segment of one,
so $\tau\in T'$.)
Therefore, $T'$ is $U$\snug-branching, so $X$ is
$U$\snug-large.

Now suppose $[T]\cap Y = \nullset $.  Let $X'_n = X_n\cap [T]$.
If $\sigma  \in  T$,
then $X_n\wkrestrict \sigma $ is $U$\snug-small by definition of
$Y$, so $X'_n\wkrestrict \sigma $ is $U$\snug-small; if $\sigma
\in  \setfuncs {<\omega}\kappa\setdiff T$, then $X'_n\wkrestrict \sigma  =
\nullset $.  Therefore, $X'_n$ is $U$\snug-null, so by Theorem 5.6
$X\cap [T] = \bigcup_{n<\omega}X'_n$ is $U$\snug-null.  By Proposition
5.2, $X \subseteq  (X\cap [T])\cup (\setfuncs {\omega}\kappa\setdiff [T])$ is
$U$\snug-small. \QED\enddemo

\proclaim{Corollary 5.7} Every Borel subset of\/ $\setfuncs
{\omega}\kappa$ is $U$\snug-measurable. \QED\endproclaim

\proclaim{Lemma 5.8} For every set $X \subseteq  \setfuncs
{\omega}\kappa$ there is an $F_{\sigma}$ set $Z \subseteq  X$ such that,
for each $\sigma \in  \setfuncs {<\omega}\kappa$, if $X\wkrestrict
\sigma $ is $U$\snug-large, then $Z\wkrestrict \sigma $ is
$U$\snug-large. \endproclaim

\demo{Proof} For each $\sigma $ such that $X\wkrestrict \sigma $ is
$U$\snug-large, choose a $U$\snug-branching tree $T_{\sigma}$ such that
$[T_{\sigma}] \subseteq  X\wkrestrict \sigma $; if $X\wkrestrict\sigma$
is not $U$\snug-large, let $T_\sigma = \nullset$.  Now let $Z = \{\sigma
\concat s\suchthat  \sigma  \in  \setfuncs {<\omega}\kappa,\, s \in
[T_{\sigma}]\}$.  Clearly $Z \subseteq  X$ and, for all $\sigma $, if
$X\wkrestrict \sigma $ is $U$\snug-large, then $Z\wkrestrict \sigma $ is
$U$\snug-large.  And since the sets $Z_n = \{\sigma \concat s\suchthat
\sigma  \in  \setfuncs n\kappa,\, s \in  [T_{\sigma}]\}$ for $n \in
\omega $ are closed, and $Z = \bigcup_{n<\omega}Z_n$, $Z$ is
$F_{\sigma}$. \QED\enddemo

\proclaim{Theorem 5.9} A set $X \subseteq  \setfuncs {\omega}\kappa$ is
$U$\snug-measurable iff there are sets $Z,Y$ such that $Z \subseteq
X \subseteq\nobreak Y$, $Z$ is $F_{\sigma}$, $Y$ is $G_\delta$, and
$Y\setdiff Z$
is $U$\snug-null. \endproclaim

\demo{Proof} If $Z$ and $Y$ are as above, then $X = Z\cup (X\setdiff
Z)$, $Z$ is Borel and hence $U$\snug-measurable, and $X\setdiff Z$ is
$U$\snug-null, so $X$ is $U$\snug-measurable.  Conversely, if $X$ is
$U$\snug-measurable, then we can find $F_{\sigma}$ sets $Z \subseteq
X$, $Z' \subseteq  \setfuncs {\omega}\kappa\setdiff X$ as in Lemma 5.8.
Let $Y = \setfuncs {\omega}\kappa\setdiff Z'$.  For each $\sigma  \in
\setfuncs {<\omega}\kappa$, either $X\wkrestrict \sigma $ or $(\setfuncs
{\omega}\kappa\setdiff X)\wkrestrict \sigma $ is $U$\snug-large, so
either $Z\wkrestrict \sigma $ or $Z'\wkrestrict \sigma $ is
$U$\snug-large, so $(Y\setdiff Z)\wkrestrict \sigma $ is $U$\snug-small.
Therefore, $Y\setdiff Z$ is $U$\snug-null. \QED\enddemo

\proclaim{Theorem 5.10} The collection of $U$\snug-measurable subsets
of\/ $\setfuncs {\omega}\kappa$ is closed under Suslin's operation
$\Suslin$. \endproclaim

\demo{Proof} The collection of $U$\snug-null sets is a $\sigma
$\snug-ideal over $\setfuncs {\omega}\kappa$.  The collection of
$U$\snug-measurable sets is a $\sigma $\snug-algebra.  For every set
$X \subseteq
\setfuncs {\omega}\kappa$, there is a $U$\snug-measurable set $Y \supseteq
X$ such that any $U$\snug-measurable subset of $Y\setdiff X$ is
$U$\snug-null.  (Let $Y$ be the complement of the set $Z \subseteq
\setfuncs {\omega}\kappa\setdiff X$ obtained by applying Lemma 5.8 to
$\setfuncs {\omega}\kappa\setdiff X$.)  By Theorem~2H.1 of
Moschovakis~\cite{\Moschovakis}, these statements imply
the desired result.  (Theorem~2H.1 is stated only for
certain spaces~$\scriptX$, but the proof of the relevant part
applies to any set~$\scriptX$.)
\QED\enddemo

     So, in the case $\kappa = \omega $, we see that all analytic and
coanalytic sets, and many others, are $U$\snug-measurable.  This does
not, however, necessarily extend to all $\DELTA^1_2$ (or even
$\Delta^1_2$) subsets of\/ $\setfuncs {\omega}\omega $ (unless $U$ is
principal, in which case every subset of\/ $\setfuncs {\omega}\omega $
is $U$\snug-measurable). If $U$~is nonprincipal, then clearly $[T]$~is a
perfect set for any $U$\snug-branching tree $T$. Therefore, any subset
of\/~$\setfuncs\omega\omega$ which is a Bernstein set (a set such that
neither it nor its complement has a perfect subset; such sets can be
constructed using a well-ordering of\/~$\setfuncs\omega\omega$) cannot be
$U$\snug-determined. Well-known results in descriptive set theory
\cite{\Jech,~\S41} show that, in the constructible universe, one can
construct a $\Delta^1_2$~Bernstein set.

Louveau's proof that all analytic sets are Ramsey is completed by
the following result.

\proclaim{Proposition 5.11} If\/~$U$ is a nonprincipal ultrafilter
over~$\omega$, and $X \subseteq \setfuncs \omega\omega$ is
$U$\snug-determined, then there an infinite set $H \subseteq \omega$
such that either all strictly increasing $\omega$\snug-sequences
from~$H$ are in~$X$ or all such sequences are in the complement of~$X$.
\endproclaim

\demo{Proof}
Let $T$ be a $U$\snug-branching tree such that $[T] \subseteq X$ or
$[T] \cap X = \nullset$.  It will suffice to construct an
infinite set~$H$ such that all strictly increasing sequences from~$H$
are in~$[T]$.  To do this, we will recursively choose
natural numbers $h_0 < h_1 < h_2 < \dotso$ such that
every finite subsequence of $\langle h_0,h_1,h_2,\dotsc\rangle$ is in~$T$.

Suppose we have $h_i$ for $i < n$.  For each subsequence~$\sigma$
of $\langle h_0,h_1,\dots,h_{n-1}\rangle$, since $\sigma \in T$
and $T$~is $U$\snug-branching, the set of~$k$ such that
$\sigma\concat\langle k\rangle \in T$ is in~$U$.  There are
$2^n$~such subsequences~$\sigma$; the intersection of the
$2^n$~corresponding sets in~$U$ is still in~$U$.  Therefore,
we can choose~$h_n$ to be any member of this intersection which
(if $n>0$) is above~$h_{n-1}$; then every subsequence of
$\langle h_0,h_1,\dots,h_n\rangle$ will be in~$T$, as desired.
\QED\enddemo

     The strong analogy between $U$\snug-measurability and ordinary
measurability suggests the following question: is there a $\sigma
$\snug-additive probability measure $m$ (on some $\sigma $\snug-algebra
of subsets of\/ $\setfuncs {\omega}\kappa$) such that all
$U$\snug-measurable sets are $m$\snug-measurable? The answer is yes if
$U$ is $\aleph_1$\snug-complete, because we can define such an $m$ by
letting $m(X) = 1$ for all $U$\snug-large sets~$X$ and $m(X) = 0$ for
all $U$\snug-small sets~$X$.  On the other hand, if $U$ is not
$\aleph_1$\snug-complete, then such an $m$ cannot exist unless its
completion is a measure on all subsets of\/ $\setfuncs {\omega}\kappa$;
this follows from the following proposition.

\proclaim{Proposition 5.12} Let $U$ be an $\aleph_1$\snug-incomplete
ultrafilter over $\kappa$.  Suppose that $m$ is a ($\sigma
$\snug-additive) probability measure on a $\sigma $\snug-algebra of subsets
of\/~$\setfuncs {\omega}\kappa$ which includes all clopen subsets
of\/~$\setfuncs {\omega}\kappa$.  Then there is a
$U$\snug-null set $X \subseteq  \setfuncs {\omega}\kappa$ such that $m(X) =
1$. \endproclaim

\demo{Proof} By Theorem 5.4, it suffices to prove that, for any $\eps >
0$, there is a $U$\snug-null set $X$ such that $m(X) \ge  1-\eps $. Let
$\{S_n\suchthat  n < \omega \}$ be a collection of sets not in $U$ such
that $\bigcup_{n<\omega}S_n = \kappa$.  We will define sets $X_n
\subseteq  \setfuncs {\omega}\kappa$ with $m(X_n) > 1-\eps $ by
recursion on $n$. Let $X_0 = \setfuncs {\omega}\kappa$.  Given $X_n$
such that $m(X_n) > 1-\eps $, let $Y_i = X_n\cap \{s\suchthat  s(n) \in
S_{<i}\}$ for $i < \omega $. Then $\langle Y_i\suchthat  i < \omega
\rangle$ is an increasing sequence of $m$\snug-measurable sets and
$\bigcup_{i<\omega}Y_i = X_n$, so there must be an $i < \omega $ such
that $m(Y_i) > 1-\eps $; let $X_{n+1}$ be $Y_i$ for the least such $i$.
This completes the definition of $\langle X_n\suchthat  n < \omega
\rangle$; it is clear that this is a decreasing sequence of sets such
that $m(X_n) > 1-\eps $ for all $n$, but $X_n\wkrestrict \sigma $ is
$U$\snug-small for all $\sigma  \in \setfuncs {<n}\kappa$.  Therefore,
if we let $X = \bigcap_{n<\omega}X_n$, then $X$ will be a $U$\snug-null
set such that $m(X) \ge  1-\eps $, as desired. \QED\enddemo

\head 6.  $(\lambda,\mu;\alpha)$\snug-nonregularity \endhead

     In this section we will apply the results of the
previous section to obtain new information about the property $\NNC$.
In particular, we will define a property of ultrafilters $U$ over $\kappa$
which implies $\NNC(\kappa,\mu,\text{$U$\snug-measurable})$, and then
give several cases in which this property is satisfied.  The definition
of this property, $(\lambda,\mu;\omega)$\snug-nonregularity, will be
given in somewhat more generality than necessary, in order to show
its relation to the usual definition of $(\lambda,\mu)$\snug-nonregularity.
We start by generalizing Definition 5.1(a).

\definition{Definition 6.1} A tree $T \subseteq  \setfuncs
{<\alpha}\kappa$ is a {\it closed $U$\snug-branching tree} of height
$\alpha $ iff:

(a)  for every successor $\beta +1 < \alpha $ and every $s \in  T\cap
\setfuncs \beta \kappa$,
     $\{\eta \in  \kappa\suchthat  s\concat \langle \eta \rangle \in
     T\} \in  U$;

(b)  for every non-successor $\beta  < \alpha $ and every $s \in
\setfuncs \beta \kappa$, $s \in  T$ iff $s\restrict \gamma  \in  T$ for
all $\gamma  < \beta $.\enddefinition

     Note that a tree $T \subseteq  \setfuncs {<\omega}\kappa$ is a
     $U$\snug-branching tree under
Definition 5.1(a) iff it is a closed $U$\snug-branching tree of height
$\omega $.

\definition{Definition 6.2} If $\kappa$, $\lambda $, and $\mu $ are
cardinals and $\alpha $ is an ordinal, then an ultrafilter $U$
over~$\kappa$ is {\it $(\lambda,\mu;\alpha)$\snug-regular} iff there is a
family $\{T_\beta\suchthat  \beta  \in  \lambda \}$ of closed
$U$\snug-branching trees of height $\alpha +1$ such that no subfamily
$\{T_\beta\suchthat  \beta  \in  S\}$ with $S \subseteq  \lambda $ of
cardinality $\mu $ has a common maximal element, i.e. a sequence $s \in
\setfuncs \alpha \kappa$ such that $s \in  T_\beta$ for all $\beta  \in
S$.  The ultrafilter $U$ is {\it $(\lambda,\mu;\alpha)$\snug-nonregular}
iff it is not $(\lambda,\mu;\alpha)$\snug-regular.\enddefinition

We start with some easy but useful results.

\proclaim{Proposition 6.3}

{\rm (a)}  If $U$ is a $(\lambda,\mu;\alpha)$\snug-regular ultrafilter,
$\lambda' \le  \lambda $, $\mu ' \ge  \mu $, and $\alpha ' \ge  \alpha
$, then $U$ is $(\lambda ',\mu ';\alpha ')$\snug-regular.

{\rm (b)}  If $U$ is $(\lambda,\lambda ';\alpha)$\snug-nonregular and
$(\lambda',\mu;\alpha')$\snug-nonregular,
then $U$ is $(\lambda,\mu;\alpha +\alpha')$\snug-nonreg\-ular.

{\rm (c)}  If $U$ is $\mu^+$\snug-complete, then the intersection of $\mu $
closed $U$\snug-branching trees of the same height is a closed
$U$\snug-branching tree of that height.

{\rm (d)}  An ultrafilter $U$ is $(\lambda,\mu;\omega)$\snug-regular iff
there is a family $\{X_\beta\suchthat  \beta  < \lambda \}$ of
$U$\snug-large sets such that any subfamily $\{X_\beta\suchthat  \beta
\in  S\}$ with $S \subseteq  \lambda $ of cardinality~$\mu $ has empty
intersection. \endproclaim

\demo{Proof} Suppose $U$ is an ultrafilter over $\kappa$.  For (a),
let $\{T_\beta\suchthat  \beta  < \lambda \}$ be a witness to the
$(\lambda,\mu;\alpha)$\snug-regularity of $U$, and let
     $$T'_\beta  = T_\beta\cup \{s \in  \setfuncs
     {\le\alpha'}\kappa\suchthat  \ell (s) > \alpha \text{ and }
     s\restrict \alpha  \in  T_\beta\};$$
then $\{T'_\beta \suchthat  \beta  < \lambda '\}$ witnesses the
$(\lambda ',\mu ';\alpha ')$\snug-regularity of $U$.  For (b), let
$\{T_\beta\suchthat \beta  < \lambda \}$ be a collection of closed
$U$\snug-branching trees of height $\alpha +\alpha' +1$.  Then
$\{T_\beta\cap \setfuncs {\le\alpha}\kappa\suchthat \beta  < \lambda \}$
is a collection of closed $U$\snug-branching trees of height $\alpha
+1$, so there is a set $S' \subseteq \lambda $ of cardinality~$\lambda
'$ such that $\{T_\beta\cap \setfuncs {\le\alpha}\kappa\suchthat  \beta
< \lambda \}$ has a common maximal element $s$.  Let $T'_\beta  =
\{t\suchthat  s\concat t \in  T_\beta\}$ for $\beta \in  S'$;
$\{T'_\beta \suchthat  \beta  \in  S'\}$ is a collection of closed
$U$\snug-branching trees of height $\alpha' +1$, and $S'$~has
cardinality $\lambda '$, so there is a set $S \subseteq  S'$ of
cardinality $\mu $ such that $\{T'_\beta \suchthat  \beta  \in  S\}$ has
a common maximal element $s'$.  Then $s\concat s'$ is a common maximal
branch of $\{T_\beta\suchthat  \beta  \in  S\}$.

     For part (c), note that if $s$ is in the intersection and is not
maximal in the original trees, then the set of $\beta $ such that
$s\concat \langle \beta \rangle$ is in the intersection is the
intersection of $\mu $ members of $U$.  Finally, for part (d), note that
$F(T) = T\cup [T]$ and $F^{-1}(T) = T\cap \setfuncs {<\omega}\kappa$
define a one-to-one correspondence between the closed $U$\snug-branching
trees of height $\omega $ and the closed $U$\snug-branching trees of
height $\omega +1$.\QED\enddemo

Trivially, any ultrafilter is $(\lambda,\mu;\alpha)$\snug-regular if
$\mu  > \lambda $.  On the other hand, if $\mu$~is finite, then every
ultrafilter~$U$ is $\mu^+$\snug-complete; hence, by Proposition~6.3(c),
$U$~is $(\lambda,\mu;\alpha)$\snug-nonregular for any $\lambda \ge \mu$
and any~$\alpha$.

     The next result gives the motivation for the
term `$(\lambda,\mu;\alpha)$\snug-regular,' and shows that
$(\lambda,\mu;\alpha)$\snug-nonregularity gives a family of properties
between $(\lambda,\mu)$\snug-nonregularity and $\mu^+$\snug-completeness.

\proclaim{Proposition 6.4} Let $U$ be an ultrafilter.

{\rm (a)}  If $\lambda  \ge  \mu $, then $U$~is
$(\lambda,\mu;0)$\snug-nonregular.

{\rm (b)}  $U$ is $(\lambda,\mu;1)$\snug-regular iff $U$~is
$(\lambda,\mu)$\snug-regular.

{\rm (c)}  If\/ $|\alpha | \ge  \lambda  \ge  \mu  \ge  \aleph_0$, then
$U$~is $(\lambda,\mu;\alpha)$\snug-nonregular iff $U$~is
$\mu^+$\snug-complete. \endproclaim

\demo{Proof} Parts (a) and (b) are easy.  For part (c), first suppose
that $U$ is $\mu^+$\snug-complete; then Proposition 6.3(c) easily
implies that $U$ is $(\lambda,\mu;\alpha)$\snug-nonregular.  To prove
the other direction of part (c), we first need a lemma.

\proclaim{Lemma 6.5} For any infinite cardinals $\lambda  \ge  \mu
$ such that $\mu $ is regular, there is a function $F\funcfrom  \lambda
\times \lambda \to \mu $ such that, for any $S \subseteq  \lambda $
of cardinality $\mu $, there is $\alpha  \in  \lambda $ such that
$\{F(\alpha,\beta)\suchthat  \beta  \in  S\}$ has cardinality $\mu
$. \endproclaim

\demo{Proof} We prove this for all {\sl ordinals} $\lambda  \ge  \mu
$, by induction on $\lambda $.  If $\lambda  = \mu $, we simply
let $F(\alpha,\beta) = \beta $.  If $\lambda $ is not a cardinal,
then $\lambda  > |\lambda | \ge  \mu $, so let $f\funcfrom  \lambda \to
|\lambda |$ be a bijection and let $F'\funcfrom  |\lambda |\times |\lambda
|\to \mu $ be obtained from the induction hypothesis; then the function
$F\funcfrom  \lambda \times \lambda \to \mu $ defined by $F(\alpha,\beta)
= F'(f(\alpha),f(\beta))$ has the required properties.  Now suppose
that $\lambda $ is a cardinal greater than $\mu $.  Let $\langle
\lambda_{\alpha}\suchthat  \alpha  < \cf \lambda \rangle$ be a strictly
increasing sequence of ordinals with limit $\lambda $, such that
$\lambda_0 \ge  \mu $.  For each $\alpha  < \cf \lambda $, obtain a
function $F_{\alpha}\funcfrom  \lambda_{\alpha}\times \lambda_{\alpha}\to
\mu $ from the induction hypothesis.  Let $\delta_{\alpha} =
1+\sum_{\beta<\alpha}\lambda_\beta$ for $\alpha  < \cf \lambda $; since
$\lambda $ is a cardinal, it is easy to see that $\delta_{\alpha} <
\lambda $.  Now define a function $F\funcfrom  \lambda \times \lambda \to
\mu $ by:
     $$\alignat2 F(\delta_{\alpha}+\gamma,\beta) &=
     F_{\alpha}(\gamma,\beta)\qquad&&\text{if $\gamma, \beta  <
     \lambda_{\alpha}$;}\\
     F(0,\delta_{\alpha}+\gamma) &= \alpha \qquad&&\text{if $\cf \lambda
     = \mu $ and $\gamma  < \lambda_{\alpha}$;}\\
     F(\gamma,\beta)    &= 0\qquad&&\text{for all other
     $(\gamma,\beta)$.}\endalignat$$
To see that this works, let $S$ be any subset of $\lambda $ of
cardinality $\mu $.  If there is $\beta  < \lambda $ such that
$|S\cap \beta | = \mu $, then there is $\alpha  < \cf \lambda $
such that $|S\cap \lambda_{\alpha}| = \mu $; if we choose $\gamma
< \lambda_{\alpha}$ such that $|\{F_{\alpha}(\gamma,\beta)\suchthat
\beta  \in  S\cap \lambda_{\alpha}\}| = \mu $, then we will have
$|\{F(\delta_{\alpha}+\gamma,\beta)\suchthat  \beta  \in  S\}| = \mu $.
Now suppose $|S\cap \beta | < \mu $ for all $\beta  < \lambda $.  This
clearly implies $\cf \lambda  = \cf \mu $, so since $\mu $ is regular
we get $|\{F(0,\beta)\suchthat  \beta  \in  S\}| = \mu $.  This completes
the induction.\QED\enddemo

     To finish the proof of Proposition 6.4, suppose that $|\alpha | \ge
\lambda  \ge  \mu  \ge  \aleph_0$ and $U$ is an ultrafilter over
$\kappa$ which is not $\mu^+$\snug-complete; we must see that $U$ is
$(\lambda,\mu;\alpha)$\snug-regular.  Let $\mu '$ be the least cardinal
such that $U$ is not ${\mu '}^+$\snug-complete; then $\mu  \ge  \mu '
\ge  \aleph_0$, so by Proposition 6.3(a) we may assume that $\mu  = \mu
'$.  It is well-known \cite{\Jech, \S27, p.~299} that $\mu '$ must be
either $\aleph_0$ or a measurable cardinal, so $\mu $ is regular.  Let
$\{W_\beta\suchthat  \beta  < \mu \}$ be a family of sets in $U$ which
has empty intersection, and let $Z_\beta =
\bigcap_{\gamma<\beta}W_\gamma$.  By our assumption, $Z_\beta \in  U$
for all $\beta  < \mu $, but, for any $S \subseteq  \mu $ of cardinality
$\mu $, $\bigcap \{Z_\beta\suchthat  \beta  \in  S\} = \nullset $.
Choose $F\funcfrom  \lambda \times \lambda \to \mu $ satisfying the
conclusion of Lemma 6.5.  Now, for each $\beta  < \lambda $, define a
closed $U$\snug-branching tree $T_\beta$ of height $\alpha +1$ as
follows: for any $s \in  \setfuncs {\le\alpha}\kappa$, $s \in  T_\beta$
iff, for each $\gamma < \min(\ell (s),\lambda)$, $s(\gamma) \in
Z_{F(\gamma,\beta)}$.  Let $S$ be any subset of $\lambda $ of
cardinality $\mu $; we must see that $\{T_\beta\suchthat  \beta  \in
S\}$ has no common maximal element. The choice of $F$ guarantees that
there is a $\gamma  < \lambda $ such that $\{F(\gamma,\beta)\suchthat
\beta  \in  S\}$ has cardinality $\mu $, and hence $\bigcap
\{Z_{F(\gamma,\beta)}\suchthat  \beta  \in  S\} = \nullset $.  But if
$s$~were a common maximal branch of $\{T_\beta\suchthat \beta  \in
S\}$, we would have $s(\gamma) \in  Z_{F(\gamma,\beta)}$ for each $\beta
\in  S$, which is impossible.  Therefore, $\{T_\beta\suchthat \beta  \in
S\}$ has no common maximal element; since $S$~was arbitrary, $U$ is
$(\lambda,\mu;\alpha)$\snug-regular.\QED\enddemo

     We now give the reason for studying
$(\lambda,\mu;\alpha)$\snug-nonregularity here.

\proclaim{Proposition 6.6} Let $U$~be an ultrafilter over~$\kappa$
such that every set in~$U$ has cardinality at least~$\lambda$.
If $U$ is $(\lambda,\mu;\omega)$\snug-nonregular, then
$\NNC(\kappa,\mu,\Umble)$.
If $U$ is $(\lambda,\mu';\omega)$\snug-nonregular for all\/ $\mu'<\mu$, then
$\NNC(\kappa,{<}\mu,\Umble)$.
\endproclaim

\demo{Proof} We prove the second implication; the proof of the first is
the same (or one can easily deduce the first from the second).  Let
$U$~be $(\lambda,\mu';\omega)$\snug-nonregular for all $\mu'<\mu$.
Suppose $\{A_n\suchthat  n < \omega \}$ is a family of
$U$\snug-measurable sets with union $\setfuncs {\omega}\kappa$; we must
show that there is an~$n$ such that, for all $\mu' < \mu$, $A_n$ is not
$\mu'$\snug-narrow in the $n$\snug'th coordinate.  By Lemma~5.3, there
exist $n < \omega $ and $\sigma  \in  \setfuncs n\kappa$ such that
$A_n\wkrestrict \sigma $ is not $U$\snug-small, and hence is
$U$\snug-large. Let $T$ be a $U$\snug-branching tree (of height $\omega
$) such that $[T] \subseteq  A_n\wkrestrict \sigma $.  The set $\{\gamma
\suchthat  \langle \gamma \rangle \in  T\}$ is in~$U$, so there exist
distinct $\gamma_\beta$, $\beta  < \lambda$, such that $\langle
\gamma_\beta\rangle \in  T$ for all $\beta $.  For each $\beta  <
\lambda $, let $T_\beta = \{\tau \suchthat  \langle
\gamma_\beta\rangle\concat \tau \in  T\}$; then $T_\beta$ is a
$U$\snug-branching tree, so $T_\beta\cup [T_\beta]$ is a closed
$U$\snug-branching tree of height $\omega +1$. For any $\mu' < \mu$,
since $U$ is $(\lambda,\mu';\omega)$\snug-nonregular, there is $S
\subseteq \lambda $ of cardinality~$\mu'$ such that $\{T_\beta\cup
[T_\beta]\suchthat \beta  \in  S\}$ has a common maximal element $z$.
For each $\beta  \in S$ we have $z \in  [T_\beta]$, so $\langle
\gamma_\beta\rangle\concat z \in  [T] \subseteq  A_n\wkrestrict \sigma
$, so $\sigma \concat \langle \gamma_\beta\rangle\concat z \in  A_n$.
The $\mu'$~points $\sigma \concat \langle \gamma_\beta\rangle\concat z$
for $\beta \in S$ are all on the same line parallel to the $n$\snug'th
coordinate axis, so $A_n$ is not $\mu'$\snug-narrow in the $n$\snug'th
coordinate. Since $\mu'$~was arbitrary, we are done.\QED\enddemo

     This proposition, together with Propositions 6.4(c) and
6.3(a), immediately gives:

\proclaim{Theorem 6.7} If $U$ is a non-principal
$\lambda$\snug-complete ultrafilter over $\kappa$, then
$\NNC(\kappa,{<}\lambda,\allowbreak\Umble)$. \QED\endproclaim

\proclaim{Corollary 6.8} If $\kappa$ is $\aleph_0$ or a measurable
cardinal, then $\NNC(\kappa,{<}\kappa,\Borel)$.\QED\endproclaim

     Corollary 6.8 for measurable cardinals also follows
from Theorem 3.1, but Theorem 6.7 gives more information for this
case.  In fact, if $U$ is a $\kappa$\snug-complete ultrafilter over
$\kappa$, then Proposition 6.3(c) easily implies that the collection
of $U$\snug-measurable sets is closed under unions and intersections
of fewer than $\kappa$ sets.  Now, any Boolean combination of certain
sets can be written as a union of intersections of these sets and
their complements; since $\kappa$ is a strong limit cardinal, we see
that any Boolean combination of fewer than $\kappa$ open subsets of
$\setfuncs {\omega}\kappa$ is $U$\snug-measurable.
So, for the measurable cardinal case, the conclusion of Theorem~6.7
subsumes that of Theorem~3.1.  In fact, the conclusion of Theorem~6.7
is strictly stronger:

\proclaim{Proposition 6.9} If $U$ is a $\kappa$\snug-complete
ultrafilter over the measurable cardinal~$\kappa$, then there are
$U$\snug-measurable subsets of\/~$\setfuncs \omega\kappa$ which cannot
be expressed as Boolean combinations of fewer than~$\kappa$ open sets.
\endproclaim

\demo{Proof} First we show that, for any $\lambda  < \kappa$, there is a
subset of\/~$\setfuncs {\omega}\kappa$ which is a Boolean combination of
fewer than~$\kappa$ open sets but not a Boolean combination of~$\lambda
$ open sets. To see this, let $\mu $ be a strong limit cardinal of
cofinality~$\omega $ such that $\lambda  < \mu  < \kappa$; then $2^\mu
= \mu^{\aleph_0}$ and $2^\lambda < \mu $.  This implies that, for any
sequence $\langle X_{\alpha}\suchthat \alpha  < \lambda \rangle$ of
subsets of\/ $\setfuncs {\omega}\mu $ and any $X \subseteq  \setfuncs
{\omega}\mu $ of cardinality~$\mu^{\aleph_0}$, there are distinct $x$
and~$y$ in~$X$ such that $\{\alpha  < \lambda \suchthat  x \in
X_{\alpha}\} = \{\alpha  < \lambda \suchthat  y \in  X_{\alpha}\}$. Let
$\delta = \mu^{\aleph_0}$. Then the number of $\lambda $\snug-sequences
of open subsets of $\setfuncs {\omega}\mu $ is $(2^\mu)^\lambda  =
2^\mu  = \mu^{\aleph_0} = \delta$, so we can
enumerate all such sequences in a sequence of length~$\delta$. Now an
easy recursive construction gives one-to-one sequences $\langle
x_\beta\suchthat  \beta  < \delta\rangle$ and $\langle y_\beta\suchthat
\beta  < \delta\rangle$ of elements of\/~$\setfuncs {\omega}\mu $ such
that $\{x_\beta\suchthat  \beta  < \delta\}\cap \{y_\beta\suchthat
\beta  < \delta\} = \nullset $ and, for any sequence $\langle
X_{\alpha}\suchthat \alpha  < \lambda \rangle$ of open subsets
of\/~$\setfuncs {\omega}\mu $, there is $\beta  < \delta$ such that $\{\alpha
< \lambda \suchthat  x_\beta \in X_{\alpha}\} = \{\alpha  < \lambda
\suchthat  y_\beta \in  X_{\alpha}\}$. So $\{x_\beta\suchthat  \beta  <
\delta\}$ is not a Boolean combination of $\lambda $~open subsets
of\/~$\setfuncs {\omega}\mu $; since the intersection of\/~$\setfuncs
{\omega}\mu $ with an open subset of\/~$\setfuncs {\omega}\kappa$ is an
open subset of\/~$\setfuncs {\omega}\mu $, $\{x_\beta\suchthat  \beta  <
\delta\}$ is not a Boolean combination of open subsets of\/~$\setfuncs
{\omega}\kappa$. But any one-element subset of\/~$\setfuncs
{\omega}\kappa$ is an intersection of $\aleph_0$~clopen subsets
of\/~$\setfuncs {\omega}\kappa$, so $\{x_\beta\suchthat \beta  < \delta\}$ is
a Boolean combination of $\delta < \kappa$ open subsets
of\/~$\setfuncs {\omega}\kappa$.

This easily implies that the collection of Boolean combinations of fewer
than~$\kappa$ open subsets of\/~$\setfuncs {\omega}\kappa$ is not closed
under the $\kappa$\snug-ary operation which takes $\langle
X_{\alpha}\suchthat  \alpha  < \kappa\rangle$ to $\{\langle \alpha
\rangle\concat s\suchthat  s \in  X_{\alpha}\}$, while the collection of
$U$\snug-measurable sets is easily seen to be closed under this
operation, so the latter collection contains sets not in the former
collection, as was to be shown.

(Another way to see that not every $U$\snug-measurable set is such a
Boolean combination is to construct a $U$\snug-measurable set which is
not $V$\snug-measurable for some other nonprincipal
$\kappa$\snug-complete ultrafilter~$V$ over~$\kappa$.  To do this,
choose $S \in  V\setdiff U$, and note that $Z = \{s \in  \setfuncs
{\omega}\kappa\suchthat (\exists m{<}\omega)(\forall n{>}m)\, s(n) \in
S\}$ is $U$\snug-null but its complement is $V$\snug-null.  By Theorems
6.7 and~2.5, there is a set $Y \subseteq  \setfuncs {\omega}\kappa$
which is not $V$\snug-measurable; then $Y\cap Z$ is $U$\snug-measurable
but not $V$\snug-measurable.) \QED\enddemo

     In the case $\kappa = \aleph_0$, Theorem~6.7 again says more than
Corollary 6.8.  For one thing, there are $2^{2^{\aleph_0}}$
$U$\snug-measurable sets but only $2^{\aleph_0}$ Borel sets.  Also,
recall the remarks after the proof of Proposition~5.10.  We now have
$\NNC(\aleph_0,{<}\aleph_0,\text{analytic})$ and more.  But we
saw that not all $\Delta^1_2$ sets are $U$\snug-measurable; in fact,
the proof of Theorem~2.5, done carefully using
a $\Sigma^1_2$\snug-good well-ordering of\/~$\setfuncs \omega\omega$
(see Moschovakis \cite{\Moschovakis,~\S5A}), shows that
$\NNC(\aleph_0,2,\Delta^1_2)$ fails in the constructible universe.

(Since this proof of $\NNC(\aleph_0,{<}\aleph_0,\text{analytic})$
uses a nonprincipal ultrafilter over~$\omega$, it would appear
to need more of the Axiom of Choice than most proofs of
similar results in descriptive set theory.  However, one can
modify the proof so that it only needs a weaker form of Choice,
such as the Axiom of Dependent Choices.  This is done by proving
versions of the results in Section~5 using the concept of
$F$\snug-measurability where $F$~is a filter rather than an
ultrafilter, and $F$~is enlarged as necessary so as to make the
relevant sets measurable.  For instance, the modified version of
Corollary~5.7 states that, for any Borel set $B \subseteq
\setfuncs\omega\kappa$ and any filter~$F$ over~$\kappa$, there
is a filter $F' \supseteq F$ such that $B$~is $F'$\snug-measurable.)

     To apply Proposition~6.6, we need to find ultrafilters which are
$(\lambda,\mu;\omega)$\snug-nonregular.  The remainder of this section
will give cases in which such ultrafilters can (or cannot) be found.

\proclaim{Proposition 6.10} If $\lambda $ is an infinite
cardinal and $U$ is an ultrafilter over $\kappa$ which
is $(\lambda,\lambda;\alpha_n)$\snug-nonregular
for each $n < \omega $, then $U$ is
$(\lambda,\aleph_0;\sum_{n<\omega}\alpha_n)$\snug-nonregular. \endproclaim

\demo{Proof} Let $\{T_\beta\suchthat  \beta  \in  \lambda \}$
be a collection of closed $U$\snug-branching trees of height
$(\sum_{n<\omega}\alpha_n)+1$.  We will recursively construct $S_n$,
$\gamma_n$, $s_n$, and $\{T^{(n)}_{\beta}\suchthat  \beta  \in  S_n\}$ for $n
< \omega $ so that $\langle \gamma_n\suchthat  n \in  \omega \rangle$ is
a sequence of distinct elements of $\lambda $ and $\{T_{\gamma_n}\suchthat
n < \omega \}$ has a common maximal element; this suffices to show that
$U$ is $(\lambda,\aleph_0;\sum_{n<\omega}\alpha_n)$\snug-nonregular.

Let $S_0 = \lambda $, $s_0 = \nullseq$, and $T^{(0)}_{\beta} = T_\beta$
for $\beta  \in  \lambda $.  Now suppose we are given $S_n \subseteq
\lambda $ of cardinality $\lambda $, a sequence $s_n$, and a collection
$\{T^{(n)}_{\beta}\suchthat \beta  \in  S_n\}$ of closed
$U$\snug-branching trees of height $(\sum_{n\le m<\omega}\alpha_m)+1$.
Let $\gamma_n$ be the least member of $S_n$.  Proposition 6.3(c) for
$\mu  = 2$ implies that, for each $\beta  \in  S_n$,
$T^{(n)}_{\beta}\cap T^{(n)}_{\gamma_n}$ is a closed $U$\snug-branching
tree of height $(\sum_{n\le m<\omega}\alpha_m)+1$, and therefore
$T^{(n)}_{\beta}\cap T^{(n)}_{\gamma_n}\cap \setfuncs
{\le\alpha_n}\kappa$ is a closed $U$\snug-branching tree of height
$\alpha_n+1$.  Since $S_n\setdiff \{\gamma_n\}$ has cardinality
$\lambda $ and $U$ is $(\lambda,\lambda;\alpha_n)$\snug-nonregular, we
can find a set $S_{n+1} \subseteq  S_n\setdiff \{\gamma_n\}$ of
cardinality $\lambda $ such that $\{T^{(n)}_{\beta}\cap
T^{(n)}_{\gamma_n}\cap \setfuncs {\le\alpha_n}\kappa\suchthat  \beta
\in  S_{n+1}\}$ has a common maximal element $s$.  For each $\beta  \in
S_{n+1}$, let $T^{(n+1)}_{\beta} = \{t\suchthat  s\concat t \in
T^{(n)}_{\beta}\cap T^{(n)}_{\gamma_n}\}$; then $T^{(n+1)}_{\beta}$ is
a closed $U$\snug-branching tree of height
$(\sum_{n<m<\omega}\alpha_m)+1$.  Let $s_{n+1} = s_n\concat s$. This
completes the recursive definition.

Clearly $s_m \subseteq  s_n$ for $m < n < \omega $; also, $\gamma_m
     \in  S_m$ but
$\gamma_n \notin  S_m$ for $n > m$, so $\gamma_m \ne  \gamma_n$ for $m
< n < \omega $.  It is easy to show by induction on $m$ that
     $$T^{(m)}_{\beta} = \{t\suchthat  s_m\concat t \in  T_\beta\cap
     \bigcap_{n<m}T_{\gamma_n}\}$$
for each $\beta  \in  S_m$; it follows immediately that $s_m \in
T_{\gamma_n}$ for all $m, n \in  \omega $.  Therefore, if we let $s =
\bigcup_{m<\omega}s_m$, then $s$ will be a common maximal element of
$\{T_{\gamma_n}\suchthat  n < \omega \}$, as desired. \QED\enddemo

\proclaim{Corollary 6.11} If $\lambda $ is infinite, then
any $(\lambda,\lambda)$\snug-nonregular ultrafilter is
$(\lambda,\aleph_0;\omega)$\snug-nonreg\-ular.\QED\endproclaim

We will see later that $(\lambda,\aleph_1;\omega)$\snug-nonregularity
does not
follow from $(\lambda,\lambda)$\snug-nonregularity.  As to the problem
of finding $(\lambda,\lambda)$\snug-nonregular ultrafilters, Silver has
shown \cite{\Jech, Ex.~34.4, p.~426} that, if $\lambda $ is regular,
any $\lambda $\snug-saturated $\lambda^+$\snug-complete ideal $I$
over $\kappa$ has the property that any collection of $\lambda $
subsets of $\kappa$ not in $I$ has a subcollection of size $\lambda
$ with nonempty intersection; this property clearly implies that any
ultrafilter over $\kappa$ disjoint from $I$ (i.e., extending the filter
dual to $I$) is $(\lambda,\lambda)$\snug-nonregular.  It follows that if
$\kappa$ is a cardinal carrying a nonprincipal $\lambda $\snug-saturated
$\lambda^+$\snug-complete ideal, then $\NNC(\kappa,\aleph_0,\Borel)$;
in particular, if there is a real-valued measurable cardinal, then
$\NNC(2^{\aleph_0},\aleph_0,\Borel)$.

To get $\NNC(\kappa,\lambda,\Borel)$ for uncountable $\lambda $ by
this method,
we need ultrafilters with stronger properties.

\proclaim{Proposition 6.12} Suppose $\kappa$, $\lambda $, and $\mu
$ are cardinals, and $U$ is an ultrafilter over $\kappa$ with
the following property: for any $S \subseteq  U$ of cardinality
at most $\lambda $, there is a set $X$ of cardinality at most
$\mu $ such that $X\cap A \ne  \nullset $ for all $A \in  S$.
Then, for any $\alpha $ such that $\mu^{|\alpha|} < \lambda $, $U$
is $(\lambda,\lambda ';\alpha)$\snug-nonregular for all $\lambda
' < \lambda $; if $\mu^{|\alpha|} < \cf \lambda $, then $U$ is
$(\lambda,\lambda;\alpha)$\snug-nonregular. \endproclaim

\demo{Proof} We may assume $\lambda $ is infinite, since otherwise
any ultrafilter is $(\lambda,\lambda;\alpha)$\snug-nonregular.
Fix $\alpha $, and let $\{T_\beta\suchthat  \beta  < \lambda \}$ be a
collection of closed $U$\snug-branching trees of height $\alpha +1$.
For each $\beta  < \lambda $ we will define a maximal branch $s_\beta$
of\/ $T_\beta$.  The definition will be by simultaneous recursion on
the length of the sequences.  So suppose $\gamma  < \alpha $ and we
have defined $s_\beta\restrict \gamma  \in  T_\beta$ for each $\beta
< \lambda $.  Fix $s \in  \setfuncs \gamma\kappa$, and let $S = \{\beta  <
\lambda \suchthat  s_\beta\restrict\gamma = s\}$.  Then $\{\{\delta\suchthat
s\concat \langle \delta\rangle \in  T_\beta\}\suchthat  \beta  \in  S\}$
is a collection of at most $\lambda $ members of $U$, so there is
a set $X \subseteq  \kappa$ of cardinality at most $\mu $ which has
nonempty intersection with each member of this collection.  For each
$\beta  \in  S$, define $s_\beta(\gamma)$ to be the least member of
$X\cap \{\delta\suchthat  s\concat \langle \delta\rangle \in  T_\beta\}$.
Do this for all $s \in  \setfuncs \gamma\kappa$ to define $s_\beta(\gamma)$
for all $\beta  < \lambda $.  This completes the recursion.

     Clearly $s_\beta$ is a maximal element of\/~$T_\beta$ for each
$\beta  < \lambda $.  It is clear from the definition of the
sequences~$s_\beta$ that, for any~$s$, there are at most $\mu$
$\delta$\snug's such that $s\concat \langle \delta\rangle$ is an
initial segment of some~$s_\beta$.  It follows easily that, for each
$\gamma \le  \alpha $, $|\{s_\beta\restrict \gamma \suchthat  \beta  <
\lambda \}| \le  \mu^{|\gamma|}$.  In particular, $|\{s_\beta\suchthat
\beta  < \lambda \}| \le  \mu^{|\alpha|}$.  Now, if $\mu^{|\alpha|} <
\lambda $ and $\lambda ' < \lambda $, then there must be an~$s$ such
that $|\{\beta < \lambda \suchthat  s_\beta = s\}| > \lambda '$, since
otherwise we would have expressed $\lambda $ as the union of at most
$\mu^{|\alpha|}$ sets each of cardinality at most $\lambda '$, which is
impossible.  Similarly, if $\mu^{|\alpha|} < \cf \lambda $, then there
must be an~$s$ such that $|\{\beta  < \lambda \suchthat  s_\beta = s\}|
= \lambda $.  Since $\{T_\beta\suchthat  \beta  < \lambda \}$ was
arbitrary, we are done. \QED\enddemo

\proclaim{Proposition 6.13} If $U$ is an ultrafilter over $\kappa$ which
is $\nu$\snug-indecomposable for all $\nu$ such that $\mu  < \nu \le
2^\lambda$, then $U$ has the property in the hypothesis of Proposition
6.12. \endproclaim

\demo{Proof} Suppose $S \subseteq  U$ and $|S| \le  \lambda $.
Define an equivalence relation $\sim $ on $\kappa$ as follows:  for any
$\beta, \gamma  < \kappa$, $\beta  \sim  \gamma $ iff, for all $A \in
S$, we have $\beta  \in  A \iff \gamma  \in  A$.  Clearly there are at
most~$2^\lambda$ $\sim $\snug-equivalence classes, and the union of these
classes is $\kappa$, so the indecomposability of~$U$ implies that there
is a set $Z$ of at most $\mu $ $\sim $\snug-equivalence classes such
that $\bigcup Z \in  U$.  Let $X \subseteq  \bigcup Z$ be a set which
contains exactly one member of each set in $Z$; then $|X| \le  \mu $.
If $A \in  S$, then $A \in  U$, so $A\cap \bigcup Z \ne  \nullset $.
If $y \in  A\cap \bigcup Z$, then there is $x \in  X$ such that $x
\sim  y$; since $y \in  A$, $x \in  A$, so $X\cap A \ne  \nullset $.
Since $A$ was arbitrary, $X$ is the desired set.\QED\enddemo

These two propositions show that $\NNC(\kappa,\lambda,\Borel)$
follows
from the existence of sufficiently indecomposable ultrafilters over
$\kappa$; in particular, if $\kappa$ is a
strong limit cardinal carrying a uniform ultrafilter $U$ which is
$\nu$\snug-indecomposable for all sufficiently large $\nu < \kappa$,
then $\NNC(\kappa,{<}\kappa,\text{$U$\snug-measurable})$ (and hence
$\NNC(\kappa,{<}\kappa,\Borel)$).

One application of these propositions is to show that certain
cardinals~$\kappa$ which satisfy $\NNC(\kappa,\lambda,\Borel)$
for all $\lambda<\kappa$ for trivial reasons (they are limits
of smaller cardinals with this property) actually satisfy
the stronger statement $\NNC(\kappa,{<}\kappa,\Borel)$ less trivially.

\proclaim{Corollary 6.14} If $\kappa$ is a limit of an $\omega
$\snug-sequence of measurable cardinals, or if $\kappa$ is the cardinal
obtained by adjoining a Prikry sequence through a measurable cardinal,
then $\NNC(\kappa,{<}\kappa,\Borel)$. \endproclaim

\demo{Proof} In each case there is a uniform ultrafilter over $\kappa$
which is $\mu $\snug-indecomposable for all~$\mu $ such that $\aleph_0 <
\mu  < \kappa$.  For the first case, let $U_n$ be a
$\kappa_n$\snug-complete nonprincipal ultrafilter over $\kappa_n$, where
$\langle \kappa_n\suchthat  n < \omega \rangle$ converges to $\kappa$,
and let $V$ be a nonprincipal ultrafilter over~$\omega $; then let $U =
\{S \subseteq  \kappa\suchthat  \{n < \omega \suchthat  S\cap \kappa_n
\in  U_n\} \in  V\}$.  To see that the ultrafilter $U$ is $\mu
$\snug-indecomposable for all $\mu $ such that $\aleph_0 < \mu  <
\kappa$, note that if $(\bigcup_{\alpha<\mu}S_{\alpha})\cap \kappa_n \in
U_n$ and $\mu  < \kappa_n$, then $S_{\alpha}\cap \kappa_n \in  U_n$ for
some $\alpha  < \mu $.  For the second case, let $U$ be any ultrafilter
extending the $\kappa$\snug-complete ultrafilter over $\kappa$ in the
ground model used to define the forcing notion; Prikry~\cite{\Prikry}
shows that $U$ is $\lambda $\snug-indecomposable for all uncountable
$\lambda  < \kappa$ (see Jech~\cite{\Jech, Ex.~37.3}).\QED\enddemo

     The cardinality hypotheses of Propositions~6.12 and~6.13 prevent us
from applying them to get new results about cardinals $\kappa \le
2^{\aleph_0}$.  The next two propositions (6.16 in particular) will show
that even the assumption that $\kappa$ is real-valued measurable is not
strong enough to get a uniform
$(\kappa,\aleph_1;\omega)$\snug-nonregular ultrafilter over $\kappa$.

\proclaim{Proposition 6.15} If $\lambda $ and $\mu $ are cardinals
and there is a set $S \subseteq  \setfuncs {\omega}\omega $ of cardinality
$\lambda $ such that, for any compact $C \subseteq  \setfuncs {\omega}\omega
$, $|S\cap C| < \mu $, then any $\aleph_1$\snug-incomplete ultrafilter
is $(\lambda,\mu;\omega)$\snug-regular. \endproclaim

\demo{Proof} Fix such an $S$, say $S = \{x_\beta\suchthat  \beta  <
\lambda \}$ with $x_\beta \ne  x_\gamma$ for $\beta  \ne  \gamma $, and
fix an $\aleph_1$\snug-incomplete ultrafilter $U$ over $\kappa$. Let
$\langle Y_n\suchthat  n < \omega \rangle$ be a sequence of sets in $U$
such that $\bigcap_{n<\omega}Y_n \notin  U$; we may assume that
$\bigcap_{n<\omega}Y_n = \nullset $ and that $Y_n \supseteq  Y_{n+1}$
for $n < \omega $.  Define a sequence $\langle T_\beta\suchthat  \beta <
\lambda \rangle$ of closed $U$\snug-branching trees of height $\omega
+1$ as follows:  for any $\beta  < \lambda $ and any $s \in \setfuncs
{\le\omega}\kappa$, put $s \in  T_\beta$ iff, for each $n < \ell (s)$,
$s(n) \in  Y_{x_\beta(n)}$.  To see that the trees $T_\beta$ have the
required properties, let $s$ be any element of\/ $\setfuncs
{\omega}\kappa$, and let $Z = \{\beta  < \lambda \suchthat  s \in
T_\beta\}$; we must see that $|Z| < \mu $.  Define $t \in  \setfuncs
{\omega}\omega $ so that $t(n)$ is the least $m$ such that $s(n) \notin
Y_m$; since $\langle Y_m\suchthat m < \omega \rangle$ is a decreasing
sequence, it is easy to see that $x_\beta(n) < t(n)$ for all $n < \omega
$ and $\beta  \in  Z$.  Let $C = \{x \in  \setfuncs {\omega}\omega
\suchthat  (\forall n{<}\omega)\, x(n) < t(n)\}$; then $C$ is a compact
set, so $|Z| = |\{x_\beta\suchthat  \beta  \in  Z\}| \le  |S\cap C| <
\mu $, as desired.\QED\enddemo

Prikry~\cite{\Prikry} (see also Jech \cite{\Jech,
pp.~425-426}) has shown that, in the model obtained by adding
$\lambda $ Cohen-generic reals to a model containing a measurable
cardinal $\kappa$, $\kappa$ carries a $\kappa$\snug-complete
$\aleph_1$\snug-saturated ideal.  However, he has also shown that
this model satisfies the hypothesis of Proposition 6.15 with $\mu
= \aleph_1$; therefore, the existence of a $\kappa$\snug-complete
$\aleph_1$\snug-saturated ideal over $\kappa$ does not imply the
existence of a $(\kappa,\aleph_1;\omega)$\snug-nonregular ultrafilter.
On the other hand, the hypothesis of Proposition 6.15 cannot hold if
there is a real-valued measurable cardinal $\kappa$ such that $\lambda
\ge  \kappa \ge  \mu  \ge  \aleph_1$; this limitation does not apply to
the following proposition.

\proclaim{Proposition 6.16} If $\kappa$, $\lambda $, and $\mu $ are
cardinals and there is a sequence $\langle A_\beta\suchthat  \beta  <
\lambda \rangle$ such that, for any infinite $S \subseteq  \kappa$,
$|\{\beta  < \lambda \suchthat  S \subseteq  A_\beta$ or $S\cap A_\beta =
\nullset \}| < \mu $, then any nonprincipal ultrafilter over $\kappa$
is $(\lambda,\mu;\omega)$\snug-regular. \endproclaim

\demo{Proof} Assume the hypothesis, and let $U$ be a nonprincipal
ultrafilter over $\kappa$.  Define closed $U$\snug-branching trees
$T_\beta$ of height $\omega +1$ for $\beta  < \lambda $ as follows:  for
any $\beta  < \lambda $ and any $s \in  \setfuncs {\le \omega}\kappa$,
put $s \in  T_\beta$ iff $s$ is one-to-one and $s(n) \in  B_\beta$ for
each $n < \ell (s)$, where $B_\beta$ is that one of $\kappa\cap A_\beta$
and $\kappa\setdiff A_\beta$ which is in $U$.  To see that these trees
have the required property, let $s$ be any element of\/ $\setfuncs
{\omega}\kappa$, and let $Z = \{\beta  < \lambda \suchthat  s \in
T_\beta\}$; we must see that $|Z| < \mu $.  We may assume that $s$ is
one-to-one, since otherwize $Z = \nullset $.  Let $S = \{s(n)\suchthat
n < \omega \}$; then $S$~is infinite and, for each $\beta  \in  Z$, $S
\subseteq  B_\beta$, so $S \subseteq  A_\beta$ or $S\cap A_\beta =
\nullset $.  Hence, $|Z| < \mu $, as desired.\QED\enddemo

     If $M$ is a model obtained by adding a sequence $\langle
x_\beta\suchthat  \beta  < \lambda \rangle$ of Cohen-generic members
of\/ $\setfuncs {\omega}2$ (or $\setfuncs {\omega}\omega $) to some
ground model, then $M$ satisfies the hypothesis of Proposition~6.16 for
$\aleph_0 \le  \kappa \le  \lambda $ and $\mu  = \aleph_1$.  To define the
sequence $\langle A_\beta\suchthat  \beta  < \lambda \rangle$, let
$f\funcfrom  \lambda \times \lambda \to \lambda $ be a bijection which
is in the ground model, and let
$$A_\beta = \{\gamma  < \lambda \suchthat  x_{f(\beta,\gamma)}(0) = 1\}.$$
To see that this works, let $S$ be an infinite subset of $\lambda $
in~$M$; since the desired property for~$S$ follows from that property
for some infinite subset of~$S$, we may assume that $S$~is countable.
Since the forcing notion has the countable chain condition,
there is a countable set $S'$ in the ground model such that $S \subseteq
S' \subseteq  \lambda $.  For each $\alpha  \in  S'$, choose a maximal
antichain (in the ground model) of conditions which decide whether
$\alpha  \in  S$; each of these antichains is countable.  Hence, if\/~$W$
is the set of $\beta  < \lambda $ such that some element of one of these
antichains gives some information about $x_\beta$, then $W$~is countable,
and so is $Z = \{\beta  < \lambda \suchthat  (\exists \gamma {<}\lambda)\,
f(\beta,\gamma) \in  W\}$.  An easy genericity argument shows that, for
any $\beta  \in  \lambda \setdiff Z$, $S\cap A_\beta \ne  \nullset $ and
$S\setdiff A_\beta \ne  \nullset $; since $Z$~is countable, we are done.

     This gives another proof that $\kappa$ can carry a
$\kappa$\snug-complete $\aleph_1$\snug-saturated ideal without carrying
a $(\kappa,\aleph_1;\omega)$\snug-nonregular ultrafilter.  In this case,
however, an analogous proof works to give a real-valued measurable
cardinal $\kappa$ carrying no $(\kappa,\aleph_1;\omega)$\snug-nonregular
ultrafilter.  The model $M$ for this case is obtained by forcing to add
$\lambda$~random reals; specifically we define this forcing notion using
as conditions the subsets of\/ $\setfuncs \lambda2$ of positive measure
in the symmetric product measure on~$\setfuncs\lambda2$.  Let $G$ be the
generic set, and find $f$ in the ground model as in the preceding
paragraph; then let
     $$A_\beta = \{\gamma  < \lambda \suchthat  \{s \in  \setfuncs
     \lambda2\suchthat s(f(\beta,\gamma)) = 1\} \in  G\}.$$
The proof that this works is the same as before, once we note that each
measurable set has countable support.  But any measurable cardinal in
the ground model is real-valued measurable in $M$.

     The preceding two propositions and the associated remarks do not
preclude the existence of a cardinal $\kappa \le 2^{\aleph_0}$ carrying
a $(\kappa,\aleph_1;\omega)$\snug-nonregular ultrafilter.  In
particular, the hypotheses of both propositions contradict Martin's
Axiom (MA), given some mild hypotheses (namely $\mu \le  \lambda  \le
2^{\aleph_0}$, $\mu  < 2^{\aleph_0}$, $\cf \mu  > \omega $, and $\kappa
\ge  \aleph_0$).  For 6.15 we recall that MA implies that for any set $S
\subseteq  \setfuncs {\omega}\omega $ of cardinality $\mu $ there is $g
\in  \setfuncs {\omega}\omega $ such that, for each $f \in  S$, $f(n)
\le  g(n)$ for all sufficiently large $n$ \cite{\Jech, p.~261}. Hence,
$S$ is contained in the union of\/ $\aleph_0$ compact sets, namely $\{f
\in  \setfuncs {\omega}\omega \suchthat  (\forall n{<}\omega)\, f(n) \le
h(n)\}$ for all $h \in  \setfuncs {\omega}\omega $ such that $h(n) =
g(n)$ for all sufficiently large $n$, so one of these compact sets must
contain $\mu$~members of $S$.  For 6.16, we use the following argument,
which Baumgartner and Hajnal \cite{\BaumgartnerHajnal, p.~196} attribute
to Solovay.  Let $\langle A_\beta\suchthat  \beta  < \lambda \rangle$ be
arbitrary.  Let $U$ be a nonprincipal ultrafilter over $\omega $; for
each $\beta  < \mu $, let $B_\beta$ be that one of $\omega \cap A_\beta$
and $\omega \setdiff A_\beta$ which is in~$U$. Define a forcing notion
$P$ as follows:  a condition is a pair $(a,b)$ where $a \subseteq
\omega $ and $b \subseteq  \mu $ are finite; $(c,d)$ is stronger than
$(a,b)$ iff $a \subseteq  c$, $b \subseteq  d$, and, for each $n \in
c\setdiff a$ and each $\beta  \in  b$, $n \in  B_\beta$. Since any two
conditions $(a,b), (c,d)$ with $a = c$ are compatible, $P$~has the
countable chain condition.  Using the definition of $B_\beta$, we easily
see that the sets $D_n = \{(a,b) \in  P\suchthat  |a| \ge  n\}$ and
$E_\beta = \{(a,b) \in  P\suchthat  \beta  \in  b\}$ are dense in $P$
(for $n < \omega $, $\beta  < \mu $).  Now apply MA to get a filter $G$
on $P$ which meets each of these dense sets.  Let $S = \bigcup
\{a\suchthat  (a,b) \in  G\}$; then $S$ is infinite and, for each $\beta
< \mu $, $S\setdiff B_\beta$ is finite.  Since $S$~has only countably
many finite subsets, there must be a finite $z \subseteq  S$ such that
$|\{\beta  < \mu \suchthat  S\setdiff z \subseteq  B_\beta\}| = \mu $,
and hence $|\{\beta  < \lambda \suchthat  S\setdiff z \subseteq
A_\beta$ or $(S\setdiff z)\cap A_\beta = \nullset \}| \ge  \mu $.

     But starting with a model containing a cardinal $\kappa$
carrying a $\kappa$\snug-complete $\aleph_1$\snug-saturated ideal, one
can obtain a model of MA${}+{}$\snug``$2^{\aleph_0}$ is large'' by a
c.c.c.\ forcing extension \cite{\Jech, \S23}, and $\kappa$ will still
carry a $\kappa$\snug-complete $\aleph_1$\snug-saturated ideal in the
extension \cite{\Jech, Ex.~34.5, p.~426}.  It is quite possible that
this model, or a model obtained by some more specialized c.c.c.\ forcing
notion over a model with a measurable cardinal, will contain a cardinal
$\kappa$ carrying a nontrivial $(\kappa,\aleph_1;\omega)$\snug-nonregular
ultrafilter.  Another possibility is that $\NNC(\kappa,\aleph_1,\Borel)$
will follow from the real-valued measurability of $\kappa$ by a
different argument.
(By the random-real case of Theorem~4.1, we know that
$\NNC(\kappa,\aleph_1,\Borel)$ is at least relatively consistent with
the real-valued measurability of $\kappa$.)

%

\head 7.  The Complexity of Narrow Clopen Partitions \endhead

In this section, we consider a slightly different question.
Let $\kappa$ and~$\lambda$ be infinite cardinals.
Suppose that there does exist a $\lambda$\snug-narrow
covering of\/~$\setfuncs \omega\kappa$ by open sets.  Must such
a covering be complicated?

Of course, we cannot ask this without a suitable measure of
complexity of open coverings of\/~$\setfuncs \omega\kappa$.
We can get such a measure by considering trees associated with the
coverings.

For any finite sequence $\sigma \in \setfuncs {<\omega}\kappa$, let
$\Nbhd(\sigma)$ be the basic open subset of\/~$\setfuncs \omega\kappa$
consisting of those infinite sequences that begin with~$\sigma$. Now,
given an open covering of\/~$\setfuncs \omega\kappa$, let $T$ be the set
of all $\sigma \in \setfuncs {<\omega}\kappa$ such that
$\Nbhd(\sigma)$~is {\sl not} a subset of any member of the covering.
Clearly $T$ is a tree, since $\Nbhd(\tau) \subseteq \Nbhd(\sigma)$
whenever $\sigma \subseteq \tau$.  Furthermore, $[T]$~is empty: any $s
\in \setfuncs \omega\kappa$ is in some member~$A$ of the covering, and
$A$~is open, so some basic neighborhood~$\Nbhd(\sigma)$ of~$s$ is
included in~$A$, which gives $\sigma \subseteq s$ and $\sigma \notin T$,
so $s \notin [T]$.

We recall some basic definitions in order to fix notation.
A tree $T \subseteq  \setfuncs {<\omega}\kappa$ is
{\it well-founded} iff\/ $[T] = \nullset $.  For any well-founded tree, we
define a rank function $\rk_T$ mapping $T$ to the ordinals by well-founded
recursion as follows:  if $\sigma  \in  T$, then $\rk_T(\sigma)$ is the
least ordinal greater than $\rk_T(\sigma \concat \langle \beta \rangle)$
for all $\beta $ such that $\sigma \concat \langle \beta \rangle \in  T$.
For $\sigma  \notin  T$ we put $\rk_T(\sigma) = -1$.  Define $\rk(T)$,
the rank of the well-founded tree~$T$,
to be $1+\rk_T(\nullseq)$.

Now we can define the complexity (or rank) of an open covering
of\/~$\setfuncs\omega\kappa$ to be the rank of the associated
well-founded tree.  This will be an ordinal less than~$\kappa^+$.

This may be slightly clearer when the open covering is actually
a partition of\/~$\setfuncs \omega\kappa$ into open sets.  In this case
the sets are necessarily clopen, since the complement of one
set is the union of the others.  And the tree~$T$ can be defined
to be the set of all~$\sigma$ such that $\Nbhd(\sigma)$~meets
more than one set in the partition.

By a standard argument, any narrow covering by open sets
can be reduced to a narrow partition:

\proclaim{Proposition 7.1} If there exists a $\lambda$\snug-narrow
covering of\/~$\setfuncs \omega\kappa$ using open sets, then there exists
a $\lambda$\snug-narrow
partition of\/~$\setfuncs \omega\kappa$ using open (and hence clopen)
sets.  \endproclaim

\demo{Proof}
Let $\langle A_n \suchthat n < \omega \rangle$ be such a covering.
Since $A_n$~is open, we have $A_n = \bigcup_{m < \omega} B_{nm}$,
where $$B_{nm} = \bigcup \{\Nbhd(\sigma) \suchthat \sigma \in
\setfuncs m\kappa,\,\Nbhd(\sigma) \subseteq A_n\}.$$
The sets~$B_{nm}$ are clopen; in fact, membership of $s \in \setfuncs
\omega\kappa$ in $B_{nm}$ depends only on $s\restrict m$.
Hence, the sets
$$B'_{nm} = B_{nm} \setminus \biggl(\bigcup_{n' < n} B_{n'm} \cup
\bigcup_{m' < m} \bigcup_{n' < \omega} B_{n'm'}\biggr)$$
are also clopen.  The sets $B'_{nm}$ are disjoint, and we have
$B'_{nm} \subseteq B_{nm}$ and $\bigcup_{n,m<\omega} B'_{nm} =
\bigcup_{n,m < \omega} B_{nm} = \setfuncs \omega\kappa$.
Therefore, the open sets $A'_n = \bigcup_{m < \omega} B'_{nm}$
form a partition of~$\setfuncs\omega\kappa$; since $A'_n \subseteq
A_n$, $A'_n$~is $\lambda$\snug-narrow, as desired.
\QED\enddemo

It is easy to see that the tree associated with the clopen partition
constructed above is the same as the tree associated with the
original open covering, so the reduction process does not
change the complexity of the covering.  Hence, we may restrict
ourselves to clopen partitions when trying to find the
minimum complexity of a $\lambda$\snug-narrow open covering
of\/~$\setfuncs \omega\kappa$.  (However, often it will be just as
convenient to work with the open coverings.)

This notion of complexity, for the case of individual
clopen subsets of the Baire space~$\setfuncs \omega\omega$,
is called the Kalmar rank; see
Barnes~\cite{\Barnes}.

If $n$~is finite, then a clopen partition of\/~$\setfuncs\omega\kappa$
has rank at most~$n$ if and only if, for every $s \in \setfuncs
\omega\kappa$, the piece of the partition that contains~$s$ is
determined by~$s \restrict n$.  Such a partition is essentially a
partition of the finite-dimensional product~$\setfuncs n\kappa$.  Also,
if a set in such a partition contains a point~$s$, it must contain {\sl
all} points on the line through~$s$ parallel to the $j$\snug'th
coordinate axis, for any $j \ge n$. Therefore, if the clopen sets~$A_j$
for $j < \omega$ form a $\lambda$\snug-narrow partition of\/~$\setfuncs
\omega\kappa$ of rank at most~$n$, where $\lambda \le \kappa$, then
necessarily $A_j = \nullset$ for $j \ge n$, and the sets~$A_j$ for $j <
n$ are determined by a $\lambda$\snug-narrow partition of\/~$\setfuncs
n\kappa$. Hence, Theorem~2.3 tells us when such partitions exist:

\proclaim{Proposition 7.2} For any natural number $n > 0$, ordinal
$\alpha$, and cardinal~$\kappa$, there exists an
$\aleph_\alpha$\snug-narrow clopen partition of\/~$\setfuncs
\omega\kappa$ of complexity at most~$n$ if and only if $\kappa <
\aleph_{\alpha+n-1}$. \QED\endproclaim

Similarly, one can translate Proposition~2.4 into a statement about
finite-rank narrow clopen partitions of\/~$\setfuncs\omega\kappa$ for
finite~$\kappa$.

The next case to consider is $\kappa=\aleph_\omega$.
Here Proposition~7.2 tells us that, for any
$\lambda < \aleph_\omega$, a $\lambda$\snug-narrow
clopen partition of\/~$\setfuncs \omega\kappa$, if it exists, must have
rank at least~$\omega$.  However, we will see that
the rank must actually be much higher than this.

Such a partition might not exist at all; see Theorem~3.5.
On the other hand, there are models in which such
partitions do exist; for instance,
Corollary~3.4(e) (along with Proposition~7.1)
tells us that, if $V=L$, then an $\aleph_1$\snug-narrow
clopen partition of\/~$\setfuncs \omega\kappa$ exists.
So, in such a model, one can try to find the least possible
complexity of such a partition.

It is convenient to reformulate this question in terms
of free subsets of algebras, as in Theorem~3.3.
Given a structure~$M$, one can form the tree~$T_M$ of
all finite sequences~$\sigma$ of members of~$M$
which are free for~$M$ (i.e., $\sigma$~is one-to-one and
the range of~$\sigma$ is a free subset of~$M$).  If $M$~has
no infinite free subset, then there can be no infinite branch
through~$T_M$, so $T_M$~is well-founded, and one can
compute its rank.

\proclaim{Proposition 7.3} Let $\kappa$ and~$\mu$ be infinite
cardinals, with $\kappa > \mu$, such that
$\NNC(\kappa,\mu^+,\allowbreak\open)$ does not hold.  Then
the least possible complexity of a $\mu^+$\snug-narrow
covering of\/~$\setfuncs \omega\kappa$ by open sets
is equal to the least possible rank of the tree~$T_M$
of finite free sequences for an algebra~$M$ of size~$\kappa$
with $\mu$~operations and no infinite free subset. \endproclaim

\demo{Proof} Given a $\mu^+$\snug-narrow covering $\langle A_n \suchthat
n < \omega \rangle$ of\/~$\setfuncs \omega\kappa$ by open sets, let $T$~be
the associated tree.  Define a structure~$M$ from the covering as in the
second part of the proof of Theorem~3.3.  If $\sigma \in \setfuncs
{<\omega}\kappa$ is not in~$T$, then $\Nbhd(\sigma) \subseteq A_m$
for some~$m$.  Let $n = \ell(\sigma)$.  Since $\kappa \ge \mu^+$,
$\Nbhd(\sigma)$ is not $\mu^+$\snug-narrow in the $j$\snug'th coordinate
for $j \ge n$, so we must have $m < n$.  Therefore, if
$\sigma'$ is $\sigma$ with the $m$\snug'th coordinate deleted, then
$\sigma(m) = f_{\alpha mn}(\sigma')$ for some $\alpha < \mu$,
so $\sigma \notin T_M$.  This proves that $T_M \subseteq T$,
so $\rk(T_M) \le \rk(T)$.

Conversely, suppose we have a structure~$M$ with universe~$\kappa$ which
has $\mu$~operations and no infinite free subset.  Define a corresponding
$\mu^+$\snug-narrow open covering $\langle A_n \suchthat n < \omega
\rangle$ of\/~$\setfuncs \omega\kappa$ as in the first part of the proof
of Theorem~3.3, and let $T$ be the associated tree.  If $\sigma \in
\setfuncs {<\omega}\kappa$ is not in~$T_M$, then, for some
$m < \ell(\sigma)$, $\sigma(m)$ is generated in~$M$ from
$\{\sigma(j) \suchthat j \ne m\}$.  This implies
$\Nbhd(\sigma) \subseteq A_m$, so $\sigma \notin T$.
Therefore, $T \subseteq T_m$, so $\rk(T) \le \rk(T_M)$.
\QED\enddemo

So we can study trees associated with coverings or trees of finite free
sequences, which\-ever is more convenient at the time.

We will now see that, when $\kappa \ge \aleph_\omega$, the
trees above must have rank much higher than the finite ranks
produced in Proposition~7.2.

\proclaim{Theorem 7.4} Let $\kappa$ be an uncountable limit cardinal. If
$M$~is an algebra with universe~$\kappa$ which has fewer than
$\kappa$~operations and no infinite free subset, then $\rk(T_M) \ge
\kappa$. \endproclaim

\demo{Proof}
First note that, if $T$~is a well-founded tree of rank~$\alpha$,
then there is a subtree $T' \subseteq T$ such that
$\rk(T') = \alpha$ and $|T'| = |\alpha|$.  This is proved by
induction on $\alpha$; it is trivial for $\alpha \le 1$.
Assume it is true for all $\beta < \alpha$,
and let $T$~be a tree of rank~$\alpha$, where $\alpha > 1$.
If $\alpha = \beta+1$, choose~$c$ such that $\langle c \rangle \in T$
and $\rk_T(\langle \rangle) = \rk_T(\langle c \rangle)+1$.  Let
$T_c = \{\sigma\suchthat \langle c \rangle \concat \sigma \in T\}$.
Then $\rk(T_c) = \beta$, so we can apply the induction hypothesis
to get $T'_c \subseteq T_c$ with $\rk(T'_c) = \beta$ and
$|T'_c| = |\beta|$.  Let $T' = \{\nullseq \} \cup \{ \langle c
\rangle \concat \sigma \suchthat \sigma \in T'_c\}$; then $T'$ is the
desired subtree of $T$.  If $\alpha$ is a limit ordinal, choose a
set~$C$ of size at most~$|\alpha|$ such that $\{ \rk_T(\langle c \rangle
\suchthat c \in C \}$ is cofinal in~$\alpha$.  Apply the induction
hypothesis to each~$T_c$ to get~$T'_c$ as above; then the tree
$T' = \{\langle \rangle \} \cup \{ \langle c
\rangle \concat \sigma \suchthat c \in C,\,\sigma \in T'_c\}$ will
be as desired.

We now prove the theorem by showing by induction on ordinals $\alpha <
\kappa$ that, if $M$~is an algebra with universe~$\kappa$ which has
fewer than $\kappa$~operations and no infinite free subset, then
$\rk(T_M) > \alpha$.  Suppose this is true for all $\alpha' < \alpha$.
Let $M$ be such an algebra.  Let $\lambda$ be an uncountable regular
cardinal less than~$\kappa$ but greater than~$\alpha$ and greater than
the number of operations of~$M$.  Let $M'$ be~$M$ with an additional
constant function~$c_\gamma$ with value~$\gamma$ for each $\gamma <
\lambda$.  By the induction hypothesis, $\rk(T_{M'})$ is greater
than~$\alpha'$ for all $\alpha' < \alpha$, so $\rk(T_{M'}) \ge \alpha$.
If $\rk(T_{M'}) > \alpha$, then $\rk(T_M) > \alpha$ as desired because
$T_{M'} \subseteq T_M$, so suppose $\rk(T_{M'}) = \alpha$.
Let $T'$ be a subtree of~$T_{M'}$ which has rank~$\alpha$ and
cardinality~$|\alpha|$.  Let $S$ be the set of members of~$\kappa$
which are mentioned in~$T'$.  Then $|S| < \lambda$ and $M$~has fewer
than~$\lambda$ operations, so $|\subalggen MS| < \lambda$ (recall that
$\subalggen MS$ is the subalgebra of~$M$ generated by~$S$).  Choose
$\gamma < \lambda$ which is not in~$\subalggen MS$.  Then, for every
$\sigma \in T'$, $\langle \gamma \rangle \concat \sigma$ is free
for~$M$.  (By choice of $\gamma$, $\gamma$ is not generated by the
members of~$\sigma$; and no member of~$\sigma$ is generated
from~$\gamma$ and the other members of~$\sigma$ because $\sigma$~is free
for~$M'$.)  Therefore, $\rk_{T_M}(\langle \gamma \rangle) \ge
\rk_{T'}(\langle \rangle)$, so $\rk(T_M) > \rk(T') = \alpha$, as
desired.  This completes the induction.
\QED\enddemo

This argument for limit cardinals produces very little when applied to
successor cardinals; in fact, the following proposition shows that the
ranks obtained from algebras of successor cardinal size are only
slightly higher than those obtained from the preceding limit cardinal.

\proclaim{Proposition 7.5} Let $\mu$ and~$\kappa$ be infinite cardinals
with $\mu \le \kappa$.  If there is an algebra on~$\kappa$
with $\mu$~operations
and no infinite free subset, then there is such an algebra on~$\kappa^+$
as well.  Furthermore, if $\alpha_0$ is the least possible rank for the
tree of finite free sequences for such an algebra on~$\kappa$, and
$\alpha_1$ is the corresponding least possible rank for~$\kappa^+$,
then $\alpha_0+1 \le \alpha_1 \le 2\cdot\alpha_0+1$.
\endproclaim

\demo{Proof} Let $M_0$ be an algebra on~$\kappa$ with $\mu$~operations
and no infinite free subset, such that $\rk(T_{M_0}) = \alpha_0$.  Also,
for each ordinal $\xi < \kappa^+$, let $g_\xi$ be a bijection between
$\xi+1$ and some ordinal $\le \kappa$.  Let $M_1$ be an algebra
on~$\kappa^+$ with $\mu$~operations which include: all of the operations
of~$M_0$, extended in some arbitrary manner to operations
on~$\kappa^+$; a binary operation~$G$ such that $G(\xi,\eta) =
g_\xi(\eta)$ whenever $\eta \le \xi$; and a binary operation~$G'$ such
that $G'(\xi,\eta) = g_\xi^{-1}(\eta)$ whenever $\eta \in \range(g_\xi)$.
We will see that $\rk(T_{M_1}) \le 2\cdot\alpha_0+1$.

Given two ordinals $\beta,\gamma < \kappa^+$, we can produce an ordinal
$\delta < \kappa$ by letting $\delta = G(\max(\beta,\gamma),
\min(\beta,\gamma))$.  On the other hand, if we are given $\delta$ and
the larger of $\beta$ and~$\gamma$, we can recover the other ordinal in
the pair $\{\beta,\gamma\}$, since $\min(\beta,\gamma) =
G'(\max(\beta,\gamma),\delta)$.  Now, given a finite sequence $\sigma
\in \setfuncs {<\omega}{\kappa^+}$ of length~$n$, we can produce a
finite sequence $h(\sigma) \in \setfuncs {<\omega}\kappa$ of
length~$n/2$ (rounded down) by applying the above procedure to the
pairs $\{\sigma(0),\sigma(1)\}$, $\{\sigma(2),\sigma(3)\}$, and so on.

If $\sigma$~is such that $h(\sigma) \notin T_{M_0}$, then there is $k <
\ell(h(\sigma))$ such that $h(\sigma)(k)$~is obtainable from the other
coordinates of~$h(\sigma)$ using the operations of~$M_0$.  Let $j$~be
whichever of $2k$ and~$2k+1$ has the smaller coordinate of~$\sigma$.
Then $\sigma(j)$~is obtainable from the other
coordinates of~$\sigma$ using the operations of~$M_1$: use~$G$
to obtain~$h(\sigma)(i)$ for $i \ne k$, then use the operations of~$M_1$
extending those of~$M_0$ to obtain~$h(\sigma)(k)$, then apply~$G'$ to
$\sigma(4k+1-j)$ and~$h(\sigma)(k)$ to get~$\sigma(j)$.  Therefore,
$\sigma \notin T_{M_1}$.

Now a straightforward induction on~$\rk_{T_{M_0}}(h(\sigma))$ shows
that, for any $\sigma \in T_{M_1}$, if $\ell(\sigma)$ is odd, then
$\rk_{T_{M_0}}(\sigma) \le 2 \cdot \rk_{T_{M_0}}(h(\sigma))$, and
if $\ell(\sigma)$ is even, then
$\rk_{T_{M_0}}(\sigma) \le 2 \cdot \rk_{T_{M_0}}(h(\sigma))+1$.
Therefore, $\rk(T_{M_1}) \le 2 \cdot \rk(T_{M_0}) + 1$, as desired.

For the other direction, let $M$ be an algebra on~$\kappa^+$ with
$\mu$~operations and no infinite free subset
such that $\rk(T_M) = \alpha_1$.  The subalgebra of~$M$
generated by the set $\kappa \subseteq \kappa^+$ has size~$\kappa$, so
we can choose $\gamma < \kappa^+$ which is not in this subalgebra.
Let~$M'$ be $M$~with an additional constant operation with
value~$\gamma$.  Now let $M''$ be an algebra on~$\kappa$ with
$\mu$~operations such that, for
each operation~$f$ on~$\kappa^+$ which is a composition of operations
of~$M'$, there is an operation~$\tilde f$ of~$M''$ such that, for any
$\beta_0,\dots,\beta_{n-1} < \kappa$, if $f(\beta_0,\dots,\beta_{n-1}) <
\kappa$, then $\tilde f(\beta_0,\dots,\beta_{n-1}) =
f(\beta_0,\dots,\beta_{n-1})$.  Any free set for~$M''$ will also be free
for~$M$, so $M''$~has no infinite free subset.  The tree~$T_{M''}$ must
have rank at least~$\alpha_0$.  But for any $\sigma \in T_{M''}$,
$\langle \gamma \rangle \concat \sigma$ must be in~$T_M$ (as in the
proof of the preceding proposition), so $\rk(T_M) > \rk(T_{M''})$,
so $\alpha_1 \ge \alpha_0+1$.
\QED\enddemo

Note that this multiplication on the left by~$2$ has no effect on the
limit part of the ordinal~$\alpha_0$.
Hence, if $\gamma$ is a limit ordinal, $\mu$~ is an infinite
cardinal less than~$\aleph_\gamma$, and the least possible rank for the
tree of finite free sequences for an algebra of size~$\aleph_\gamma$
with $\mu$~operations
is~$\delta+m$ where $\delta$~is a limit ordinal and $m$~is finite, then,
for any finite~$n$, one can apply Proposition~7.5 $n$~times to show that
the least possible rank~$\alpha$ for the
tree of finite free sequences for an algebra of size~$\aleph_{\gamma+n}$
with $\mu$~operations must satisfy $\delta+m+n \le \alpha \le
\delta + (m+1)2^n - 1$.  One can instead use a direct argument, rather than
an $n$\snug-fold iteration, to reduce this upper bound to
$\delta + (m+1)(n+1) - 1$.  This will suffice to determine~$\alpha$
completely if $m$~happens to be~$0$.

This shows that the main case of interest for the problem of
free-sequence tree ranks, or for complexity of open narrow coverings, is
the case of limit cardinals~$\kappa$.  Here Theorem~7.4 gives a lower
bound of~$\kappa$, but it is quite possible that this bound can be
improved; the only obvious upper bound is~$\kappa^+$ (assuming that a
suitable algebra or narrow covering exists at all).  In the rest of this
section, we will see that, for the particular case where $\kappa$~is an
uncountable strong limit cardinal of cofinality~$\omega$, the lower
bound can indeed be substantially improved.

For the rest of this section, we will make the following definitions and
assumptions:

{\narrower\narrower
Let $\kappa$ be an strong limit cardinal of cofinality~$\omega$.
Assume that we have (not necessarily fixed) sequences $\langle \kappa_n
\suchthat n < \omega \rangle$ and $\langle \lambda_n
\suchthat n < \omega \rangle$ of infinite cardinals such that $\kappa_n
\le \lambda_n$, $\kappa_{n+1} = (2^{\lambda_n})^+$, and $\lim_{n \to
\infty} \kappa_n = \kappa$.  Also, in order to make~$\kappa_0$ have the
same properties as the other cardinals~$\kappa_n$, assume that we have
infinite cardinals $\kappa_{-1} \le \lambda_{-1}$ such that $\kappa_0 =
(2^{\lambda_{-1}})^+$.

For each~$n$, let $P_n$ be the $n$\snug-fold Cartesian product $\prod_{i
= 0}^{n-1} \kappa_i$ (not the cardinal product, which would just be
$\kappa_{n-1}$).

}
We will be using primarily the cardinals~$\kappa_n$; the separate
cardinals~$\lambda_n$ are only needed in order to allow the sequence
$\langle \kappa_n \suchthat n < \omega \rangle$ to be cofinal
in~$\kappa$ even when $\kappa$~is a limit of strong limit cardinals.  If
$\kappa = \aleph_\omega$, we can just let $\lambda_n = \kappa_n$.

We will show that any narrow open covering of~$\setfuncs \omega\kappa$
must have high complexity by establishing two facts: the tree associated
with a narrow open covering must meet all `large' subproducts of the
product sets~$P_n$, and a tree of small rank cannot meet all such
subproducts.

\definition{Definition 7.6}
A finite sequence $\Yv = \langle \Yv(i) \suchthat i < n \rangle$ with
$\Yv(i) \subseteq \kappa_i$ for all $i < n$ is a {\it large sequence} if
$|\Yv(i)| = \kappa_i$ for all~$i$.

If $T \subseteq \setfuncs {<\omega}\kappa$ is a tree and $\Yv$~is a
large sequence, then $T \treerestrict \Yv$ is the subtree of~$T$
consisting of all $\sigma \in T$ such that $\sigma(i) \in \Yv(i)$ for
all $i < \min(\ell(\Yv),\ell(\sigma))$.  Also, we say that {\it $T$
avoids~$\Yv$} if $T \cap \prod_{i < \ell(\Yv)} \Yv(i) = \nullset$.

If $\Yv$ and~$\Zv$ are large sequences, then $\Zv \preceq \vec
Y$ means that $\ell(\Zv) \ge \ell(\Yv)$ and $\Zv(i) \subseteq
\Yv(i)$ for all $i < \ell(\Yv)$.
\enddefinition

Easily, if $\Zv \preceq \Yv$ and the tree~$T$ avoids~$\Yv$,
then $T$ avoids~$\Zv$.  Also, $T$ avoids~$\Yv$ if and only if
$\rk(T \treerestrict \Yv) \le \ell(\Yv)$.

\proclaim{Lemma 7.7} Let $F$ be a function from~$P_{n+1}$ to~$S$, where
$|S| < \kappa_n$.  Then there is a set $Z \subseteq \kappa_n$ of
size~$\kappa_n$ such that $F(\sigma)$~depends only on~$\sigma\restrict
n$ if $\sigma(n) \in Z$ (i.e., if $\sigma \restrict n = \sigma'
\restrict n$ and $\sigma(n),\sigma'(n) \in Z$, then $F(\sigma) =
F(\sigma')$).  Furthermore, if\/ $Y$~is a given subset of~$\kappa_n$ of
size~$\kappa_n$, then $Z$~can be taken to be a subset of\/~$Y$.
\endproclaim

\demo{Proof} For each $\beta < \kappa_n$, define $f_\beta\funcfrom
P_n\to S$ by $f_\beta(\sigma) = F(\sigma \concat \langle \beta
\rangle)$. Since $|S| < \kappa_n$, $|S| \le  2^{\lambda_{n-1}}$, so the
number of possible functions~$f_\beta$ is at most
     $$|S|^{\kappa_0\cdot\kappa_1\cdot\dots\cdot\kappa_{n-1}} \le
     (2^{\lambda_{n-1}})^{\lambda_{n-1}} = 2^{\lambda_{n-1}}.$$
Since there are $\kappa_n = (2^{\lambda_{n-1}})^+$ ordinals~$\beta$
in~$Y$ (let~$Y$ be~$\kappa_n$ if no~$Y$ is given),
there must be a set $Z \subseteq Y$ of cardinality~$\kappa_n$
such that $f_\beta = f_\gamma$ for all $\beta,\gamma \in Z$.
This~$Z$ satisfies the conclusion of the lemma. \QED\enddemo

If we have a function $F\funcfrom P_{n+m} \to S$ where $|S| < \kappa_n$,
then we can apply Lemma~7.7 repeatedly to restrict~$F$ to a subdomain on
which $F(\sigma)$~depends only on~$\sigma\restrict n$.  This can be
stated in terms of large sequences as follows:

\proclaim{Lemma 7.8} Let $F$ be a function from~$P_{n+m}$ to~$S$, where
$|S| < \kappa_n$, and let~$\Yv$ be a large sequence. Then there is a
large sequence $\Zv \preceq \Yv$ of length at least~$n+m$ such that
$\Zv(i) = \Yv(i)$ for $i < \min(n,\ell(\Yv))$ and, for $\sigma \in
\prod_{i < n+m} \Zv(i)$, $F(\sigma)$~depends only on~$\sigma \restrict
n$. \endproclaim

\demo{Proof} If the given~$\Yv$ has length less than~$n+m$, then extend
it to length~$n+m$ by letting $\Yv(i) = \kappa_i$ for larger values
of~$i$.  We now define $\Zv \preceq \Yv$ of the same length as~$\Yv$ as
follows.  Let $\Zv(i) = \Yv(i)$ if $i < n$ or $i \ge n+m$.
Also, let $F_m = F$.  If $i < m$ and we have a function $F_{i+1}\funcfrom
P_{n+i+1}\to S$, then by Lemma~7.7 we can find $\Zv(i) \subseteq
\Yv(i)$ and $F_i\funcfrom P_{n+i}\to S$ such that $|\Zv(i)| =
\kappa_{n+i}$ and $F_i(\sigma)
= F_{i+1}(\sigma \concat \langle \beta \rangle)$ for all $\sigma  \in
P_{n+i}$ and $\beta \in \Zv(i)$. Do this successively for~$i$ from~$m-1$
down to~$0$ to finish defining the required large sequence~$\Zv$.
\QED\enddemo


This argument applies just as well if $F$~is not defined on all
of~$P_{n+m}$, but only on $\prod_{i < n+m}\Yv(i)$, assuming $\ell(\Yv)
\ge n+m$.  Or one can extend~$F$ trivially to a function from all
of~$P_{n+m}$ to~$S$ and then apply the lemma as stated.

In the case $n=0$, the conclusion of Lemma~7.8 is that $F$~is constant
on the part of its domain specified by the large sequence~$\Zv$.

Using Lemma~7.8, we can prove one of the two facts mentioned earlier:

\proclaim{Proposition 7.9} If $T$~is the tree associated with a
$\kappa_0$\snug-narrow open covering of\/~$\setfuncs \omega\kappa$, then
$T$~does not avoid any large sequence. \endproclaim

\demo{Proof}
Let $\langle A_n \suchthat n < \omega \rangle$ be the narrow open
covering, and suppose that $\Yv$~is a large sequence which is avoided
by~$T$.  Then, for each $\sigma \in \prod_{i < \ell(\Yv)} \Yv(i)$, since
$\sigma \notin T$, there
exists $n < \omega$ such that $\Nbhd(\sigma) \subseteq A_n$; let
$F(\sigma)$ be the least such~$n$.  This defines a function $F \funcfrom
\prod_{i < \ell(\Yv)} \Yv(i) \to \omega$.  Apply Lemma~7.8 to get $\Zv
\preceq \Yv$ such that $F$~is constant on $\prod_{i < \ell(\Yv)}
\Zv(i)$, say with value~$\bar n$.  This means that any $s \in \setfuncs
\omega\kappa$ such that $s(i) \in \Zv(i)$ for all $i < \ell(\Zv)$ is
in~$A_{\bar n}$.  But clearly we can fix all coordinates of such an~$s$
other than the $\bar n$\snug'th coordinate, which we allow to vary, to
get $\kappa_{\bar n}$ points in~$A_{\bar n}$ on the same line parallel
to the $\bar n$\snug'th coordinate axis.  Therefore, $A_{\bar n}$ is not
$\kappa_0$\snug-narrow in the $\bar n$\snug'th coordinate, which is a
contradiction.
\QED\enddemo

It now remains to prove the other fact, that a tree of low rank must
avoid some large sequence.  This will be proved by induction on the rank
of the tree. We will give two versions of the inductive argument; the
second version will be more complicated, but will attain a better
result.

\proclaim{Proposition 7.10}
If $T \subset \setfuncs\omega\kappa$ is a well-founded tree of rank less
than~$\kappa \cdot \kappa_0$ (ordinal multiplication), and $\Yv$~is a
large sequence, then there is a
large sequence~$\Zv \preceq \Yv$ such that $T$ avoids~$\Zv$.
\endproclaim

\demo{Proof}
By induction on $\rk(T)$.  Suppose that the result is already known for
trees of rank less than~$\rk(T)$.  We consider three cases.

Case 1: $\rk(T) < \kappa$.  Choose~$n$ such that $\rk(T) < \kappa_n$.
It follows that the range of the function~$\rk_T$ has size less
than~$\kappa_n$.

Find large sequences $\Yv_0 \succeq \Yv_1 \succeq \Yv_2 \succeq \dotsb$
as follows.  Let $\Yv_0$ be~$\Yv$, extended arbitrarily if necessary so
as to have length at least~$n$.  Given $\Yv_{m-1}$, apply Lemma~7.8 to the
function $\rk_T\restrict P_{n+m}$ to get $\Yv_m \preceq \Yv_{m-1}$
such that $\Yv_m\restrict n = \Yv_{m-1} \restrict n$ and, for $\sigma
\in \prod_{i < n+m} \Yv_m(i)$, $\rk_T(\sigma)$~depends only
on~$\sigma \restrict n$.

Now, for any $\tau \in \prod_{i < n} \Yv_0(i)$ and any $m < \omega$, let
$F_m(\tau)$ be the common value of $\rk_T(\sigma)$ for $\sigma \in
\prod_{i < n+m} \Yv_m(i)$ extending~$\tau$.  We also have $F_{m-1}(\tau)
= \rk_T(\sigma \restrict (n+m-1))$ for such~$\sigma$; hence, either
$F_m(\tau) < F_{m-1}(\tau)$ or $F_m(\tau) = F_{m-1}(\tau) = -1$.
Since there is no infinite descending sequence of ordinals, for
each~$\tau$ there must be an~$m$ such that $F_m(\tau) = -1$; let
$G(\tau)$ be the least such~$m$.

Apply Lemma~7.8 again to get $\Zv_0 \preceq \Yv_0$ such that $G$~is
constant on $\prod_{i < n} \Zv_0(i)$; let $\bar m$ be the constant value
of~$G$ on this set.  Define the large sequence~$\Zv$ of
length~$\ell(\Yv_{\bar m})$ by letting $\Zv(i) = \Zv_0(i)$ for $i < n$
and $\Zv(i) = \Yv_{\bar m}(i)$ for $i \ge n$.  Then we have
$\rk_T(\sigma) = F_{\bar m}(\sigma \restrict n) = -1$ for all $\sigma
\in \prod_{i < n+\bar m} \Zv(i)$, so $T$ avoids~$\Zv$.

Case 2: $\rk(T)$ is of the form $\kappa\cdot\alpha + \beta$, where
$\alpha > 0$ and $\omega \le \beta < \kappa$. Let $T' = \{\sigma  \in
T\suchthat  \rk_T(\sigma) \ge  \kappa\cdot \alpha \}$.  Then $T'$ is a
subtree of~$T$, and an easy induction shows that $\rk_T(\sigma) =
\kappa\cdot \alpha +\rk_{T'}(\sigma)$ for any $\sigma  \in  T'$.  Hence,
$\rk(T') = \beta < \rk(T)$, so, by the induction hypothesis, there
exists a large sequence $\Zv' \preceq \Yv$ such that $T'$ avoids~$\Zv'$.
Now let $T'' = T\treerestrict \Zv'$; it is easy to see that $\rk(T'') <
\kappa\cdot \alpha +\ell(\Zv')$, so we can again apply the induction
hypothesis to get $\Zv \preceq \Zv'$ such that $T''$ avoids~$\Zv$.  It
follows that $T$ avoids~$\Zv$.

Case 3: $\rk(T)$ is of the form $\alpha+n$, where $\alpha $ is a limit
ordinal, $n < \omega$, and $\cf  \alpha  < \kappa_n$.

Let $\langle \alpha_\beta\suchthat  \beta  < \delta\rangle$ be
a strictly increasing sequence of ordinals which converges to $\alpha
$, where $\delta < \kappa_n$.  For each $\sigma  \in  P_{n+1}$, we
must have $\rk_T(\sigma) < \alpha $; hence, we can define a function
$F\funcfrom  P_{n+1}\to \delta$ by:  $F(\sigma)$ is the least $\beta $ such
that $\rk_T(\sigma) < \alpha_\beta$.  By Lemma~7.8, there is a large
sequence $\Yv' \preceq \Yv$ of length at least~$n+1$ such that,
for $\sigma \in
\prod_{i < n+1} \Yv'(i)$, $F(\sigma)$~depends only on~$\sigma \restrict
n$.  Let $T' = T \treerestrict \Yv'$.
Clearly $\rk_{T'}(\sigma) \le  \rk_T(\sigma)$
for all $\sigma \in \setfuncs {<\omega}\kappa$.  If $\sigma \in T'$ is
of length~$n$, and $\beta$ is the common value of $F(\sigma \concat \langle \gamma \rangle)$ for
$\gamma \in \Yv'(n)$, then $\rk_{T'}(\sigma \concat \langle \gamma \rangle)
< \alpha_\beta$ for all $\gamma  \in  \kappa_n$, so $\rk_{T'}(\sigma) \le
\alpha_\beta < \alpha $; this implies that $\rk_{T'}(\nullseq) \le  \alpha
+n-1$, so $\rk(T') < \rk(T)$.  Apply the induction hypothesis to~$T'$ to
get $\Zv \preceq \Yv'$ such that $T'$ avoids~$\Zv$; then $T$~also
avoids~$\Zv$.

     It is not hard to see that any value for~$\rk(T)$ less
than $\kappa\cdot \kappa_0$ falls under at least one of these three
cases, so the induction is complete.
\QED\enddemo

\proclaim{Corollary 7.11}
If $\kappa$~is an uncountable strong limit cardinal of cofinality~$\omega
$, and $\lambda < \kappa$, then any $\lambda$\snug-narrow covering
of\/~$\setfuncs \omega\kappa$ using open sets must have complexity at
least~$\kappa \cdot \kappa$.
\endproclaim

\demo{Proof}
Let $T$~be the tree associated with such a covering, and suppose that
$\rk(T) < \kappa\cdot\kappa$; then there is $\alpha < \kappa$ such that
$\rk(T) < \kappa\cdot\alpha$.  Choose the cardinals $\kappa_n$
and~$\lambda_n$ as specified in the global assumptions, so that
$\kappa_0$ is greater than $\alpha$ and~$\lambda$.  Then
Proposition~7.10 (with $\Yv = \nullseq$) states that there is a large
sequence~$\Zv$ such that $T$ avoids~$\Zv$, while Proposition~7.9 states
that there is no such~$\Zv$, so we have a contradiction.
\QED\enddemo

     Now we give the second version of the inductive argument.  In order
to reach higher tree ranks, we work with an entire collection of trees
simultaneously.  We will show not only that each tree in the collection
avoids some large sequence, but that one can find a relatively small
number of large sequences such that each tree in the collection avoids
at least one of them.

\proclaim{Proposition 7.12}
Let $\mu $ be an infinite cardinal less than~$\kappa_0$.  Suppose that
$\scriptX $~is a collection of well-founded trees $T \subseteq
\setfuncs {<\omega}\kappa$ such that $\sup\{\rk(T)\suchthat  T \in
\scriptX \} < \kappa^{\mu^+}$ (ordinal exponentiation) and\/ $|\scriptX
| \le \kappa_{-1}$.  Finally, suppose that $\Yv$~is a large sequence.
Then there is a collection~$\scriptC$ of large sequences $\Zv \preceq
\Yv$ such that $|\scriptC| \le \mu$ and, for every $T \in \scriptX$,
there exists $\Zv \in \scriptC$ such that $T$ avoids~$\Zv$.
\endproclaim

\demo{Proof}
Let $\alpha_0$ be the least ordinal which is greater than~$\rk(T)$ for
all $T \in  \scriptX $; then $\alpha_0 < \kappa^{\mu^+}$. If $\alpha_0 <
\omega$, the conclusion is trivial:  just let $\scriptC = \{\Zv\}$ where
$\Zv$~is any large sequence of length at least~$\alpha_0$ such that
$\Zv \preceq \Yv$. So suppose $\alpha_0
\ge \omega$.  Then there is a unique ordinal $\theta < \mu^+$ such that
$\omega \cdot \kappa^{\theta} \le  \alpha_0 < \omega \cdot
\kappa^{\theta +1}$. The proof will be by induction on $\theta$,
simultaneously for all sequences of cardinals $\kappa_n$ and~$\lambda_n$
satisfying the global assumptions.  (However, $\kappa$ and~$\mu$ will be
fixed.)

Suppose the statement is true for all $\theta' < \theta$.  For
convenience, we divide the induction step into two cases.

Case 1: $\omega
\cdot \kappa^{\theta} \le  \alpha_0 < \omega \cdot
\kappa^{\theta}\cdot \kappa_0$.  Let $\delta = \omega
\cdot \kappa^{\theta}$, and choose a strictly increasing
and continuous sequence $\langle \delta_\beta\suchthat  \beta  < \cf
\delta\rangle$ converging to $\delta$ such that $\delta_0 = 0$ and
$\delta_1 \le  \omega $.  Note that $\cf \delta$ is either $\cf \omega
$, $\cf \kappa$, or $\cf \theta $, so $\cf \delta \le  \mu $.

We will construct a sequence $\langle \scriptC_k \suchthat k < \omega
\rangle$ of sets of large sequences with the following properties:
\roster
\item"$\bullet$" $\scriptC_0 = \{\Yv\}$;
\item"$\bullet$" $|\scriptC_k| \le \mu$ for all~$k$;
\item"$\bullet$" if $\Zv' \in \scriptC_{k+1}$, then there is $\Zv \in
\scriptC_k$ such that $\Zv' \preceq \Zv$; and
\item"$\bullet$" if $\Zv \in \scriptC_k$ and $T \in \scriptX$, then
there is $\Zv' \preceq \Zv$ in~$\scriptC_{k+1}$ such that either $T$
avoids~$\Zv'$ or $$\rk(T \treerestrict \Zv') < \rk(T \treerestrict
\Zv).$$
\endroster
Once we have this sequence, we can let $\scriptC = \bigcup_{k < \omega}
\scriptC_k$.  Then $\scriptC$~will be a set of large sequences $\Zv
\preceq \Yv$, with $|\scriptC| \le \mu$.  For every $T \in \scriptX$,
there will be $\Zv \in \scriptC$ such that $T$ avoids~$\Zv$; if this
were not so, one could start with $\Zv_0 = \Yv$, find $\Zv_1 \in
\scriptC_1$ such that $\rk(T \treerestrict \Zv_1) < \rk(T \treerestrict
\Zv_0)$, then find $\Zv_2 \in
\scriptC_2$ such that $\rk(T \treerestrict \Zv_2) < \rk(T \treerestrict
\Zv_1)$, and so on, thus producing an infinite descending sequence of
ordinal ranks, which is impossible.  Therefore, $\scriptC$~will be as
desired.

Given~$\scriptC_k$, we will construct~$\scriptC_{k+1}$ by examining each
large sequence $\Zv \in \scriptC_k$ and thereby producing a collection
of at most~$\mu$ new large sequences to be put into~$\scriptC_{k+1}$.
So let $\Zv$~be an arbitrary member of~$\scriptC_k$, and proceed as
follows.

For each $T \in \scriptX$, we can express $\rk(T\treerestrict \Zv)$ in
the form $\delta\cdot \gamma_0+\gamma_1+n$ where $n < \omega $ and
$\gamma_1$~is zero or a limit ordinal less than~$\delta$, and this
expression is unique. Note that the number of possibilities
for~$\gamma_0$ is less than~$\kappa_0$, since $\alpha_0 <
\delta\cdot\kappa_0$. Let $f(T) = (n,\beta,c)$, where $\beta$~is the
unique~$\beta$ such that $\delta_\beta \le \gamma_1 < \delta_{\beta+1}$,
and $c$~is $0$~if $\gamma_0 = 0$, $1$ otherwise.  Note that the number
of possible values for $f(T)$ is at most $|\omega \times (\cf
\delta)\times 2| \le  \mu $.  We consider each possible triple
$(n,\beta,c)$ separately. Fix $(n,\beta,c)$ with $n < \omega$, $\beta <
\cf \delta$, and $c < 2$, and let $$\scriptX_{n\beta c} = \{T
\treerestrict \Zv \suchthat T \in
\scriptX\text{ and }f(T) = (n,\beta,c)\}.$$  We now consider several
subcases.

Subcase 1: $c = 0$.  Then the trees in~$\scriptX_{n\beta c}$ all have
rank less than~$\delta_{\beta+1}+n$, which is below~$\delta$, so we can
apply the induction hypothesis to get a collection of at most~$\mu$
large sequences $\Zv' \preceq \Zv$ such that every $T' \in
\scriptX_{n\beta c}$ avoids at least one of the sequences.  It follows
that every tree $T \in \scriptX$ such that $f(T) = (n,\beta,c)$
must avoid one of these sequences.  Add all of
these large sequences to~$\scriptC_{k+1}$.

Subcase 2: $c = 1$ and $\beta > 0$.  For each $T' \in \scriptX_{n\beta
c}$, express~$\rk(T')$ in the form $\delta\cdot \gamma_0+\gamma_1+n$
as above (where $\gamma_0$ and~$\gamma_1$ depend on~$T'$),
and let $T'_* = \{\sigma \in T'\suchthat
\rk_{T'}(\sigma) \ge  \delta\cdot \gamma_0\}$.  It is
easy to see that $\rk(T'_*) = \gamma_1+n < \delta_{\beta+1}+n < \delta$ for
each $T' \in \scriptX_{n\beta
c}$.  Therefore, we can apply the induction hypothesis to
the set $\{T'_*\suchthat  T' \in \scriptX_{n\beta
c}\}$ to get a collection of at most~$\mu$
large sequences $\Zv' \preceq \Zv$ such that every
such tree~$T'_*$ avoids at least one of the sequences.
If $T'_*$ avoids~$\Zv'$, then $\rk(T' \treerestrict \Zv')
< \delta\cdot \gamma_0+\ell(\Zv') <
\rk(T')$.  (Note that, if $T' = T \treerestrict \Zv$ and $\Zv' \preceq
\Zv$, then $T' \treerestrict \Zv' = T \treerestrict \Zv'$.)
Again, add all of
these large sequences~$\Zv'$ to~$\scriptC_{k+1}$.

Subcase 3: $c = 1$ and $\beta = 0$. Then, for each $T' \in
\scriptX_{n\beta c}$, we can express~$\rk(T')$ in the form $\delta\cdot
\gamma_0+n$, where $\gamma_0 < \kappa_0$ since $\alpha_0 < \delta\cdot
\kappa_0$. Since $\cf \delta \le  \mu  < \kappa_0$, we have $\cf
(\delta\cdot \gamma_0) < \kappa_0$.  Let $\nu$ be the predecessor
cardinal of $\kappa_0$, i.e., $2^{\lambda_{-1}}$.  Then we can
partition~$\delta\cdot \gamma_0$ into sets~$B^{T'}_{\xi}$, $\xi < \nu$,
none of which is cofinal in~$\delta\cdot \gamma_0$.  Now, for each
$\sigma \in P_{n+1}$, let $F(\sigma)$ be the function
from~$\scriptX_{n\beta c}$ to~$\nu$ defined by:  $F(\sigma)(T')$ is the
unique $\xi < \nu$ such that $\rk_{T'}(\sigma) \in B^{T'}_{\xi}$.  The
number of possible functions~$F(\sigma)$ is at
most~$\nu^{|\scriptX_{n\beta c}|}$; since $|\scriptX_{n\beta c}| \le
|\scriptX| \le \kappa_{-1} \le \lambda_{-1}$ and $\nu =
2^{\lambda_{-1}}$, we have $\nu^{|\scriptX_{n\beta c}|} \le \nu$, so
there are fewer than~$\kappa_0$ possible values of~$F(\sigma)$.  We can
now apply Lemma~7.8 to get a large sequence $\Zv' \preceq \Zv$ of length
at least~$n+1$ such that $F$~is constant on $\prod_{i < n+1} \Zv'(i)$.
Let $h\funcfrom \scriptX_{n\beta c}\to \nu$ be the constant value of~$F$
on this set.  For any $T' \in \scriptX_{n\beta c}$, if $\xi = h(T')$,
then, since $B^{T'}_{\xi}$ is a non-cofinal subset of the limit ordinal
$\delta\cdot \gamma_0$ and $\rk_{T'}(\sigma) \in B^{T'}_{\xi}$ for any
$\sigma \in \prod_{i < n+1} \Zv'(i)$, we have $\rk({T'}\treerestrict
\Zv') < \delta\cdot \gamma_0 \le \rk({T'})$.  Add this~$\Zv'$
to~$\scriptC_{k+1}$.

Once the relevant subcase step has been performed for each $\Zv \in
\scriptC_k$ and each $(n,\beta,c)$, the construction of~$\scriptC_{k+1}$
is complete.  We have ensured that $\scriptC_{k+1}$ has all of the
required properties.  This finishes Case~1 of the induction.

Case 2: $\omega \cdot \kappa^{\theta}\cdot \kappa_0
\le \alpha_0 < \omega \cdot \kappa^{\theta +1}$.  Fix~$n$
such that $\alpha_0 < \omega \cdot \kappa^{\theta}\cdot
\kappa_n$.  If we let $\kappa'_m = \kappa_{n+m}$ and $\lambda'_m =
\lambda_{n+m}$ for all~$m$ (including~$-1$), and define~$P'_m$ accordingly, then
the global assumptions will be satisfied for these new values.  Consider the
collection $$\scriptX' = \{T_\sigma\suchthat \sigma \in P_n,\, T \in
\scriptX \},$$
where $T_\sigma = \{\tau \suchthat  \sigma \concat \tau  \in
T\}$.
Clearly
$\scriptX'$~is a collection of trees $T'\subseteq \setfuncs
{<\omega}\kappa$, and $|\scriptX'|
\le  |\scriptX |\cdot |P_n| \le  \kappa_{n-1} = \kappa'_{-1}$.  Also, for
each $T \in \scriptX$ and each $\sigma \in P_n$, $$\rk(T_\sigma)
\le \rk(T) < \alpha_0 < \omega \cdot
\kappa^{\theta}\cdot \kappa'_0.$$  Therefore, we can apply Case~1
to get a collection~$\scriptC'$ of at most~$\mu$ large (for the
cardinals~$\kappa'_m$) sequences
$\Zv' \preceq \langle \Yv(n+i) \suchthat i < \ell(\Yv)-n \rangle$
such that, for each $T' \in \scriptX'$, there exists $\Zv' \in
\scriptC'$ such that $T'$ avoids~$\Zv'$.

For each $\sigma \in P_n$, define $F(\sigma)\funcfrom \scriptX \to
\scriptC'$
so that $F(\sigma)(T)$ is some $\Zv' \in \scriptC'$ such that $T_\sigma$
avoids~$\Zv'$.  Since $|\scriptC'| \le  \mu  < \kappa_0$
and $|\scriptX| \le \kappa_{-1}$, there are fewer than~$\kappa_0$
functions from~$\scriptX$ to~$\scriptC'$.  Therefore, by Lemma~7.8, there
is a large sequence~$\Yv_*$ of length at least~$n$ such that
$\Yv_* \preceq \Yv\restrict n$ and
$F$~is constant on $\prod_{i < n} \Yv_*(i)$.
Let $h = F(\sigma)$ for $\sigma $ in this set, and let
$\scriptC = \{(\Yv_*\restrict n) \concat \Zv' \suchthat \Zv' \in
\scriptC'\}$.  Then $\scriptC$ is a collection of at most~$\mu$ large
sequences $\Zv \preceq \Yv$, and each $T \in \scriptX$ avoids some $\Zv
\in \scriptC$, namely $(\Yv_*\restrict n) \concat h(T)$.  This completes
the induction.
\QED\enddemo

Just as for Corollary~7.11, we can apply Proposition~7.12 and
Proposition~7.9 to a single given tree~$T$ (i.e., let $\scriptX =
\{T\}$) with the cardinals $\mu$ and~$\kappa_n$ chosen as large as
necessary below~$\kappa$ to get:

\proclaim{Corollary 7.13}
If $\kappa$~is an uncountable strong limit cardinal of cofinality~$\omega
$, and $\lambda < \kappa$, then any $\lambda$\snug-narrow covering
of\/~$\setfuncs \omega\kappa$ using open sets must have complexity at
least~$\kappa^\kappa$.
\QED\endproclaim

There is no reason to believe that the lower bound~$\kappa^\kappa$
obtained here is optimal; improvements in the argument might yield
better results.  The obvious way to provide an upper limit on this
ordinal would be to produce an explicit open narrow covering, or an
algebra with no infinite free subset, and compute the rank of the
corresponding tree.  For instance, in the constructible universe, with
$\kappa = \aleph_\omega$, one can consider the algebra consisting of all
operations (unary, binary, etc.) on~$\kappa$ which are definable
in~$L_{\kappa^+}$; there are countably many of these.  Devlin and
Paris~\cite{\DevlinParis} have shown that this algebra has no infinite
free subset.  However, their proof gives no information about the rank
of the tree of finite free sequences.  I do not know of any upper bound
for this rank beyond the obvious fact that it is less
than~$\kappa^+$.


There are other families of sets besides the open sets for which one
could make a similar study of complexity of narrow coverings.  For
instance, one could consider the case of $\SIGMA^0_2$~sets.  A reduction
argument similar to that of Proposition~7.1 shows that, if there is a
narrow covering using $\SIGMA^0_2$~sets, then there is a narrow
partition using~$\SIGMA^0_2$, and hence~$\DELTA^0_2$, sets.  One can
assign an ordinal rank to such a partition in various ways, such as the
first level in the difference hierarchy which includes all of the
individual $\DELTA^0_2$~sets, and then ask what the smallest possible
rank for the partition is.  However, it does not seem useful to study
this question yet, since no case is currently known where a
$\DELTA^0_2$~narrow partition exists and a clopen narrow partition does
not.

\head 8.  Open Problems \endhead

There are many open questions related to the concepts studied in
this paper; here are some of the more interesting ones.

1.  Does $\NNC(\kappa,\aleph_1,F_{\sigma})$ (or even
$\NNC(\kappa,{<}\kappa,\Borel)$) follow from
the real-valued measurability of $\kappa$?

2.  What is the exact consistency strength of
$\NNC(\kappa,\aleph_1,\Borel)$?  In particular, does it imply the
existence of\/~$0^\sharp$?

3.  Must the least $\kappa$ satisfying $\NNC(\kappa,\lambda,\Borel)$
for a given $\lambda$ actually satisfy $\NNC(\kappa,\allowbreak{<}\kappa,
\Borel)$?

4.  Does $\NNC(\kappa,\lambda,\open)$ always imply
$\NNC(\kappa,\lambda,\Borel)$?

5. Is $\NNC(\kappa,\lambda,\Borel)$ preserved by any forcing
with the countable chain condition?

6.  Can a cardinal $\kappa \le  2^{\aleph_0}$ carry a uniform
$(\kappa,\aleph_1;\omega)$\snug-nonregular ultrafilter?

7.  Does Projective Determinacy imply that all projective
subsets of\/ $\setfuncs \omega\omega$ are $U$-measur\-able, where
$U$ is a nonprincipal ultrafilter over~$\omega$?
(Louveau~\cite{\Louveau} mentions that, if a measurable cardinal exists
or if\/ $\text{MA}{+}\neg\text{CH}$ holds, then there are many
ultrafilters $U$ such that all $\SIGMA^1_2$ sets are
$U$-measurable, but that it is open whether this is so
for {\sl all}~$U$ in these cases.)  Does Projective Determinacy imply
$\NNC(\aleph_0,{<}\aleph_0,\text{projective})$?

8. What is the least possible rank for the tree of all finite free
sequences obtained from an algebra of size~$\aleph_\omega$ with fewer
than~$\aleph_\omega$ operations and no infinite free subset?

9.  Mr\'owka~\cite{\MrowkaSID} gives some ostensibly weaker variants of
$\NNC(2^{\aleph_0},\aleph_1,F_\sigma)$ which would still suffice for his
metric space constructions.  One such variant is: $\setfuncs
\omega{(\setfuncs \omega2)}$ cannot be written as a union of sets $A_n$
($n < \omega$) where $A_n$ is $\aleph_1$\snug-narrow in the $n$\snug'th
coordinate and $A_n$ is $F_\sigma$ in the product topology on $\setfuncs
\omega{(\setfuncs \omega2)}$ where the $n$\snug'th factor $\setfuncs
\omega2$ is given the usual Cantor topology while the other factors are
given the discrete topology.  Are these variants actually weaker? Can
they be attained using weaker large cardinals (or none at all)?
Can such a metric space be constructed at all without large cardinals?

10.  Mr\'owka~\cite{\MrowkaSID} also mentions the statement
$\NNC(2^\kappa,\kappa^+,F_\sigma)$ where $\kappa$ is an uncountable
strong limit cardinal of cofinality~$\omega$.  Is this consistent
with~ZFC?  If so, it will require stronger large cardinals than
the ones used in this paper: the statement clearly implies
$2^\kappa > \kappa^+$, so $\kappa$ must be a counterexample
to the Singular Cardinals Hypothesis, and this entails the
consistency of measurable cardinals of high
order~\cite{\Gitik}.

11.  What happens if one considers products of sets of different
sizes?  That is, when can one express the infinite product
$\Pi_{n < \omega} \kappa_n$ as a union of `nice' sets
$A_n$ ($n < \omega$) such that $A_n$ is $\lambda_n$\snug-narrow
in the $n$\snug'th coordinate?
This is of interest even in the finite-dimensional version
with no restrictions on the sets~$A_n$; Simms~\cite{\Simms}
lists this as Open Problem~2, and cites results of
Ristow~\cite{\Ristow} that settle it assuming a weak form of~GCH
(every limit cardinal is a strong limit).

In fact, one can ask the same question about finite products of finite
sets, but this question has been settled.  If\/ $|X_j| = k_j$ for each
$j < n$, then $\prod_{j < n} X_j$ can be expressed as the union of sets
$A_j$ ($j < n$), where $A_j$ is $l_j$\snug-narrow in the $j$\snug'th
coordinate, if and only if\/ $\sum_{j < n} (l_j-1)/k_j \ge 1$.  The
necessity of this inequality is a simple counting argument.  Conversely,
if the inequality holds, then one can canstruct suitable sets~$A_j$ by
the following modification of the proof of Proposition~2.4, due to
J.~Rickard (personal communication): partition the half-open
interval~$[0,1)$ into intervals~$[y_j,z_j)$ for $j < n$ so that $z_j -
y_j \le (l_j-1)/k_j$, let $X_j = \{0,1,\dots,k_j-1\}$, and define sets
$A_j \subseteq \prod_{j < n} X_j$ for $j < n$ as follows: $$x \in A_j
\iff y_j \le \biggl(\sum_{i=0}^{n-1} x(i)/k_i\biggr) \bmod 1 < z_j.$$


\Refs

\ref \no \Bagemihl \by F. Bagemihl \paper A decomposition of an
infinite dimensional space \jour Z. Math. Logik Grundlag. Math.
\vol 31 \yr 1985 \pages 479--480 \endref

\ref \no \Barnes \by R.
Barnes, Jr. \paper The classification of the closed-open and the
recursive sets of number-theoretic functions \paperinfo Doctoral
Dissertation, University of California, Berkeley \yr 1966 \endref

\ref \no \BaumgartnerHajnal \by J. Baumgartner and A. Hajnal \paper A
proof (involving Martin's Axiom) of a partition relation \jour Fund. Math.
\vol 78 \yr 1973 \pages 193--203 \endref

\ref \no \Devlin \by K. Devlin \paper Some weak versions of large cardinal
axioms \jour Ann. Math. Logic \vol 5 \yr 1973 \pages 291--325 \endref

\ref \no \DevlinParis \by K. Devlin and J. Paris \paper More on the free
subset problem \jour Ann. Math. Logic \vol 5 \yr 1973 \pages 327--336
\endref

\ref \no \Dougherty \by R. Dougherty \paper Narrow coverings
of $\omega$\snug-product spaces \paperinfo Doctoral
Dissertation, University of California, Berkeley \yr 1985 \endref

\ref \no \Gitik \by M. Gitik \paper The strength of the failure of
the Singular Cardinal Hypothesis \jour Ann. Pure Appl. Logic
\vol 51 \yr 1991 \pages 215--240 \endref

\ref \no \Jech \by T. Jech \book Set theory \publ Academic Press \publaddr
New York \yr 1978 \endref

\ref \no \Koepke \by P. Koepke \paper The consistency strength of the
free-subset property for ${}^{\omega}\omega $ \jour J. Symbolic
Logic \vol 49 \yr 1984 \pages 1199--1204 \endref

\ref \no \KuratowskiP
\by C. Kuratowski
\paper Sur une caract\'erisation des alephs \jour Fund. Math. \vol 38
\pages 14--17 \endref


\ref \no \Louveau \by A. Louveau \paper Une m\'ethode topologique pour l'etude
de la propri\'et\'e de Ramsey \jour Israel J. Math. \vol 23 \yr 1976
\pages 97--116 \endref

\ref \no \Maharam \by D. Maharam \paper On homogeneous measure algebras
\jour Proc. Nat. Acad. Sci. U.S.A. \vol 28 \yr 1942 \pages 108--111 \endref

\ref \no \Moschovakis \by Y. Moschovakis \book Descriptive set theory
\publ North-Holland \publaddr Amsterdam \yr 1980 \endref

\ref \no \MrowkaNMC \by S. Mr\'owka \paper $N$\snug-compactness,
metrizability, and covering dimension \inbook Rings of continuous
functions \ed C.~Aull \publ Marcel Dekker \publaddr New York and
Basel \yr 1985 \pages 247--275 and 312--314 \endref

\ref \no \MrowkaSID \bysame \paper Small inductive dimension
of completions of metric spaces \paperinfo preprint \endref

\ref \no \Prikry \by K. Prikry \paper Changing measurable into accessible
cardinals \jour Dissertationes Math. (Rozprawy Mat.) \vol 68
\yr 1970 \pages 5--52 \endref

\ref \no \Ristow \by A. Ristow \paper The existence of certain partitions
on Cartesian products \jour Z. Math. Logik Grundlag. Math.
\vol 24 \yr 1978 \pages 325--333 \endref

\ref \no \Shelah \by S. Shelah \paper Independence of strong partition
relation for small cardinals, and the free-subset problem \jour J.
Symbolic Logic \vol 45 \yr 1980 \pages 505--509 \endref
      

\ref \no \Silver \by J. Silver \paper Every analytic set is Ramsey
\jour J. Symbolic Logic \vol 35 \yr 1970 \pages 60--64 \endref

\ref \no \Simms \by J. Simms \paper Sierpi\'nski's theorem
\jour Simon Stevin \vol 65 \yr 1991 \pages 69--163 \endref

\ref \no \Solovay \by R. Solovay \paper Real-valued measurable cardinals
\inbook Axiomatic set theory \vol I \ed D. Scott \publ Amer. Math. Soc.
\publaddr Providence, Rhode Island \yr 1971 \pages 397--428 \bookinfo
Proc. Sympos. Pure Math. 13 \endref

\endRefs
\enddocument